\documentclass[preprint,3p,times]{elsarticle}

\usepackage{array}
\usepackage{color}
\usepackage{tabularx}
\usepackage{graphicx}
\usepackage{amsmath}
\usepackage{amssymb}
\usepackage{amsfonts}	
\usepackage{moreverb}
\usepackage{dsfont}
\usepackage{bm}
\usepackage{multirow}
\usepackage{soul}
\usepackage{txfonts} % REAL SET  \varmathbb{R}
\usepackage{tikz}
\newcommand\BibTeX{{\rmfamily B\kern-.05em \textsc{i\kern-.025em b}\kern-.08em
T\kern-.1667em\lower.7ex\hbox{E}\kern-.125emX}}

\newcommand{\reds}[1]{{#1}}  %\color{red} #1}}  % adding red color for the final revision: put 'red' for the reviewers; put 'black' for the final draft.

\newcommand{\x}{\mbf{x}}

\newcommand{\mbf}[1]{\mathbf{#1}}			%
\newcommand{\p}{\textnormal{P}}

\newcommand{\halb}{\frac{1}{2}}

\newcommand{\be}{\begin{equation}}
\newcommand{\ee}{\end{equation}}
\newcommand{\bdm}{\begin{displaymath}}
\newcommand{\edm}{\end{displaymath}}

\newfont{\numerikEleven}{ecrm1000}
\newfont{\numerikTen}{cmss10}
\newfont{\numerikNine}{cmss9}
\newfont{\numerikEight}{cmss8}
\newcommand{\Mass}[1]{\mathbf{M}^{\textsl{#1}} }
\newcommand{\iMass}[1]{\mathbf{M}^{\mathbf{-1}\,\textsl{#1}} }
\newcommand{\barr}[2]{\mathbf{#1}_{\text{\textbf{\textsl{p}}}}^{\textsl{#2}}  }
\newcommand{\tild}[2]{\mathbf{#1}_{\text{\textbf{\textsl{v}}}}^{\textsl{#2}}  }
\newcommand{\operator}[2]{\mathbf{#1}^{\textsl{#2}}  }
\newcommand{\Hoperator}[2]{\mathbb{#1}^{\textsl{#2}}  }
\newcommand{\dof}[1]{\widehat{\mathbf{#1}}}

\DeclareMathAlphabet{\mathpzc}{OT1}{pzc}{m}{it}

\journal{Applied Numerical Mathematics}

\begin{document} 
%!=========================================================================
%!      F R O N T    M A T T E R 
\begin{frontmatter}
%-------------------------------------------------------
% TITLE
\title{Spectral semi-implicit and space-time discontinuous Galerkin methods for the incompressible Navier-Stokes equations on staggered Cartesian grids} 
%-------------------------------------------------------
% AUTHORS
\author[UniTN]{Francesco Fambri}
\ead{francesco.fambri@unitn.it}

\author[UniTN]{Michael Dumbser$^{*}$}
\ead{michael.dumbser@unitn.it}
\cortext[cor1]{Corresponding author}

%-------------------------------------------------------
% INSTITUTIONS
\address[UniTN]{Laboratory of Applied Mathematics,\\Department of Civil, Environmental and Mechanical Engineering, \\University of Trento, Via Mesiano, 77 - 38123 Trento, Italy.}  
%-------------------------------------------------------
% ABSTRACT
\begin{abstract}
In this paper two new families of arbitrary high order accurate spectral discontinuous Galerkin (DG) finite element methods are derived on \textit{staggered} Cartesian grids for the solution of the incompressible Navier-Stokes (NS) equations in two and three space dimensions. 
% high-order in space (theta-method)
The discrete solutions of pressure and velocity are expressed in the form of piecewise polynomials along \emph{different} meshes. While the pressure is defined on the control volumes of the main grid, the velocity components are defined on edge-based dual control volumes, leading to a \emph{spatially staggered} mesh. 
Thanks to the use of a nodal basis on a tensor-product domain, all discrete operators can be written efficiently as a combination of simple one-dimensional operators in a dimension-by-dimension fashion. 

In the first family, high order of accuracy is achieved only in space, while a simple semi-implicit time discretization is derived by introducing an implicitness factor $\theta\in [0.5,1]$ for the pressure gradient in the momentum equation. 
The real advantages of the \emph{staggering} arise after substituting the discrete momentum equation into the weak form of the continuity equation. In fact, the resulting linear system for the pressure is \emph{symmetric} and \emph{positive definite} and either block penta-diagonal (in 2D) or block hepta-diagonal (in 3D). 
As a consequence, the pressure system can be solved very efficiently by means of a classical matrix-free conjugate gradient method. From our numerical experiments we find that the pressure system appears to be reasonably 
\textit{well-conditioned}, since in all test cases shown in this paper the use of a preconditioner was \textit{not} necessary. 
This is a rather unique feature among existing implicit DG schemes for the Navier-Stokes equations. In order to avoid a stability restriction due to the viscous terms, the latter are discretized implicitly using again
a staggered mesh approach, where the viscous stress tensor is also defined on the dual mesh. 
%high-order in time

The second family of staggered DG schemes proposed in this paper achieves high order of accuracy also in time by expressing 
the numerical solution in terms of piecewise \emph{space-time} polynomials. In order to circumvent the low order of accuracy 
of the adopted fractional stepping, a simple iterative Picard procedure is introduced, which leads to a space-time 
pressure-correction algorithm. In this manner, the symmetry and positive definiteness of the pressure system are not compromised. 
The resulting algorithm is stable, computationally very efficient, and at the same time arbitrary high order accurate 
in both space and time. These features are typically not easy to obtain all at the same time for a numerical method applied 
to the incompressible Navier-Stokes equations. 
The new numerical method has been thoroughly validated for approximation polynomials of degree up to $N=11$, using a large set 
of non-trivial test problems in two and three space dimensions, for which either analytical, numerical or experimental reference 
solutions exist. 
\end{abstract}
%-------------------------------------------------------
% KEY WORDS
%-------------------------------------------------------
\end{frontmatter}
%!======================================================
% 		A L L     I N C L U D E S 
%-------------------------------------------------------

\section{Introduction} \label{sec:introduction}

In this paper two novel families of efficient arbitrary high order accurate discontinuous Galerkin (DG) methods are presented for the solution of the two- and 
three-dimensional  incompressible Navier-Stokes (NS) equations on \textit{staggered} Cartesian meshes. The governing partial differential equations (PDE) read  
\begin{align}
&\frac{\partial \mathbf{v}}{\partial t} + \nabla \cdot \mathbf{F} + \nabla p = 0 \label{eq:NSmom}\\
&\nabla \cdot \mathbf{v} = 0 \label{eq:NSinc}
\end{align}
where $\mathbf{v}=\mathbf{v}(\mathbf{x},t)=(u(\mathbf{x},t),v(\mathbf{x},t),w(\mathbf{x},t))$ is the velocity vector in three space dimensions, 
$p=p(\mathbf{x},t)$ is the normalized fluid pressure, $\mathbf{x}=(x,y,z)$ is the vector of the spatial-coordinates and $\mathbf{F}=(\mathbf{F}_u, \mathbf{F}_v, \mathbf{F}_w)$ 
is the flux tensor that contains both, nonlinear convection $\mathbf{F}_c = \mathbf{v} \otimes \mathbf{v}$ and diffusion $\mathbf{F}_d = - \nu \nabla \mathbf{v}$, and which 
therefore reads 
\begin{align}
\mathbf{F} = \mathbf{F}_c + \mathbf{F}_d = \mathbf{v} \otimes \mathbf{v} - \nu \nabla \mathbf{v}, % = \left(\begin{array}{c|c|c} \mathbf{F}_u & \mathbf{F}_v & \mathbf{F}_w\end{array} \right),
\label{eqn:fluxtensor} 
\end{align}
where $\nu$ is the kinematic viscosity. 

The incompressible Navier-Stokes equations \eqref{eq:NSmom} and \eqref{eq:NSinc} are of great interest for practical applications concerning the simulation of fluid flow in hydraulics, mechanical and naval engineering, 
oceanography and geophysics, physiological fluid flow in the human cardiovascular and human respiratory system, just to mention a few, but also in astrophysics or high-energy 
physics when the compressibility of high-density plasma becomes negligible. 
Because of this great interest across many different scientific disciplines, many attempts in resolving these equations have been made in the past, but research on numerical schemes
for the Navier-Stokes equations remains an important research topic even nowadays. For many decades either finite-difference schemes \cite{HarlowWelch,Patankar1972,patankar,vanKan1986} 
or continuous finite element methods \cite{Taylor1973,Brooks1982,Hughes1986,Fortin1981,Verfuerth,Heywood1982,Heywood1988,Arnold1984,Brezzi1989} were the state of the art. Only more recently, the discontinuous  
Galerkin (DG) finite-element method is used for the solution of the incompressible Navier-Stokes equations. 

Reed and Hill were the first in introducing the DG finite-element discretization \cite{Reed:1973} for the solution of neutron-transport equations. Later, Cockburn and Shu 
extended the DG framework to the general case of non-linear systems of hyperbolic conservation laws in a series of well-known fundamental papers 
\cite{Cockburn1989a,Cockburn1989b,Cockburn1990,CockburnShu98}. Further to that, the nonlinear $L_2$ stability of DG methods has been proven by Jiang and Shu \cite{Jiang1994} 
by demonstrating the validity of a cell entropy inequality for semi-discrete DG schemes, and then the proof has been extended to the case of systems in \cite{BarthCharrier,HouLiu}. 
%The time-discretization of  DG finite element methods is still mostly resolved through 
Initially, DG schemes were only used as higher-order spatial discretization, while time discretization was done with standard TVD Runge-Kutta schemes, leading to the 
family of classical Runge-Kutta-DG (RKDG) schemes. For alternative Lax-Wendroff-type or ADER-type time discretizations in the DG context, see \cite{QiuDumbserShu,dumbser_jsc,taube_jsc}. 
A review of DG finite element methods is provided in \cite{cockburn_2000_dg,cockburn_2001_rkd}. 
Even if higher-order DG schemes became more and more attractive and popular in recent years, probably the major drawback of explicit DG methods consists in the severe CFL stability 
condition that make the time step proportional to $1/(2N+1)$, where $N$ is the degree of the approximation polynomials used in the DG scheme. 
The DG method has been also extended to a uniform space-time formalism by Van der Vegt et al. \cite{spacetimedg1,spacetimedg2,KlaijVanDerVegt}, resulting in a fully implicit 
discretization. On the counterpart, a fully implicit DG formulation leads to a globally coupled nonlinear system for the degrees of freedom of the space-time DG polynomials, 
the solution of which can become computationally very demanding at every single time-step. 
The first DG method for the compressible Navier-Stokes equations has been presented by Bassi and Rebay in \cite{BassiRebay} and and Baumann and Oden \cite{BaumannOden1,BaumannOden2}. 
Notice that the DG  finite-element formulation of the parabolic (second order) terms in the equations, or for even higher order spatial derivatives, is not straightforward 
\cite{CockburnShu1998,yan2002,levy2004}. A unified analysis of DG schemes for elliptic problems is outlined in \cite{Arnold2001}. 
Many other DG methods have been presented for the Navier-Stokes equations in the meantime, see for example \cite{MunzDiffusionFlux,stedg2,DumbserNSE,HartmannHouston1,HartmannHouston2,Crivellini2013,KleinKummerOberlack2013} for a non-exhaustive overview of the ongoing research in this very active field. 

Moreover, the elliptic character of the incompressible Navier-Stokes equations introduces an important difficulty in in their numerical solution: whenever the smallest physical 
or numerical perturbation arises in the fluid flow then it will instantaneously affect the entire computational domain. Thus, in principle, the most natural way would be a 
fully implicit discretization of the governing equations. 
%From one hand the higher order resolution and suppleness of DG methods became even more attractive, on the other the CFL restriction of explicit methods remained an unresolved problem in spite of %computational efficiency.
The elliptic behaviour of the pressure can be avoided by weakening the incompressibility condition, i.e. by introducing the so-called method of artificial compressibility, 
see \cite{chorin1967,chorin1968}, which was also used in the DG finite element framework by Bassi et al. in \cite{Bassi2007}.

%Obviously, this procedure will produce a hardly-coupled non-linear system in the velocity and pressure field whose resolution is highly computationally demanding. \cite{}

%Further to that a general unification of finite-volume schemes and DG finite-element methods have been introduced by Dumbser et al. \cite{Dumbser2008} giving a general framework for the construction of 

%This fact led the numerical community to formulate new effective numerical methods by weakening the incompressibility condition \cite{...}

It has to be noticed that a family of very efficient semi-implicit finite difference methods for staggered structured and unstructured grids has been developed by Casulli et al. 
in the field of hydrostatic and non-hydrostatic gravity-driven free-surface and sub-surface flows, see 
\cite{CasulliCompressible,Casulli1990,CasulliCheng,Casulli1999,CasulliWalters,CasulliZanolli2002,CasulliStelling2011,Casulli2009,CasulliZanolli2010,CasulliVOF}. These methods have been theoretically 
analyzed, for example, in \cite{CasulliCattani,BrugnanoCasulli,BrugnanoCasulli2,BrugnanoSestini,CasulliZanolli2012}. 
In the above-mentioned semi-implicit framework, the schemes ensure exact mass conservation thanks to a conservative finite-volume formulation of the continuity equation 
and a rigorous nonlinear treatment of its implicit discretization. Moreover, numerical stability is ensured for large Courant numbers \reds{ (for free-surface  hydrodynamics or for compressible gas dynamics) and is independent of the kinematic viscosity}.
 The main advantage of making use 
of a semi-implicit discretization is that the numerical stability can be obtained for large time-steps without leading to an excessive computational demand. 

% by means of a weakly non-linear system for the pressure. A rigorous study of these special non-linear 
%
Thanks to their computational efficiency, these semi-implicit methods have been later also extended to the simulation of hydrostatic and non-hydrostatic blood flow in the human 
arterial system in two and three space dimensions \cite{CasulliDumbserToro,Blood3D2014}, but also to the simulation of the flow of compressible fluids in compliant tubes 
\cite{DumbserIbenIoriatti}. A generalization to the compressible Navier-Stokes equations with general equation of state has been introduced in \cite{DumbserCasulli2016}. 

Very recently, the aforementioned family of efficient semi-implicit finite-difference methods has been extended to a higher-order DG formulation for the shallow water 
equations, originally on staggered Cartesian grids \cite{DumbserCasulli2013} and then also on general unstructured meshes \cite{TavelliDumbser2014}. Based on the same ideas, 
a high order staggered DG scheme for the two-dimensional incompressible Navier-Stokes equations has been presented in \cite{TavelliDumbser2014b} and \cite{TavelliDumbser2015}, 
while the extension to three-dimensional unstructured meshes was achieved in \cite{TavelliDumbser2016}. 
Several alternative attempts of combining the stability properties of semi-implicit methods with the higher-order of accuracy of DG methods have been made in 
\cite{Dolejsi2004,Dolejsi2007,Dolejsi2008} for compressible flows and for nonlinear convection diffusion equations, and more recently in \cite{GiraldoRestelli,Tumolo2013} 
for the shallow water equations. In all these methods, a \textit{collocated grid} was used. A novel family of DG schemes on \textit{edge-based staggered grids} has been presented 
by Chung et al. in \cite{chung2012staggered,StaggeredDG2,StaggeredDG3,ChungNS}, while an interesting analysis of DG methods on \textit{vertex-based staggered grids} has been outlined 
in \cite{Liu2007,Liu2008}. For a review of \textit{spectral} DG FEM schemes on \textit{collocated grids}, the reader is referred to the work of Kopriva and Gassner et al. 
\cite{Kopriva2006,Kopriva2010,Becketal,Gassner2011,Gassner2013,Gassner2016}, and references therein, while classical spectral element methods for the Navier-Stokes equations can 
be found in the work of Canuto et al. \cite{Canuto1985,Canuto1984,Canuto1996,Canuto1998}. 

In this paper, two new families of \textit{spectral} semi-implicit and spectral space-time DG methods for the solution of the two and three dimensional Navier-Stokes equations on 
\textit{edge-based staggered} Cartesian grids are presented and discussed, following the ideas outlined in \cite{DumbserCasulli2013} for the shallow water equations. 
In the resulting schemes, all discrete operators can 
be written as a combination of simple one-dimensional operators, applied in a dimension-by-dimension fashion, thanks to the use of tensor-product control volumes. 
In this paper, we show numerical results using approximation polynomials of degree up to $N=11$ in both space and time. To the knowledge of the authors, such a high order of accuracy 
in space and time has \textit{never} been reached before with any DG scheme applied to the incompressible Navier-Stokes equations. 

The rest of the paper is organized as follows:  Section \ref{sec:LO} is dedicated to staggered semi-implicit DG schemes that achieve high order of accuracy only in space, 
while Section \ref{sec:HO} is devoted to high order staggered space-time DG schemes, which achieve arbitrary high order of accuracy in both space and time. The paper is rounded-off
by some concluding remarks in Section \ref{sec:conclusion}. 

 %The pressure gradients in the momentum equations and the incompressibility condition are discretized implicitly in order to obtain a very well conditioned algebraic linear system for the pressure. By doing this, the incompressibility condition is ensured up to machine precision in the weak sense of the governing equations. From one hand the nonlinear advective terms are discretized through a 

\section{Spectral semi-implicit DG schemes on staggered Cartesian grids}
 \label{sec:LO}

\subsection{Numerical method}
The staggered DG approach \cite{DumbserCasulli2013,TavelliDumbser2014,TavelliDumbser2014b,TavelliDumbser2015} is based on a weak formulation of the governing partial differential equations integrated along \emph{different sets} of \emph{overlapping} control volumes $d\Omega$, $d\Omega^*_x$, $\partial\Omega^*_y$, $d\Omega^*_z$ that define the \emph{main grid} and the three different \emph{edge-based staggered} (dual) grids respectively, 
\begin{align}
%&\int \limits_{\Omega_x} \mathit{f_u} \frac{\partial u}{\partial t} + \nabla \cdot \mathbf{F_u} + \partial_x p  = 0 \label{eq:NSmom}\\
%&\int \limits_{\Omega_y} \mathit{f_v} \frac{\partial v}{\partial t} + \nabla \cdot \mathbf{F_v} + \partial_y p = 0 \label{eq:NSmom}\\
%&\int \limits_{\Omega_z} \mathit{f_w} \frac{\partial w}{\partial t} + \nabla \cdot \mathbf{F_w} + \partial_z p = 0 \label{eq:NSmom}\\
\int \limits_{d\Omega^*_x}  \mathit{f_u} \left( \frac{\partial u}{\partial t} + \nabla \cdot \mathbf{F}_u + \partial_x p \right) d\x = 0, \hspace{0.2cm}
&\int \limits_{d\Omega^*_y} \mathit{f_v} \left( \frac{\partial v}{\partial t} + \nabla \cdot \mathbf{F}_v + \partial_y p \right) d\x = 0,& \hspace{0.2cm}
\int \limits_{d\Omega^*_z} \mathit{f_w} \left( \frac{\partial w}{\partial t} + \nabla \cdot \mathbf{F}_w + \partial_z p \right)  d\x = 0, \label{eq:wNSmom}\\
&\int \limits_{d\Omega} \mathit{f} \, \left( \nabla \cdot \mathbf{v}\right) d\x = 0,& \label{eq:wNSinc}
\end{align}
where $\mathit{f_u}$, $\mathit{f_v}$, $\mathit{f_w}$ and $\mathit{f}$ are the so called test-functions and $d\x = dx dy dz$.
If $\Omega$ is the computational domain, then the following properties hold for the staggered grids: 
\begin{align*}
&\Omega = \bigcup_i d\Omega_i =  \bigcup_i d\Omega^*_{x,i}=  \bigcup_i d\Omega^*_{y,i}=  \bigcup_i d\Omega^*_{z,i},\\
& \varnothing = \bigcup_{i\neq j} \left(  d\Omega^{\circ}_i \cap d\Omega^{\circ}_j \right) =  \bigcup_{i\neq j}  \left(  d\Omega^{*\circ}_{x,i}\cap  d\Omega^{*\circ}_{x,j} \right)=  \bigcup_{i\neq j}   \left(  d\Omega^{*\circ}_{y,i} \cap d\Omega^{*\circ}_{y,j} \right)=  \bigcup_{i\neq j}   \left(  d\Omega^{*\circ}_{z,i} \cap  d\Omega^{*\circ}_{z,j} \right).
\end{align*}
where $\circ$ denotes the interior of the cell without the boundary and the indices $i$ and $j$ run over all spatial elements of the corresponding main or dual mesh, respectively. 
Then, the discretization of the PDE restricts the solution of the physical variables $(u_h,v_h,w_h,p_h)$ to belong to the spaces 
%$\mathpzc{U}_N$, $\mathpzc{V}_N$, $\mathpzc{W}_N$, $\mathpzc{P}_N$  
of tensor products of piecewise polynomials of maximum degree $N$ with respect to the corresponding 
main or dual mesh. By dividing the domain in $N_x$, $N_y$ and $N_z$ elements in $x$, $y$ and $z$ direction, the control volumes for the pressure on the main grid are  
given by 
\begin{equation*}
d\Omega_{i,j,k}=[x_{i-\frac{1}{2}},x_{i+\frac{1}{2}}] \times [y_{j-\frac{1}{2}},y_{j+\frac{1}{2}}] \times [z_{k-\frac{1}{2}},z_{k+\frac{1}{2}}] \equiv T_{i,j,k},
\end{equation*} 
while the corresponding \emph{edge-based staggered} dual control volumes are 
\begin{equation*}
d\Omega^*_{i+\frac{1}{2},j,k}=[x_{i},x_{i+1}] \times [y_{j-\frac{1}{2}},y_{j+\frac{1}{2}}] \times [z_{k-\frac{1}{2}},z_{k+\frac{1}{2}}] \equiv T_{i+\frac{1}{2},j,k},
\end{equation*} 
for the velocity component $u$, 
\begin{equation*}
d\Omega^*_{i,j+\frac{1}{2},k}=[x_{i-\frac{1}{2}},x_{i+\frac{1}{2}}] \times [y_{j},y_{j+1}] \times [z_{k-\frac{1}{2}},z_{k+\frac{1}{2}}] \equiv T_{i,j+\frac{1}{2},k},
\end{equation*} 
for the velocity component $v$, 
\begin{equation*}
d\Omega^*_{i,j,k+\frac{1}{2}}=[x_{i-\frac{1}{2}},x_{i+\frac{1}{2}}] \times [y_{j-\frac{1}{2}},y_{j+\frac{1}{2}}] \times [z_{k},z_{k+1}] \equiv T_{i,j,k+\frac{1}{2}}, 
\end{equation*} 
for the velocity component $w$, respectively. The chosen staggered mesh corresponds to the one used in \cite{DumbserCasulli2013} for the shallow water equations, 
where each velocity component is defined on a different staggered mesh. In alternative, the entire velocity vector can also be defined on a single edge-based
dual grid, according to the choice made in \cite{Bermudez1998,Bermudez2014,USFORCE,TavelliDumbser2014,TavelliDumbser2014b}. 

\begin{figure} 
\centering 
\begin{tabular}{cc}
			\includegraphics[width=0.45\textwidth]{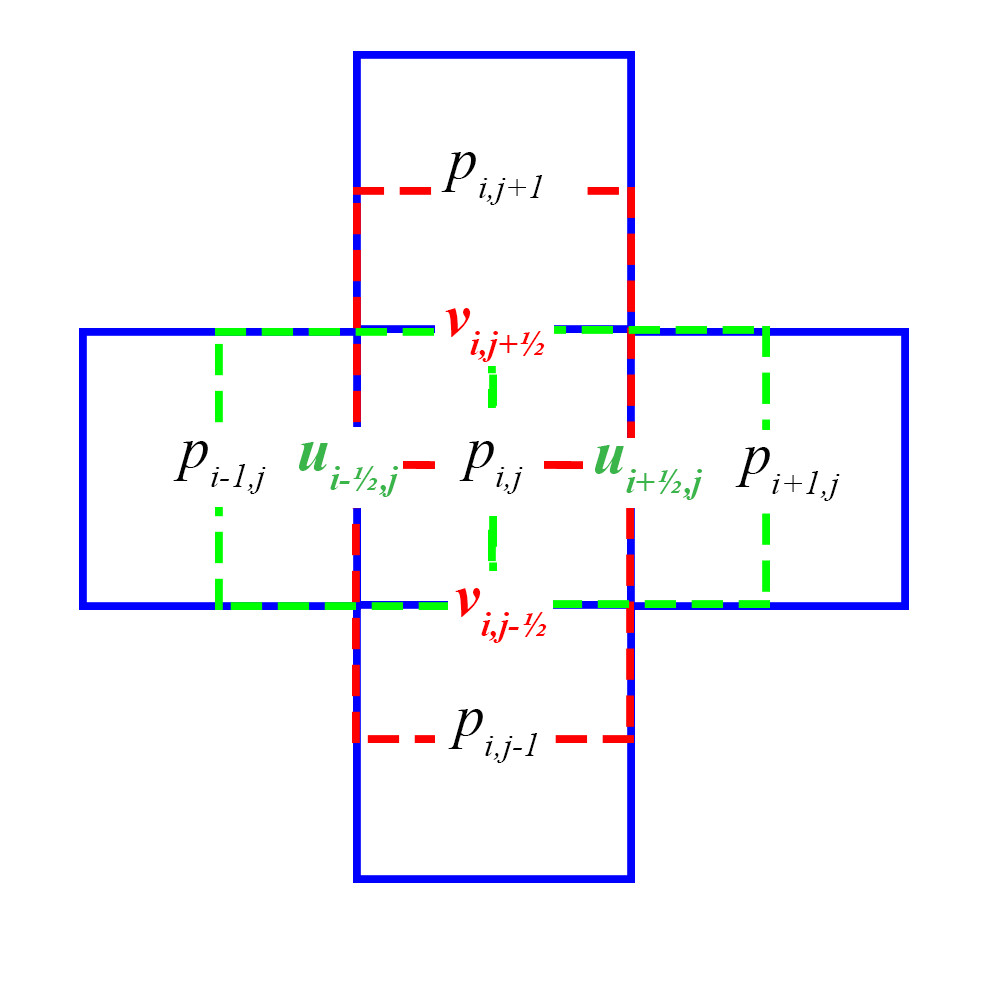} & 
			\includegraphics[width=0.45\textwidth]{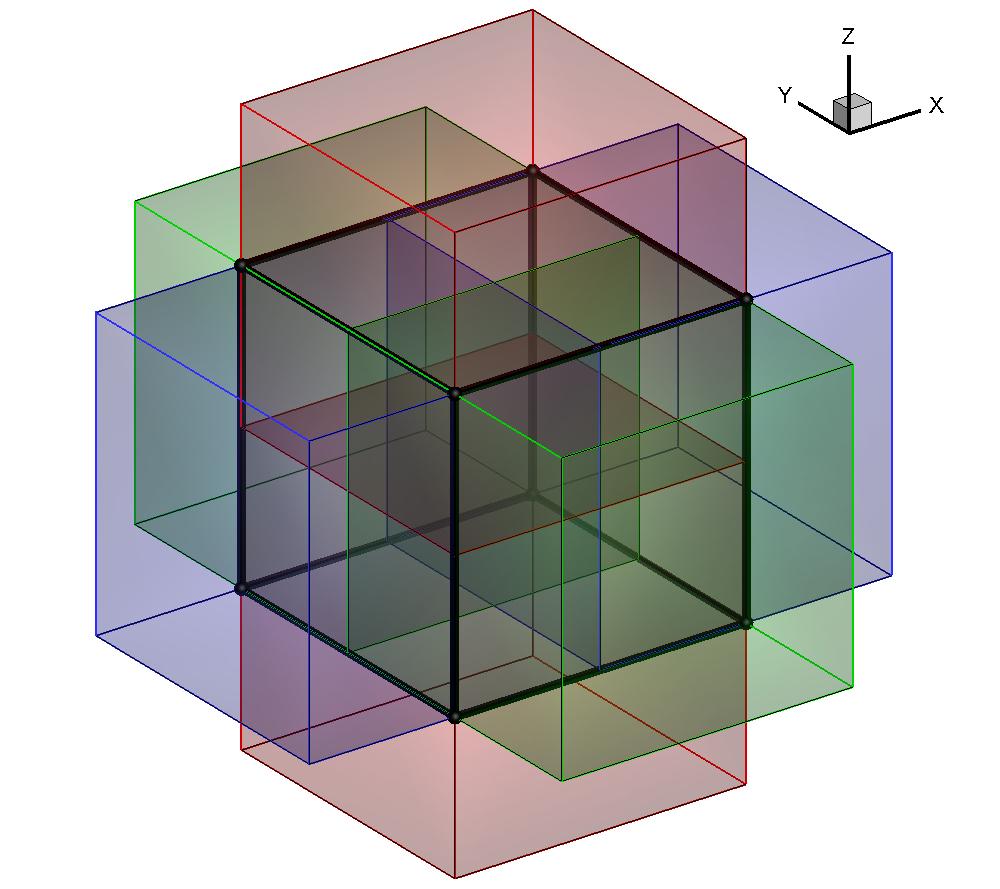}  
\end{tabular} 			
\caption{Mesh-staggering for the two dimensional case (left) and for the three-dimensional case (right).}\label{fig:staggering}
\end{figure}

The discrete solution is defined with respect to the same but shifted polynomial basis $\mathpzc{B}_N=\left\{ \varphi_l(x) \right\}_{l=0,N}$ along its own control volume 
(main or dual) for \emph{each} spatial dimension, having 
\begin{align}
 u_{{h}}(\mathbf{x},t)|_{T_{i+\frac{1}{2},j,k}} = u_{i+\frac{1}{2},j,k}(\mathbf{x},t) = \sum\limits_{l_1,l_2,l_3=0}^N \psi_{l_1}(x)\omega_{l_2}(y)\omega_{l_3}(z) \, \hat{u}_{l,i+\frac{1}{2},j,k}(t) \;\;\; \text{for}\; 
\mathbf{x}\in T_{i+\frac{1}{2},j,k} \label{eq:uDGpoly}\\  %x\in T_{i+\frac{1}{2}},\;y\in T_{j},\;z\in T_{k}\\
\text{for}\;\;\; i=0,...,N_x\;\;j=1,...,N_y\;\;k=1,...,N_z\;\;\nonumber
\end{align}
\begin{align}
 v_{{h}}(\mathbf{x},t)|_{T_{i,j+\frac{1}{2},k}} = v_{i,j+\frac{1}{2},k}(\mathbf{x},t) = \sum\limits_{l_1,l_2,l_3=0}^N \omega_{l_1}(x)\psi_{l_2}(y)\omega_{l_3}(z) \, \hat{v}_{l,i,j+\frac{1}{2},k}(t)\;\;\; \text{for}\; 
\mathbf{x}\in T_{i,j+\frac{1}{2},k}  \label{eq:vDGpoly}\\
\text{for}\;\;\; i=1,...,N_x\;\;j=0,...,N_y\;\;k=1,...,N_z\;\;\nonumber
\end{align}
\begin{align}
 w_{{h}}(\mathbf{x},t)|_{T_{i,j,k+\frac{1}{2}}} = w_{i,j,k+\frac{1}{2}}(\mathbf{x},t) = \sum\limits_{l_1,l_2,l_3=0}^N \omega_{l_1}(x)\omega_{l_2}(y)\psi_{l_3}(z) \, \hat{w}_{l,i,j,k+\frac{1}{2}}(t)\;\;\; \text{for}\; 
\mathbf{x}\in T_{i,j,k+\frac{1}{2}} \label{eq:wDGpoly}\\
\text{for}\;\;\; i=1,...,N_x\;\;j=1,...,N_y\;\;k=0,...,N_z\;\; \nonumber
\end{align}
\begin{align}
 p_{{h}}(\mathbf{x},t)|_{T_{i,j,k}} = p_{i,j,k}(\mathbf{x},t) = \sum\limits_{l_1,l_2,l_3=0}^N \omega_{l_1}(x)\omega_{l_2}(y)\omega_{l_3}(z) \, \hat{p}_{l,i,j,k}(t)\;\;\; 
\text{for}\; \mathbf{x}\in T_{i,j,k}  \label{eq:pDGpoly} \\
\text{for}\;\;\; i=1,...,N_x\;\;j=1,...,N_y\;\;k=1,...,N_z\;\;\nonumber
\end{align}
with the multi-index $l=(l_1,l_2,l_3)$ and where $\hat{u}_{l,i+\frac{1}{2},j,k}(t)$, $\hat{v}_{l,i,j+\frac{1}{2},k}(t)$, $\hat{w}_{l,i,j,k+\frac{1}{2}}(t)$ 
and $\hat{p}_{l,i,j,k}(t)$ (for $0 \leq \max{(l)} \leq N$) are called 
\emph{degrees of freedom} of the corresponding physical variables; as already defined above, $N_x$, $N_y$, and $N_z$ are the number of elements on the main 
grid in the $x$, $y$ and $z$ direction, respectively. The polynomials $\psi$ and $\omega$ are generated from the basis functions $\varphi$ with the rule 
\begin{align*}
 & \psi(s) = \varphi(\xi),    & \text{with} & \qquad s = s_i+\xi \Delta s,                & 0 \leq \xi \leq 1, \\
 & \omega(s) = \varphi(\xi),  & \text{with} & \qquad s = s_{i-\frac{1}{2}}+\xi \Delta s,  & 0 \leq \xi \leq 1,
\end{align*}
where $s$ stands for $x$, $y$ or $z$. A simplified picture of the resulting mesh-staggering is depicted in Figure \ref{fig:staggering} for the two dimensional and 
for the three dimensional case. 
In our particular implementation, the $\varphi(\xi)$ are defined by the Lagrange interpolation polynomials passing through the Gauss-Legendre quadrature points on 
the unit interval $[0,1]$, see \cite{DumbserCasulli2013}, hence leading to an \textit{orthogonal} nodal basis. As a result, all element mass matrices are \textit{diagonal}. 

By direct substitution of the definitions (\ref{eq:uDGpoly}-\ref{eq:pDGpoly}) into (\ref{eq:wNSmom}-\ref{eq:wNSinc}) and by using the same basis functions also as test functions, one obtains 
the following \emph{semi-discrete} staggered DG discretiation of the incompressible Navier-Stokes equations: 
\begin{align}
&\int \limits_{T_{i+\frac{1}{2},j,k}}  \psi_{m_1}(x)\omega_{m_2}(y)\omega_{m_3}(z) \left( \frac{\partial u_{h}}{\partial t} + \nabla \cdot \mathbf{F}_{u} + \partial_x p_{h} \right) d\x = 0,\label{eq:uDGNSmom} \\
&\int \limits_{T_{i,j+\frac{1}{2},k}} \omega_{m_1}(x)\psi_{m_2}(y)\omega_{m_3}(z)  \left( \frac{\partial v_{h}}{\partial t} + \nabla \cdot \mathbf{F}_{v} + \partial_y p_{h} \right) d\x = 0,\label{eq:vDGNSmom}\\
&\int \limits_{T_{i,j,k+\frac{1}{2}}} \omega_{m_1}(x)\omega_{m_2}(y)\psi_{m_3}(z)  \left( \frac{\partial w_{h}}{\partial t} + \nabla \cdot \mathbf{F}_{w} + \partial_z p_{h} \right) d\x = 0, \label{eq:wDGNSmom}\\
&\int \limits_{T_{i,j,k}} \omega_{m_1}(x)\omega_{m_2}(y)\omega_{m_3}(z) \, \left( \nabla \cdot \mathbf{v}_{h}\right) d\x = 0.\label{eq:DGNSinc}
\end{align}
Integration of \eqref{eq:DGNSinc} by parts yields 
\begin{equation}
  \int \limits_{\partial T_{i,j,k}} \omega_{m_1}(x)\omega_{m_2}(y)\omega_{m_3}(z) \, \mathbf{v}_{h} \cdot \vec{n} \, dS 
	- \int \limits_{T_{i,j,k}} \nabla \left(\omega_{m_1}(x)\omega_{m_2}(y)\omega_{m_3}(z) \right) \cdot \mathbf{v}_{h} \, d\x = 0, \label{eq:DGNSinc2}
\end{equation} 
which is well defined, since the velocity vector $\mathbf{v}_h$ is continuous across the element boundary $\partial T_{i,j,k}$, thanks to the use of a staggered grid approach.  
However, because of the staggering, $p_h$ is discontinuous inside the domains of integration of the momentum equations (\ref{eq:uDGNSmom}-\ref{eq:wDGNSmom}) and the following \emph{jump} 
contributions arise % from the integration
\begin{align}
\int\limits_{x_i}^{x_{i+1}} \int \limits_{y_{j-\frac{1}{2}}}^{y_{j+\frac{1}{2}}} \int \limits_{z_{k-\frac{1}{2}}}^{z_{k+\frac{1}{2}}}  
\psi_{m_1}(x) \omega_{m_2}(y) \omega_{m_3}(z) \partial_x p_h(\mathbf{x},t) \, d\x = && \nonumber \\   
\int\limits_{x_i}^{x_{i+\frac{1}{2}}} \int \limits_{y_{j-\frac{1}{2}}}^{y_{j+\frac{1}{2}}} \int \limits_{z_{k-\frac{1}{2}}}^{z_{k+\frac{1}{2}}}  
\psi_{m_1}(x) \omega_{m_2}(y) \omega_{m_3}(z) \partial_x p_{i,j,k}(\mathbf{x},t) \, d\x + 
\int\limits_{x_{i+\frac{1}{2}}}^{x_{i+1}} \int \limits_{y_{j-\frac{1}{2}}}^{y_{j+\frac{1}{2}}} \int \limits_{z_{k-\frac{1}{2}}}^{z_{k+\frac{1}{2}}}  
\psi_{m_1}(x) \omega_{m_2}(y) \omega_{m_3}(z) \partial_x p_{i+1,j,k}(\mathbf{x},t) \, d\x + && \nonumber \\
+ \int \limits_{y_{j-\frac{1}{2}}}^{y_{j+\frac{1}{2}}} \int \limits_{z_{k-\frac{1}{2}}}^{z_{k+\frac{1}{2}}}  
\psi_{m_1}(x_{i+\frac{1}{2}}) \omega_{m_2}(y) \omega_{m_3}(z)  \left(p_{i+1,j,k}(x_{i+\frac{1}{2}},y,z,t) - p_{i,j,k}(x_{i+\frac{1}{2}},y,z,t)\right) dy dz, && \label{eq:pressurejump}%\\
%\int\limits_{y_j}^{y_{j+1}} dy\,\psi_{m'}(y) \partial_y p = \int\limits_{y_j}^{y_{j+\frac{1}{2}}}dy\, \psi_{m'}(y) \partial_y p^j + \int\limits_{y_{j+\frac{1}{2}}}^{y_{j+1}}dy\, \psi_{m'}(y) \partial_y p^{j+1}   + \psi_{m'}(y_{j+\frac{1}{2}})\left(p^{j+1}(y_{j+\frac{1}{2}}) - p^j(y_{j+\frac{1}{2}})\right)\\
%\int\limits_{z_k}^{z_{k+1}}dz\, \psi_{m''}(z) \partial_z p = \int\limits_{z_k}^{z_{k+\frac{1}{2}}}dz\, \psi_{m''}(z) \partial_z p^k + \int\limits_{z_{k+\frac{1}{2}}}^{z_{k+1}}dz\, \psi_{m''}(z) \partial_z p^{k+1}   + \psi_{m''}(z_{k+\frac{1}{2}})\left(p^{k+1}(z_{k+\frac{1}{2}}) - p^k(z_{k+\frac{1}{2}})\right)\\
\end{align}
with similar expressions also in the $y$- and $z$-momentum equations, respectively. 
Thus, an efficient \textit{semi-implicit} time discretization of the governing PDE system is obtained by introducing an \textit{explicit} discretization of the nonlinear convective and viscous terms 
and an \textit{implicit} discretization of the pressure gradients in the momentum equations (\ref{eq:uDGNSmom}-\ref{eq:wDGNSmom}) and of the incompressibility condition (\ref{eq:DGNSinc}). 
After evaluating the integrals and via some manipulations one can obtain the following coupled system of equations for the vectors of the degrees of freedom of velocity 
$\dof{U}^{n+1}_{i+\frac{1}{2},j,k}=\hat{u}^{n+1}_{l,i+\frac{1}{2},j,k}$, $\dof{V}^{n+1}_{i,j+\frac{1}{2},k}=\hat{v}^{n+1}_{l,i,j+\frac{1}{2},k}$, 
$\dof{W}^{n+1}_{i,j,k+\frac{1}{2}}=\hat{w}^{n+1}_{l,i,j,k+\frac{1}{2}}$ and pressure $\dof{P}^{n+1}_{i,j,k}=\hat{p}^{n+1}_{l,i,j,k}$, respectively: 
\begin{align}
\Mass{xyz} \cdot\left( \dof{U}^{n+1}_{i+\frac{1}{2},j,k}- \dof{Fu}^{n}_{i+\frac{1}{2},j,k}  \right) + \frac{\Delta t}{\Delta x}\Mass{yz}\cdot \left( \tild{R}{x} \cdot\dof{P}^{n+\theta}_{i+1,j,k} - \tild{L}{x} \cdot\dof{P}^{n+\theta}_{i,j,k} \right)&=0,\label{eq:uSDG1} \\
\Mass{xyz}\cdot \left( \dof{V}^{n+1}_{i,j+\frac{1}{2},k}- \dof{Fv}^{n}_{i,j+\frac{1}{2},k}  \right) + \frac{\Delta t}{\Delta y}\Mass{zx} \cdot\left( \tild{R}{y} \cdot\dof{P}^{n+\theta}_{i,j+1,k} - \tild{L}{y}\cdot \dof{P}^{n+\theta}_{i,j,k} \right)&=0,\label{eq:vSDG1} \\
\Mass{xyz} \cdot\left( \dof{W}^{n+1}_{i,j,k+\frac{1}{2}}- \dof{Fw}^{n}_{i,j,k+\frac{1}{2}}  \right) + \frac{\Delta t}{\Delta z} \Mass{xy}\cdot\left( \tild{R}{z}\cdot \dof{P}^{n+\theta}_{i,j,k+1} - \tild{L}{z} \cdot\dof{P}^{n+\theta}_{i,j,k} \right)&=0,\label{eq:wSDG1} 
\end{align}
\begin{equation}
 \frac{\Mass{yz} \left(\barr{R}{x} \cdot\dof{U}^{n+1}_{i+\frac{1}{2},j,k} \!-\! \barr{L}{x} \cdot\dof{U}^{n+1}_{i-\frac{1}{2},j,k}\right)}{\Delta x} + 
 \frac{\Mass{zx} \left(\barr{R}{y} \cdot\dof{V}^{n+1}_{i,j+\frac{1}{2},k} \!-\! \barr{L}{y}\cdot \dof{V}^{n+1}_{i,j-\frac{1}{2},k}\right)}{\Delta y} +
 \frac{\Mass{xy} \left(\barr{R}{z} \cdot\dof{W}^{n+1}_{i,j,k+\frac{1}{2}} \!-\! \barr{L}{z} \cdot\dof{W}^{n+1}_{i,j,k-\frac{1}{2}}\right)}{\Delta z} = 0.\label{eq:pSDG1}
\end{equation}
Here, the following matrices have been used
\begin{align*}
&\Mass{} \equiv \left\{ M_{pq} \right\}_{p,q=0,...,N} \equiv \left\{\int \limits_0^1 \,\varphi_p\left(\xi\right) \varphi_q\left(\xi\right)  d\xi \right\}_{p,q=0,...,N} \nonumber \\
&\tild{R}{}  \equiv \left\{ R_{pq} \right\}_{p,q=0,...,N} \equiv \left\{\varphi_p(\frac{1}{2})\varphi_q(0) +\frac{1}{2}\int \limits_0^1  \,\varphi_p\left(\frac{1}{2}+\frac{\xi}{2}\right) \varphi_q'\left(\frac{\xi}{2}\right)d\xi \right\}_{p,q=0,...,N} \nonumber\\
&\tild{L}{}   \equiv \left\{ L_{pq} \right\}_{p,q=0,...,N} \equiv \left\{\varphi_p(\frac{1}{2})\varphi_q(1) -\frac{1}{2}\int \limits_0^1 \,\varphi_p \left( \frac{\xi}{2}\right)\varphi_q '\left(\frac{1}{2}+\frac{\xi}{2}\right) d\xi\right\}_{p,q=0,...,N} \nonumber\\
\end{align*}
\begin{align}
&\barr{R}{}   \equiv \left\{ L^T_{pq} \right\}_{p,q=0,...,N} \equiv \left\{\varphi_p(1)\varphi_q(\frac{1}{2}) - \frac{1}{2}\int \limits_0^1 \,\varphi_p'\left(\frac{1}{2}+\frac{\xi}{2}\right) \varphi_q\left( \frac{\xi}{2}\right)  d\xi \right\}_{p,q=0,...,N} \nonumber\\
&\barr{L}{} \equiv \left\{ R^T_{pq} \right\}_{p,q=0,...,N} \equiv \left\{\varphi_p(0)\varphi_q(\frac{1}{2})  + \frac{1}{2}\int \limits_0^1 \,\varphi_p'\left( \frac{\xi}{2}\right) \varphi_q \left(\frac{1}{2}+\frac{\xi}{2}\right) d\xi\right\}_{p,q=0,...,N}, \label{eq:matrices}
\end{align}
which operate along a generic vector of degrees of freedom 
$%\begin{align}
\dof{X} = \left\{\hat{x}_{mm'm''}\right\}_{m,m',m''=0,..,N},
$%\end{align} 
via the tensor products 
\begin{align}
\operator{Z}{x}\cdot \dof{X} = Z_{ml}I_{m'l'}I_{m''l''} \hat{x}_{l\,l'l''}, \qquad 
\operator{Z}{y}\cdot \dof{X} = I_{ml}Z_{m'l'}I_{m''l''} \hat{x}_{l\,l'l''}, \qquad 
\operator{Z}{z} \cdot\dof{X} = I_{ml}I_{m'l'}Z_{m''l''} \hat{x}_{l\,l'l''} \label{eq:tensorP1}
\end{align} 
and
%\begin{align} 
%\operator{Z}{x$_i$ x$_j$}  \equiv \operator{Z}{x$_i$}\cdot\operator{Z}{x$_j$},
%\end{align} 
\begin{align}
\begin{array}{rl} 
\operator{Z}{x$_i$ x$_j$}  &  \equiv \operator{Z}{x$_i$}\cdot\operator{Z}{x$_j$},\\
\operator{Z}{x$_i$ x$_j$ x$_k$}  &\equiv \operator{Z}{x$_i$}\cdot\operator{Z}{x$_j$}\cdot\operator{Z}{x$_k$}, \end{array}
 \;\;\; \text{with}\; x_i,x_j,x_k \in \left\{x,y,z\right\}  \label{eq:tensorP2}
\end{align}
where $\mathbf{Z}$ is a real square matrix, $I$ is the identity operator and the Einstein convention of summation over repeated indexes is assumed. Note that for the pressure gradients an implicitness factor $\theta \in [\frac{1}{2},1]$ 
has been introduced, by defining $\dof{P}^{n+\theta}=\theta\dof{P}^{n+1} + \left(1-\theta\right)\dof{P}^n$.  By choosing $\theta=\left.1\middle/2\right.$, the time discretization of 
(\ref{eq:uSDG1})-(\ref{eq:pSDG1}) is equivalent to a Crank-Nicolson scheme, which is \textit{second-order accurate} in time.

$\dof{Fu}^n$, $\dof{Fv}^n$ and $\dof{Fw}^n$ can be computed with any suitable explicit discretization for advection and diffusion. An insight into these terms will be given later in the text.

The coupled system of equations (\ref{eq:uSDG1})-(\ref{eq:pSDG1}) \reds{has a typical saddle point structure that arises naturally from the discretization 
of the incompressible Navier-Stokes equations}. Its direct solution can be cumbersome, since it involves four unknown quantities: three velocity components and the scalar 
pressure. The complexity of the problem can be considerably reduced with a very simple manipulation. After multiplying the momentum equations by the inverse of the mass matrix $\Mass{xyz}$, the 
discrete velocity equations can be substituted 
% 
%following  equations for the velocity components are obtained
%\begin{align}
  %\dof{U}^{n+1}_{i+\frac{1}{2},j,k} & = \dof{Fu}^{n}_{i+\frac{1}{2},j,k} - \frac{\Delta t}{\Delta x} \iMass{x}\cdot \left( \tild{R}{x} \cdot\dof{P}^{n+\theta}_{i+1,j,k} - \tild{L}{x} \cdot\dof{P}^{n+\theta}_{i,j,k} \right),\label{eq:uSDG2} \\
%%
  %\dof{V}^{n+1}_{i,j+\frac{1}{2},k} & =  \dof{Fv}^{n}_{i,j+\frac{1}{2},k}   - \frac{\Delta t}{\Delta y} \iMass{y} \cdot\left( \tild{R}{y} \cdot\dof{P}^{n+\theta}_{i,j+1,k} - \tild{L}{y}\cdot \dof{P}^{n+\theta}_{i,j,k} \right),\label{eq:vSDG2} \\
%%
 %\dof{W}^{n+1}_{i,j,k+\frac{1}{2}} & =  \dof{Fw}^{n}_{i,j,k+\frac{1}{2}}   - \frac{\Delta t}{\Delta z} \iMass{z} \cdot\left( \tild{R}{z}\cdot \dof{P}^{n+\theta}_{i,j,k+1} - \tild{L}{z} \cdot\dof{P}^{n+\theta}_{i,j,k} \right),\label{eq:wSDG2} 
%\end{align}
into the discrete incompressibility condition (\ref{eq:pSDG1}). As a result, one obtains one single linear system for the degrees of freedom of the unknown scalar pressure $\dof{P}^{n+1}$ only, i.e. 
\begin{align}
 \frac{\theta \Delta t}{\Delta x^2}\left(\Mass{yz}\cdot \Hoperator{R}{x} \right)\cdot\dof{P}^{n+1}_{i+1,j,k} + \frac{\theta \Delta t}{\Delta y^2}\left(\Mass{zx} \cdot \Hoperator{R}{y} \right)\cdot\dof{P}^{n+1}_{i,j+1,k}+ \frac{\theta \Delta t}{\Delta z^2}\left(\Mass{xy}\cdot  \Hoperator{R}{z} \right)\cdot\dof{P}^{n+1}_{i,j,k+1} + \nonumber \\
+ \left( \frac{\theta \Delta t}{\Delta x^2}\Mass{yz}\cdot  \Hoperator{C}{x} + \frac{\theta\Delta t}{\Delta y^2}\Mass{zx} \cdot \Hoperator{C}{y}+ \frac{\theta \Delta t}{\Delta z^2}\Mass{xy} \cdot \Hoperator{C}{z}\right)\cdot \dof{P}^{n+1}_{i,j,k} + \nonumber \\
 + \frac{\theta \Delta t}{\Delta x^2}\left(\Mass{yz}\cdot  \Hoperator{L}{x} \right)\cdot\dof{P}^{n+1}_{i-1,j,k} + \frac{\theta\Delta t}{\Delta y^2}\left(\Mass{zx} \cdot \Hoperator{L}{y} \right)\cdot\dof{P}^{n+1}_{i,j-1,k}+ \frac{\theta \Delta t}{\Delta z^2}\left(\Mass{xy} \cdot \Hoperator{L}{z} \right)\cdot\dof{P}^{n+1}_{i,j,k-1} & = \dof{b}^n_{i,j,k},\label{eq:pSyst}\\
\text{for}\;\;\;i=2,...,N_x-1;\;\;\;j=2,...,N_y-1;\;\;\,k=2,...,N_z-1 \nonumber
\end{align}
Here, the following new tensors have been defined 
\begin{align}
\Hoperator{R}{} = - \left(\barr{R}{}\cdot\iMass{}\cdot\tild{R}{} \right),\;\;\; %\nonumber\\
\Hoperator{L}{} = - \left(\barr{L}{}\cdot\iMass{}\cdot\tild{L}{}\right),\;\;\; % \nonumber\\
\Hoperator{C}{} = + \left(\barr{R}{}\cdot\iMass{}\cdot\tild{L}{}\right) + \left(\barr{L}{}\cdot\iMass{}\cdot\tild{R}{} \right).
\end{align}
System (\ref{eq:pSyst}) can be written in compact form as 
$%\begin{align}
\mathbb{H}\cdot \mathpzc{P}^{n+1} = \mathpzc{b}^n %\\
%\text{with}\;\;\; \mathbb{H}\in \left[\varmathbb{R}^{d(N+1)}\right]^{N_e}\times\left[\varmathbb{R}^{d(N+1)}\right]^{N_e}\;\;\;\text{and}\;\;\;\mathbb{H}\in \left[\varmathbb{R}^{d(N+1)}\right]^{N_e}
$%\end{align}
, where $\mathbb{H}$ is the block coefficient matrix, $\mathpzc{P}^{n+1}$ collects all the unknown pressure degrees of freedom of the computational domain at the new time level 
and $\mathpzc{b}^n$ collects all the known terms of the equations. All the real advantages of the chosen mesh-staggering and the semi-implicit discretization arise in the particular 
features of the resulting linear system (\ref{eq:pSyst}). 
\reds{ 
The substitution of the discrete velocity equation into the discrete divergence condition can be seen as the application of the Schur complement to the saddle point problem
(\ref{eq:uSDG1})-(\ref{eq:pSDG1}).   
The particular grid staggering used in this paper, i.e. the so-called C-grid according to the nomenclature of Arakawa \& Lamb \cite{Arakawa}, has been selected to be the one that minimizes the stencil  size of the resulting pressure system.}\footnote{\reds{ Without staggering (A-grid case), the integral of the pressure gradients in the momentum  equations (\ref{eq:uDGNSmom}-\ref{eq:wDGNSmom}), after integrating by parts, would generate a three-point stencil of dependence between the elements by means of some numerical flux functions that are necessary for approximating the pressure at the element  interfaces, i.e. $p(x_{i+1/2})=\mathcal{G}(p_i,p_{i+1})$. With the same argument, further flux functions are needed also in the incompressibility condition (\ref{eq:DGNSinc}) for evaluating the velocities at the interfaces and the resulting discrete pressure system  would become: block $5$\emph{-diagonal} for the $1$d case, instead of being block $3$\emph{-diagonal}; block $9$\emph{-diagonal} for the two  dimensional case, versus our block $5$\emph{-diagonal} system; block $13$\emph{-diagonal} for the three dimensional case, versus our block $7$\emph{-diagonal} system. Concerning the vertex-based staggered grids (B-grid), Riemann solvers or numerical flux functions are not necessary. However, with a vertex based staggering, a block $9$\emph{-diagonal system} or a block $27$\emph{-diagonal system}  are  obtained for the two and for the three dimensional case, respectively, see also Table \ref{tab:StencilSize}}}  
In fact, $\mathbb{H}$ is only \textbf{\emph{block}} \textbf{\emph{hepta}}-\textbf{\emph{diagonal}} for the three 
dimensional case, and only \textbf{\emph{block}} \textbf{\emph{penta}}-\textbf{\emph{diagonal}} for the two dimensional case.  
\reds{ Table \ref{tab:StencilSize} shows the stencil-sizes (number of non-zero blocks) of the resulting algebraic systems for the pressure, varying for different choices of the grid type and for 
different numbers of space dimensions.}
In particular, the \textbf{symmetry} of $\mathbb{H}$ 
can be easily proven by showing directly from the definition in eq. (\ref{eq:pSyst}) that $\mathbb{H} = \mathbb{H}^T$. The key point of the demonstration is that the next three 
equivalences are true by construction of (\ref{eq:matrices}) 
\begin{align}
\Mass{T}	&= \Mass{}\nonumber\\
\Hoperator{R}{T} &= - \left(\barr{R}{}\cdot\iMass{}\cdot\tild{R}{}\right)^T =  - \left(\tild{R}{T}\cdot \left(\iMass{}\right)^{T}\cdot\barr{R}{T}\right)=  - \left(\barr{L}{}\cdot\iMass{}\cdot\tild{L}{}\right) \equiv \Hoperator{L}{},  \nonumber\\
\Hoperator{C}{T} &= \left(\barr{R}{}\cdot\iMass{}\cdot\tild{L}{}\right)^T + \left(\barr{L}{}\cdot\iMass{}\cdot\tild{R}{} \right)^T =  \nonumber\\
&=  \left(\tild{L}{T}\cdot\left(\iMass{}\right)^{T}\cdot\barr{R}{T}\right) + \left(\tild{R}{T}\cdot\left(\iMass{}\right)^{T}\cdot\barr{L}{T} \right) =  \nonumber\\
&=  \left(\barr{R}{}\cdot\iMass{}\cdot\tild{L}{}\right) + \left(\barr{L}{}\cdot\iMass{}\cdot\tild{R}{} \right) \equiv \Hoperator{C}{}.
\end{align}
Further to that, it can be shown that $\mathbb{H}$ is also \textbf{\emph{positive semi-definite}} in the general case, i.e. 
\begin{align}
\mathpzc{v}^T\mathbb{H}\mathpzc{v} \geq 0, \qquad \forall \mathpzc{v} \neq 0.  %\left[\varmathbb{R}^{N+1}\times\varmathbb{R}^{N+1}\times\varmathbb{R}^{N+1}\right]^{N_e}
\end{align}
%In particular, the the validity of (\ref{eq:posdef}) is verified if
%\begin{align}
%\mathpzc{v}^T\mathbb{H}^{i,j,k}\mathpzc{v} \geq 0 \;\;\; \text{for every} \;\;\; \mathpzc{v}\in \left[\varmathbb{R}^{d(N+1)}\right]^{N_e} \label{eq:posdef2} %\left[\varmathbb{R}^{N+1}\times\varmathbb{R}^{N+1}\times\varmathbb{R}^{N+1}\right]^{N_e}
%\end{align}
%where $\mathbb{H}^{i,j,k}$ is the $(i,j,k,)$-th block-row of $\mathbb{H}$.
Note that in this notation, equation (\ref{eq:pSyst}) can be written as $\mathbb{H}^{i,j,k}\cdot \mathpzc{P}^{n+1} = \mathpzc{b}^n_{i,j,k}$. %Moreover,  (\ref{eq:posdef2}) is valid if a stronger inequality holds, i.e.
%\begin{align}
%\mathpzc{v}^T_{p,q,r}\mathbb{H}^{i,j,k}_{p,q,r}\mathpzc{v}_{p,q,r} \geq 0 \;\;\; \text{for every} \;\;\; \mathpzc{v}_{p,q,r}\in \varmathbb{R}^{d(N+1)}\label{eq:posdef3} 
%\end{align}
%that means showing the validity of (\ref{eq:posdef3}) for the three cases
%\begin{align}
%\mathbb{H}^{i,j,k}_{p,q,r} = \left\{ \begin{array}{l} \Hoperator{R} \\ \Hoperator{L} \\ \Hoperator{C} \right\. \label{eq:posdefcases}
%\end{align}
%Note that in (\ref{eq:posdefcases}) the mass matrix tensors, that are present in equation (\ref{eq:pSyst}), is omitted because of obvious  motivations.
Matrix $\mathbb{H}$ can be written in the form of a tensor product of the matrices $\mathbb{H}=\mathbb{H}^{xyz} =\mathbb{H}^x\mathbb{H}^y\mathbb{H}^z$. Next, the positive semi-definiteness 
is shown to be valid for the one-dimensional case $\mathbb{H}=\mathbb{H}^{x}$, then the extension to $\mathbb{H}=\mathbb{H}^{xyz}$ is straightforward. 
If $d=1$ and periodic boundary conditions are assumed the left hand side of (\ref{eq:pSyst}) can be written as 
\begin{align}
\mathbb{H}^{i}\cdot \mathpzc{P}^{n+1} = & - \left[\left(\barr{R}{}\cdot\iMass{}\cdot\tild{R}{} \right)  \cdot\dof{P}^{n+1}_{i+1} -  \left(\barr{R}{}\cdot\iMass{}\cdot\tild{L}{}\right) \dof{P}^{n+1}_{i}\right] + \nonumber \\
& + \left[ \left(\barr{L}{}\cdot\iMass{}\cdot\tild{R}{} \right) \dof{P}^{n+1}_{i}    - \left(\barr{L}{}\cdot\iMass{}\cdot\tild{L}{}\right)\dof{P}^{n+1}_{i-1}\right] ,\label{eq:pSyst1d}\\
&\text{for}\;\;\;i=1,...,N_x, \;\;\;\text{with}\;P_{N_x+1} =P_1 ;\nonumber
%
%&\Hoperator{R}{} = - \left(\barr{R}{}\cdot\iMass{}\cdot\tild{R}{} \right)\nonumber\\
%&\Hoperator{L}{} = - \left(\barr{L}{}\cdot\iMass{}\cdot\tild{L}{}\right) \nonumber\\
%&\Hoperator{C}{} = + \left(\barr{R}{}\cdot\iMass{}\cdot\tild{L}{}\right) + \left(\barr{L}{}\cdot\iMass{}\cdot\tild{R}{} \right).
\end{align}
where the mass matrix tensors and the discretization constants can be removed after multiplication with suitable factors from the left. 
Now a new nomenclature is introduced to emphasize the features of the system,  
\begin{align}
&\mathcal{R} =   \tild{R}{} \equiv \barr{L}{T} ,\;\;\;\mathcal{L} =   \tild{L}{} \equiv \barr{R}{T} ,\label{eq:matricesLR}
%&\mathcal{R} = \mathbf{M}^{-\frac{1}{2}} \tild{R}{} \equiv \barr{L}{T}\left(\mathbf{M}^{-\frac{1}{2}}\right)^T,\nonumber\\
%&\mathcal{L} = \mathbf{M}^{-\frac{1}{2}} \tild{L}{} \equiv \barr{R}{T} \left(\mathbf{M}^{-\frac{1}{2}}\right)^T,\nonumber
\end{align}
from the definitions (\ref{eq:matrices}). Now, relation (\ref{eq:pSyst1d}) can be written as 
\begin{align}
\mathbb{H}^{i}\cdot \mathpzc{P}^{n+1} = & - \left[\mathcal{L}^T\mathbf{M}^{-1}\mathcal{R}  \cdot\dof{P}^{n+1}_{i+1} -  \mathcal{L}^T\mathbf{M}^{-1}\mathcal{L} \cdot\dof{P}^{n+1}_{i}\right] + %\nonumber \\
   \left[\mathcal{R}^T\mathbf{M}^{-1}\mathcal{R} \cdot\dof{P}^{n+1}_{i}    - \mathcal{R}^T\mathbf{M}^{-1}\mathcal{L} \cdot\dof{P}^{n+1}_{i-1}\right] ,\nonumber \\ %\label{eq:pSyst}\\
%\mathpzc{P}^{n+1}_i\cdot\mathbb{H}^{i}\cdot \mathpzc{P}^{n+1} = & - \dof{P}^{n+1}_i\cdot\left[\mathcal{L}^T\mathcal{R}  \cdot\dof{P}^{n+1}_{i+1} -  \mathcal{L}^T\mathcal{L} \cdot\dof{P}^{n+1}_{i}\right] + \nonumber \\
%& +  \dof{P}^{n+1}_i\cdot\left[\mathcal{R}^T\mathcal{R} \cdot\dof{P}^{n+1}_{i}    - \mathcal{R}^T\mathcal{L} \cdot\dof{P}^{n+1}_{i-1}\right] ,\label{eq:pSyst1d}\\
&\text{for}\;\;\;i=1,...,N_x, \;\;\;\text{with}\;P_{N_x+1} =P_1 ;\nonumber
\end{align}
Then, the global system can be written as 
\begin{align}
\mathbb{H} \cdot \mathpzc{P}^{n+1}& \equiv \left(\begin{array}{ccccc} 
-\mathcal{L}^T & 0 & \cdots 			& 0      			& \mathcal{R}^T	\\
\mathcal{R}^T 	& -\mathcal{L}^T & 0		& 0 			& 0     			 	\\
0           	& \mathcal{R}^T & -\mathcal{L}^T & \ddots  			& \vdots     	\\
\vdots  			& 0 & \ddots 			   & \ddots 				& 0	    	\\
0 &  \cdots    &  0	 	 			& \mathcal{R}^T	& -\mathcal{L}^T     	
\end{array}\right)\cdot \mathbb{M}^{-1} \cdot \left(\begin{array}{ccccc} 
-\mathcal{L} & \mathcal{R}           & 0 			& \cdots      			& 0 	\\
0 & -\mathcal{L} & \mathcal{R}      			& 0  			& 0     			 	\\
\vdots           & 0 & -\mathcal{L} & \ddots 			& 0     	\\
0  		& 0 			& \ddots 			& \ddots 				& \mathcal{R}    	\\
\mathcal{R} 					&  0     & \cdots 		 			& 0	& -\mathcal{L}     	
\end{array}\right)\cdot \left(\begin{array}{c} 
\dof{P}^{n+1}_{1} \\   	
\dof{P}^{n+1}_{2} \\   	
\vdots \\   	
\dof{P}^{n+1}_{N_x-1} \\   	
\dof{P}^{n+1}_{N_x}   	
\end{array}\right) = \mathpzc{b}^{n},\label{eq:PerPrSyst}
\end{align}
where the \emph{diagonal} mass matrix $\mathbb{M}$ has been introduced.
Matrix $\mathbb{H}$ is proved to be \textbf{\emph{positive semi-definite}} because it can be decomposed into the matrix product
\begin{align}
%&\mathpzc{v}^T\mathbb{H}\mathpzc{v} \equiv \mathpzc{v}^T\mathbb{A}^T\mathbb{A} \mathpzc{v}  = \mathpzc{w}^T  \mathpzc{w} \geq 0, \;\;\; \text{with}\;\;\;  \mathpzc{w} = \mathbb{A} \mathpzc{v}  \\ 
&\mathpzc{v}^T\mathbb{H}\mathpzc{v} \equiv \mathpzc{v}^T\mathcal{D} ^T \mathbb{M}^{-1}\mathcal{D}  \mathpzc{v}  = \mathpzc{w}^T   \mathbb{M}^{-1} \mathpzc{w} \geq 0, \;\;\; \text{with}\;\;\;  
\mathpzc{w} = \mathcal{D}  \mathpzc{v} \qquad  % \\ 
%&\mathbb{H}=\mathbb{A}^T\cdot\mathbb{A},\\
 \forall \; \mathpzc{v}  \nonumber
\end{align}
because the mass matrix is positive definite. Notice that
\begin{align}
\mathcal{D} =\left(\begin{array}{ccccc} 
-\mathcal{L} & \mathcal{R} & 0 			& \cdots      			& 0 	\\
0 & -\mathcal{L} & \mathcal{R}      			& 0  			& 0     			 	\\
\vdots           & 0 & -\mathcal{L} & \ddots 			& 0     	\\
0  		& 0 			& \ddots 			& \ddots 				& \mathcal{R}    	\\
\mathcal{R} 					&  0     & \cdots 		 			& 0	& -\mathcal{L}    	
\end{array}\right)\label{eq:weakgradientop}
\end{align}
is precisely the weak form of the gradient operator, in fact
\begin{align}
\partial_x p \xrightarrow{\;\;\;\int_{d\Omega^*} \psi \cdot\;\;\;} \mathcal{D} ^i \,\mathpzc{P} = \mathcal{R} \cdot\dof{P}_{i+1} - \mathcal{L}\cdot \dof{P}_i \label{eq:weakgradient}
\end{align}
This is an interesting property because the problem of the uniqueness of the solutions of the pressure system is shifted to the uniqueness of the solutions of 
\begin{align}
\mathcal{D} \, \mathpzc{P}  = \text{right hand side}.
\end{align}
that in ensured in general \emph{up to the solutions of}
$%\begin{align}
\mathcal{D} \, \mathpzc{P}  = 0 %\mathcal{D}_{\textsl{w}}
$. %\end{align}
This means that (for periodic boundaries) the discrete pressure $\mathpzc{P}$ is defined up to weak solutions of $\partial_x p = 0$, which is exactly what one could 
expect from a discrete formulation of the incompressible Navier-Stokes equations. If pressure boundary conditions are specified, it can be verified easily that the resulting
system for the pressure is indeed \textit{symmetric} and \textit{positive-definite}. We further observe all in our numerical experiments that the pressure system seems to be
reasonably well conditioned, since the conjugate gradient method converges in rather few iterations even \textit{without} the use of \textit{any preconditioner}. This is a rather unique feature among existing implicit DG schemes. 

\reds{Future research is concerned with the theoretical analysis of the condition number of the resulting linear systems of our method and the design of specific 
preconditioners for Krylov subspace solvers, using the theory of matrix-valued symbols and Generalized Locally Toeplitz (GLT) algebras, see \cite{serra1998,GSz,glt,TyZ}. }

Finally, once the pressure $\mathpzc{P}^{n+1}$ has been computed, the velocities can be updated directly accordingly to equations (\ref{eq:uSDG1}-\ref{eq:wSDG1}).

\subsection{Explicit discretization of the nonlinear convective and viscous terms}
For an explicit discretization of the nonlinear convective and viscous terms, a standard DG scheme based on the Rusanov flux (local Lax-Friedrichs flux 
\cite{Rusanov1961a,Toro99}) can be adopted on the main grid, see also \cite{MunzDiffusionFlux,DumbserNSE,Hidalgo2011} for numerical flux functions in 
the presence of physical viscosity: 
\begin{align}
\dof{Fu}^{n}_{i,j,k} &= \dof{U}^n_{i,j,k} - \frac{\Delta t}{\Delta x \Delta y \Delta z} \left(\Mass{xyz}\right)^{-1} \cdot\left(\, \int \limits_{\partial T_{i,j,k}} \boldsymbol{\omega} \mathbf{F}_u\cdot \vec{n} dS- \int \limits_{T_{i,j,k}} \nabla \boldsymbol{\omega} \cdot \mathbf{F}_u \, d\x \right) \label{eq:Fu}\\
\dof{Fv}^{n}_{i,j,k} &= \dof{V}^n_{i,j,k} - \frac{\Delta t}{\Delta x \Delta y \Delta z} \left(\Mass{xyz}\right)^{-1} \cdot\left(\, \int \limits_{ \partial T_{i,j,k}} \boldsymbol{\omega} \mathbf{F}_v\cdot \vec{n} dS - \int \limits_{T_{i,j,k}} \nabla \boldsymbol{\omega} \cdot \mathbf{F}_v \, d\x \right) \label{eq:Fv} \\
\dof{Fw}^{n}_{i,j,k} &= \dof{W}^n_{i,j,k} - \frac{\Delta t}{\Delta x \Delta y \Delta z} \left(\Mass{xyz}\right)^{-1} \cdot\left(\, \int \limits_{ \partial T_{i,j,k}} \boldsymbol{\omega} \mathbf{F}_w\cdot \vec{n} dS- \int \limits_{T_{i,j,k}} \nabla \boldsymbol{\omega} \cdot \mathbf{F}_w \, d\x \right)  \label{eq:Fw},
\end{align}
where the numerical flux has the following simple form 
\begin{align}
\mathbf{F}_q\cdot\vec{n} &= \frac{1}{2}\left(\mathbf{F}^+_q + \mathbf{F}^-_q \right)\cdot\vec{n} - \frac{1}{2}s_q \left(q^+-q^-\right)\;\;\; \text{with}\;q=u,v,w. %,\\
%\mathbf{F}_v\cdot\vec{n} &= \frac{1}{2}\left(\mathbf{F}^+_v+\mathbf{F}^-_v\right)\cdot\vec{n} - \frac{1}{2}s_v\left(v^+-v^-\right),\\
%\mathbf{F}_w\cdot\vec{n} &= \frac{1}{2}\left(\mathbf{F}^+_w+\mathbf{F}^-_w\right)\cdot\vec{n} - \frac{1}{2}s_w\left(w^+-w^-\right). 
\end{align}
where $s_q$ are the maximum eigenvalues of the Jacobian of the convective and viscous flux tensor 
\begin{align}
s_{q} &= 2 \, \text{max}\left(|q^+|,|q^-|\right) + 2\nu\frac{2N+1}{\Delta x_q\sqrt{\left.{\pi}\middle/{2}\right. }}\;\;\;\;\;\; \text{with} \; (q,\Delta x_q) = (u,\Delta x), (v,\Delta y),(w,\Delta z). %\\
%s_{v} &= 2 \,\text{max}\left(v^+,v^-\right) + 2\nu\frac{2N+1}{\Delta y\sqrt{\left.{\pi}\middle/{2}\right. }},\\
%s_{w} &= 2 \,\text{max}\left(w^+,w^-\right) + 2\nu\frac{2N+1}{\Delta z\sqrt{\left.{\pi}\middle/{2}\right. }}.
\end{align}
A linear transformation that allows to compute (\ref{eq:Fu}-\ref{eq:Fw}) with the velocity polynomials centered in the main grid is the $L_2$-projection
\begin{align}
\dof{U}^{n}_{i,j,k} &= \left(\Mass{x}\right)^{-1} \cdot\left( \mathbf{M}^{Lx} \cdot \dof{U}^{n}_{i-\frac{1}{2},j,k} + \mathbf{M}^{Rx}\cdot \dof{U}^{n}_{i+\frac{1}{2},j,k} \right), \\
\dof{V}^{n}_{i,j,k} &= \left(\Mass{y}\right)^{-1} \cdot\left( \mathbf{M}^{Ly} \cdot \dof{V}^{n}_{i,j-\frac{1}{2},k} + \mathbf{M}^{Ry}\cdot \dof{V}^{n}_{i,j+\frac{1}{2},k} \right), \\
\dof{W}^{n}_{i,j,k} &= \left(\Mass{z}\right)^{-1} \cdot\left( \mathbf{M}^{Lz} \cdot \dof{W}^{n}_{i,j,k-\frac{1}{2}} + \mathbf{M}^{Rz}\cdot \dof{W}^{n}_{i,j,k+\frac{1}{2}} \right),
\label{eqn.l2forward} 
\end{align}
with
\begin{align}
\mathbf{M}^L = \left\{ M^L_{pq}\right\}_{p,q=0,..,N} = \left\{\frac{1}{2}\int \limits_0^1 \varphi_p\left(\frac{\xi}{2}\right)\varphi_q\left(\frac{1}{2}+\frac{\xi}{2}\right) \,d\xi\right\}_{p,q=0,..,N} \nonumber\\
\mathbf{M}^R = \left\{ M^R_{pq}\right\}_{p,q=0,..,N} = \left\{\frac{1}{2}\int \limits_0^1 \varphi_p\left(\frac{1}{2}+\frac{\xi}{2}\right)\varphi_q\left(\frac{\xi}{2}\right) \,d\xi\right\}_{p,q=0,..,N}\nonumber
\end{align}
Once the advection-diffusion terms have been computed on the main grid, they are projected back to the dual grid with 
\begin{align}
\dof{Fu}^{n}_{i+\frac{1}{2},j,k} &= \left(\Mass{x}\right)^{-1} \cdot\left( \mathbf{M}^{Lx} \cdot \dof{Fu}^{n}_{i,j,k} + \mathbf{M}^{Rx}\cdot \dof{Fu}^{n}_{i+1,j,k} \right), \\
\dof{Fv}^{n}_{i,j+\frac{1}{2},k} &= \left(\Mass{y}\right)^{-1} \cdot\left( \mathbf{M}^{Ly} \cdot \dof{Fv}^{n}_{i,j,k} + \mathbf{M}^{Ry}\cdot \dof{Fv}^{n}_{i,j+1,k} \right), \\
\dof{Fw}^{n}_{i,j,k+\frac{1}{2}} &= \left(\Mass{z}\right)^{-1} \cdot\left( \mathbf{M}^{Lz} \cdot \dof{Fw}^{n}_{i,j,k} + \mathbf{M}^{Rz}\cdot \dof{Fw}^{n}_{i,j,k+1} \right).
\label{eqn.l2backward} 
\end{align}
Since a simple first order Euler time discretization is likely to become linearly unstable, a classical third order TVD Runge-Kutta scheme is used \cite{cockburn_2001_rkd,DumbserCasulli2013,Shu88,TavelliDumbser2014,TavelliDumbser2014b,TavelliDumbser2015}. The explicit discretization has to satisfy a CFL-type time step restriction 
\begin{align}
\Delta t = \text{CFL}\left[\left(2N+1\right)\left(\frac{|u_\text{max}|}{\Delta x} + \frac{|v_\text{max}|}{\Delta y} + \frac{|w_\text{max}|}{\Delta z}\right) + \left(2N+1\right)^2\left(\frac{2\nu}{\Delta x^2} + \frac{2\nu}{\Delta y^2}+\frac{2\nu}{\Delta z^2}\right)\right]^{-1}, \label{eq:CFL}
\end{align}
with $0<$CFL$< 1$. %\left.1\middle/d\right.$.

\subsection{Implicit diffusion}
The discrete formulation of advection and diffusion (\ref{eq:Fu}-\ref{eq:Fw}) is an \emph{explicit} discretization of the equation 
\begin{align}
&\frac{\partial \mathbf{v}}{\partial t} + \nabla \cdot \mathbf{F} = 0.\label{eq:NSmom2}
\end{align}
The time-step restriction (\ref{eq:CFL}) can become rather severe, in particular for highly refined meshes and large values of the kinematic viscosity. An important improvement that allows 
the time-step restriction to become \emph{independent} of the kinematic viscosity is achieved by taking advantage from an \emph{implicit discretization} of the viscous terms. 
In the following, a novel \emph{semi-implicit} numerical method for the advection-diffusion problem is described. 
An efficient semi-implicit discretization of equation (\ref{eq:NSmom2}) is obtained by considering the velocity polynomials to be centered in the \emph{main} grid, and the velocity gradient 
(i.e. the stress tensor $\nu \nabla \mathbf{v}$) to be defined on the edge-based \emph{dual} grid. 
The use of the staggered control volumes for the stress tensor leads to a continuous function $\mathbf{F}_{d}$ across the cell interfaces of the main grid. 
Hence, integration over the control volume $T_{i,j,k}$ yields 
\begin{align}
\int \limits_{T_{i,j,k}} \omega_{k_1}(x) \omega_{k_2}(y) \omega_{k_3}(z) \left( \nabla \cdot \mathbf{F}_{d}\right) d\x = 
\int \limits_{\partial T_{i,j,k}} \omega_{k_1}(x)\omega_{k_2}(y)\omega_{k_3}(z) \mathbf{F}_{d}\cdot\vec{n} - 
\int \limits_{T_{i,j,k}} \nabla  \left(\omega_{k_1}(x)\omega_{k_2}(y)\omega_{k_3}(z)\right) \cdot \mathbf{F}_{d}  d\x,
\end{align}
and then
\begin{align}
&\Mass{xyz} \cdot\left( \dof{U}^{n+\halb}_{i,j,k}- \dof{Fu}^{n}_{i,j,k}  \right) -  \nu \frac{\Delta t}{\Delta x} \Mass{yz}\cdot\left(\barr{R}{x} \cdot\dof{U}^{\textsl{(x)}}_{i+\frac{1}{2},j,k} - \barr{L}{x} \cdot\dof{U}^{\textsl{(x)}}_{i-\frac{1}{2},j,k}\right) +\nonumber\\
&- \nu \frac{\Delta t}{\Delta y}\Mass{zx}\cdot\left(\barr{R}{y} \cdot\dof{U}^{\textsl{(y)}}_{i,j+\frac{1}{2},k} - \barr{L}{y}\cdot \dof{U}^{\textsl{(y)}}_{i,j-\frac{1}{2},k}\right)    - 
\nu \frac{\Delta t}{\Delta z} \Mass{xy}\cdot\left(\barr{R}{z} \cdot\dof{U}^{\textsl{(z)}}_{i,j,k+\frac{1}{2}} - \barr{L}{z} \cdot\dof{U}^{\textsl{(z)}}_{i,j,k-\frac{1}{2}}\right) = 0,\label{eq:uSDG12}  %\\
%&\Mass{xyz}\cdot \left( \dof{V}^{n+1}_{i,j,k}- \dof{Fv}^{n}_{i,j,k}  \right) - \frac{\Delta t}{\Delta x} \Mass{yz}\cdot\left(\barr{R}{x} \cdot\dof{V}^{\textsl{(x)}\,n+1}_{i+\frac{1}{2},j,k} - \barr{L}{x} \cdot\dof{V}^{\textsl{(x)}\,n+1}_{i-\frac{1}{2},j,k}\right) +\nonumber\\
%&- \frac{\Delta t}{\Delta y}\Mass{zx}\cdot\left(\barr{R}{y} \cdot\dof{V}^{\textsl{(y)}\,n+1}_{i,j+\frac{1}{2},k} - \barr{L}{y}\cdot \dof{V}^{\textsl{(y)}\,n+1}_{i,j-\frac{1}{2},k}\right)    - \frac{\Delta t}{\Delta z} \Mass{xy}\cdot\left(\barr{R}{z} \cdot\dof{V}^{\textsl{(z)}\,n+1}_{i,j,k+\frac{1}{2}} - \barr{L}{z} \cdot\dof{V}^{\textsl{(z)}\,n+1}_{i,j,k-\frac{1}{2}}\right) = 0,\label{eq:vSDG12} \\
%&\Mass{xyz} \cdot\left( \dof{W}^{n+1}_{i,j,}- \dof{Fw}^{n}_{i,j,k}  \right) - \frac{\Delta t}{\Delta x} \Mass{yz}\cdot\left(\barr{R}{x} \cdot\dof{W}^{\textsl{(x)}\,n+1}_{i+\frac{1}{2},j,k} - \barr{L}{x} \cdot\dof{W}^{\textsl{(x)}\,n+1}_{i-\frac{1}{2},j,k}\right) +\nonumber\\
%&- \frac{\Delta t}{\Delta y}\Mass{zx}\cdot\left(\barr{R}{y} \cdot\dof{W}^{\textsl{(y)}\,n+1}_{i,j+\frac{1}{2},k} - \barr{L}{y}\cdot \dof{W}^{\textsl{(y)}\,n+1}_{i,j-\frac{1}{2},k}\right)    - \frac{\Delta t}{\Delta z} \Mass{xy}\cdot\left(\barr{R}{z} \cdot\dof{W}^{\textsl{(z)}\,n+1}_{i,j,k+\frac{1}{2}} - \barr{L}{z} \cdot\dof{W}^{\textsl{(z)}\,n+1}_{i,j,k-\frac{1}{2}}\right) = 0,\label{eq:uSDG12} 
%& .\label{eq:pSDG1}
\end{align}
where the velocity derivatives can be computed with 
%\begin{align}
%\
%\end{align}
%that can be written exactly as %in the same notation of the pressure gradients (\ref{eq:weakgradient}) exactly as
\begin{align}
\dof{U}^{\textsl{(x)}}_{i+\frac{1}{2},j,k} = \frac{1}{\Delta x} \left(\Mass{x}\right)^{-1} \cdot \left( \tild{R}{x}\cdot \dof{U}^{n+\halb}_{i+1,j,k} -\tild{L}{x}\cdot \dof{U}^{n+\halb}_{i,j,k}  \right),\label{eq:Ux}\\
\dof{U}^{\textsl{(y)}}_{i,j+\frac{1}{2},k} = \frac{1}{\Delta y} \left(\Mass{y}\right)^{-1} \cdot \left( \tild{R}{y}\cdot \dof{U}^{n+\halb}_{i+1,j,k} -\tild{L}{y}\cdot \dof{U}^{n+\halb}_{i,j,k}  \right),\label{eq:Uy}\\
\dof{U}^{\textsl{(z)}}_{i,j,k+\frac{1}{2}} = \frac{1}{\Delta z} \left(\Mass{z}\right)^{-1} \cdot \left( \tild{R}{z}\cdot \dof{U}^{n+\halb}_{i+1,j,k} -\tild{L}{z}\cdot \dof{U}^{n+\halb}_{i,j,k}  \right),\label{eq:Uz}
\end{align}
and analogous equations for $\dof{V}^{n+\halb}$ and $\dof{W}^{n+\halb}$. By a formal substitution of equations (\ref{eq:Ux}-\ref{eq:Uz}) into (\ref{eq:uSDG12}), the following systems for the velocity components can be written in a very compact form as 
\begin{align}
\left(\mathbb{M} + \nu  \mathbb{H}\right) \cdot \mathpzc{U}^{n+\halb} = \mathbb{M}  \cdot \mathpzc{Fu}^n, \;\;\;
\left(\mathbb{M} + \nu  \mathbb{H}\right) \cdot \mathpzc{V}^{n+\halb} = \mathbb{M}  \cdot \mathpzc{Fv}^n,\;\;\; %\mathpzc{b}_v \label{eq:implicitviscosity} \\
\left(\mathbb{M} + \nu  \mathbb{H}\right) \cdot \mathpzc{W}^{n+\halb} = \mathbb{M}  \cdot \mathpzc{Fw}^n,\label{eq:implicitviscosity} %\mathpzc{b}_w\nonumber
\end{align}
by means of the piecewise polynomials $\mathpzc{U}^{n+\halb}$, $\mathpzc{V}^{n+\halb}$, $\mathpzc{W}^{n+\halb}$ on the \emph{main grid}. Here, $\mathbb{M} > 0$ is the diagonal tensor of the 
\emph{element mass matrices}, $\mathbb{H} \geq 0$ is \emph{exactly the same} operator that has been obtained for the pressure system in (\ref{eq:pSyst1d}). Then, the tensor coefficient matrices of systems 
(\ref{eq:implicitviscosity}) are \emph{all} \emph{\textbf{positive definite}} because the sum of a positive semi-definite matrix and a positive definite matrix is positive definite. 
The right hand sides of the system of equations (\ref{eq:implicitviscosity}) contains only the fully explicit discretization of the \emph{nonlinear convective terms} 
$\mathpzc{Fu}^n=\left\{\dof{Fu}^n\right\}$,  $\mathpzc{Fv}^n=\left\{\dof{Fv}^n\right\}$ and $\mathpzc{Fw}^n=\left\{\dof{Fw}^n\right\}$ multiplied by the mass matrix.  
This means that the CFL-type restriction on the time-step (\ref{eq:CFL}) looses the dependency on the viscosity and relaxes to 
\begin{align}
\Delta t = \frac{\text{CFL}}{\left(2N+1\right)}\left(\frac{|u_\text{max}|}{\Delta x} + \frac{|v_\text{max}|}{\Delta y} + \frac{|w_\text{max}|}{\Delta z}\right)^{-1}.\label{eq:CFL1}
\end{align}
If the solutions % $\mathpzc{U}\equiv\left\{\dof{U}\right\}$, $\mathpzc{V}\equiv\left\{\dof{V}\right\}$ and $\mathpzc{W}\equiv\left\{\dof{W}\right\}$  
of the semi-implicit formulation of the advection-diffusion  system (\ref{eq:implicitviscosity}) substitute the fully explicit terms  $\dof{Fu}^n$, $\dof{Fv}^n$ and $\dof{Fw}^n$ in (\ref{eq:Fu}-\ref{eq:Fw}), a coherent DG scheme is obtained by means of a \emph{fractional time-stepping} approach. % that is first order in time. 
%following coupled system would occur
%\begin{align}
%&\left(\mathbb{M} + \nu  \mathbb{H}\right) \cdot \mathpzc{U}^{n+1} = \mathbb{M}\mathpzc{F}_u^n - \frac{\Delta t }{\Delta x} \mathbb{M}^{yz}\mathpzc{D}_{\textsl{w}}^x \mathpzc{P}^{n+1}\nonumber\\
%&\left(\mathbb{M} + \nu  \mathbb{H}\right) \cdot \mathpzc{V}^{n+1}  = \mathbb{M}\mathpzc{F}_v^n - \frac{\Delta t }{\Delta y}\mathbb{M}^{zx} \mathpzc{D}_{\textsl{w}}^y \mathpzc{P}^{n+1} 
%\label{eq:implicitviscosity2}\\
%&\left(\mathbb{M} + \nu  \mathbb{H}\right) \cdot \mathpzc{W}^{n+1}  = \mathbb{M}\mathpzc{F}_w^n - \frac{\Delta t }{\Delta z}\mathbb{M}^{xy} \mathpzc{D}_{\textsl{w}}^z \mathpzc{P}^{n+1}\nonumber
%\end{align}
%with $\mathpzc{D}_{\textsl{w}}^{\xi}$ is the \emph{bi}-block diagonal tensor matrix (\ref{eq:weakgradientop}) that represents the weak formulation of the gradient in the $\xi$-direction. In principle a similar resolution can be performed after multiplication of the inverse tensor matrix $\left(\mathbb{M} + \nu  \mathbb{H}\right)^{-1}$ and substitution in the incompressibility condition (\ref{}), i.e.
%\begin{align}
%\mathbb{M}^{yz} \left(\mathpzc{D}_{\textsl{w}}^x\right)^T \mathpzc{U} + 
%\mathbb{M}^{zx} \left(\mathpzc{D}_{\textsl{w}}^y\right)^T \mathpzc{V} + 
%\mathbb{M}^{xy} \left(\mathpzc{D}_{\textsl{w}}^z\right)^T \mathpzc{W} = 0
%\end{align}
 %but the resulting algorithm could become expensive from a computational point of view even if theoretically \emph{\textbf{second order accurate in time}}.
The resulting numerical scheme can be finally written in compact form as
\begin{align}
&\left(\mathbb{M} + \nu  \mathbb{H}\right) \cdot \mathpzc{U}^{n+\frac{1}{2}}_{\#} =  \mathbb{M}  \cdot\mathpzc{Fu}^n_{\#} \label{eq:LOdiffusionU}, \\
&\left(\mathbb{M} + \nu  \mathbb{H}\right) \cdot \mathpzc{V}^{n+\frac{1}{2}}_{\#} =  \mathbb{M}  \cdot\mathpzc{Fv}^n_{\#} \label{eq:LOdiffusionV}, \\
&\left(\mathbb{M} + \nu  \mathbb{H}\right) \cdot \mathpzc{W}^{n+\frac{1}{2}}_{\#} =  \mathbb{M}  \cdot\mathpzc{Fw}^n_{\#} \label{eq:LOdiffusionW}, 
\end{align}
\begin{gather}
\mathbb{H}\cdot \mathpzc{P}^{n+\theta}_{\#}   = \mathpzc{b}^{n+\frac{1}{2}}_{\#}, \label{eq:LOpressure}\\
\text{with} \;\;\; \mathpzc{b}^{n+\frac{1}{2}}_{\#}\equiv\mathbb{M}^{yz} \left(\mathcal{D}^{x}\right)^T  \cdot \mathpzc{U}^{n+\frac{1}{2}}_{*}    +
\mathbb{M}^{zx} \left(\mathcal{D}^{y}\right)^T \cdot \mathpzc{V}^{n+\frac{1}{2}}_{*}    + \mathbb{M}^{xy} \left(\mathcal{D}^{z}\right)^T \cdot \mathpzc{W}^{n+\frac{1}{2}}_{*},  \nonumber
\end{gather}
\begin{align}
\mathpzc{U}^{n+1}_{*}& =  \mathpzc{U}^{n+\frac{1}{2}}_{*}   - \frac{\Delta t}{\Delta x}  \left(\mathbb{M}^{-1}  \mathcal{D}\right)^x \cdot\mathpzc{P}^{n+\theta}_{{\#}}, \label{eq:LOvelocityU}\\
\mathpzc{V}^{n+1}_{*} &= \mathpzc{V}^{n+\frac{1}{2}}_{*} - \frac{\Delta t}{\Delta y}  \left(\mathbb{M}^{-1}  \mathcal{D}\right)^y \cdot\mathpzc{P}^{n+\theta}_{{\#}},   \label{eq:LOvelocityV} \\
\mathpzc{W}^{n+1}_{*}& = \mathpzc{W}^{n+\frac{1}{2}}_{*}  - \frac{\Delta t}{\Delta z} \left(\mathbb{M}^{-1}  \mathcal{D}\right)^z \cdot\mathpzc{P}^{n+\theta}_{{\#}},   \label{eq:LOvelocityW} %\\
%&\left(
%\begin{array}{ccc} 
%\mathbb{M}^{yz} \mathcal{D}^x &
%\mathbb{M}^{zx} \mathcal{D}^y &
%\mathbb{M}^{xy} \mathcal{D}^z 
%\end{array}
%\right) \cdot
%\left(
%\begin{array}{ccc} 
%\mathpzc{U}^{n+1} \\
%\mathpzc{V}^{n+1} \\
%\mathpzc{W}^{n+1} 
%\end{array}
%\right) = 0 \nonumber\\
\end{align}
where $\mathpzc{Fu}^n  = \left\{\dof{Fu}^n \right\}$, $\mathpzc{Fv}^n  = \left\{\dof{Fv}^n  \right\}$ and $\mathpzc{Fw}^n  = \left\{\dof{Fw}^n  \right\}$ are the purely explicit 
discretization of the nonlinear convective terms outlined in the previous section; '$*$' is used for the field variables defined along the dual grids; '$\#$' for those variables 
defined along the main grid; notice that the projection of a field variable from the dual grid to the main grid (\reds{$\pi_{:*\rightarrow \# }$}) and vice versa  (\reds{$\pi_{:\#\rightarrow *}$}) is simply performed by the 
$L_2$ projections defined in \eqref{eqn.l2forward} and \eqref{eqn.l2backward}, respectively.  
Once the nonlinear convective terms have been computed with respect to the field variables of the old time step $t^n$, then the viscous terms are computed implicitly at 
a \emph{fictitious} fractional time-step $t^{n+\frac{1}{2}}$ (\ref{eq:LOdiffusionU}-\ref{eq:LOdiffusionW}). Finally the pressure forces and the incompressibility condition are 
solved implicitly (\ref{eq:LOpressure}) and the field variables at the future time $t^{n+1}$ are worked out (\ref{eq:LOvelocityU}-\ref{eq:LOvelocityW}). The fractional time 
discretization $t^{n+\frac{1}{2}}$ is only an auxiliary notation emphasizing that it is a intermediate stage. In fact, the real time evolution of the discrete diffusion 
equations (\ref{eq:LOdiffusionU}-\ref{eq:LOdiffusionW}) is from $t^n$ to $t^{n+1}=t^n+\Delta t$. 

\subsection{Numerical validation}

In order to check the ability of the new method in solving the governing equations accurately, some different numerical test problems in two and three space dimensions have been chosen, 
for which an analytical or other numerical reference solutions exist.

\subsubsection{Blasius boundary layer}

\begin{figure} 
\centering %subfloat
			\includegraphics[width=0.45\textwidth]{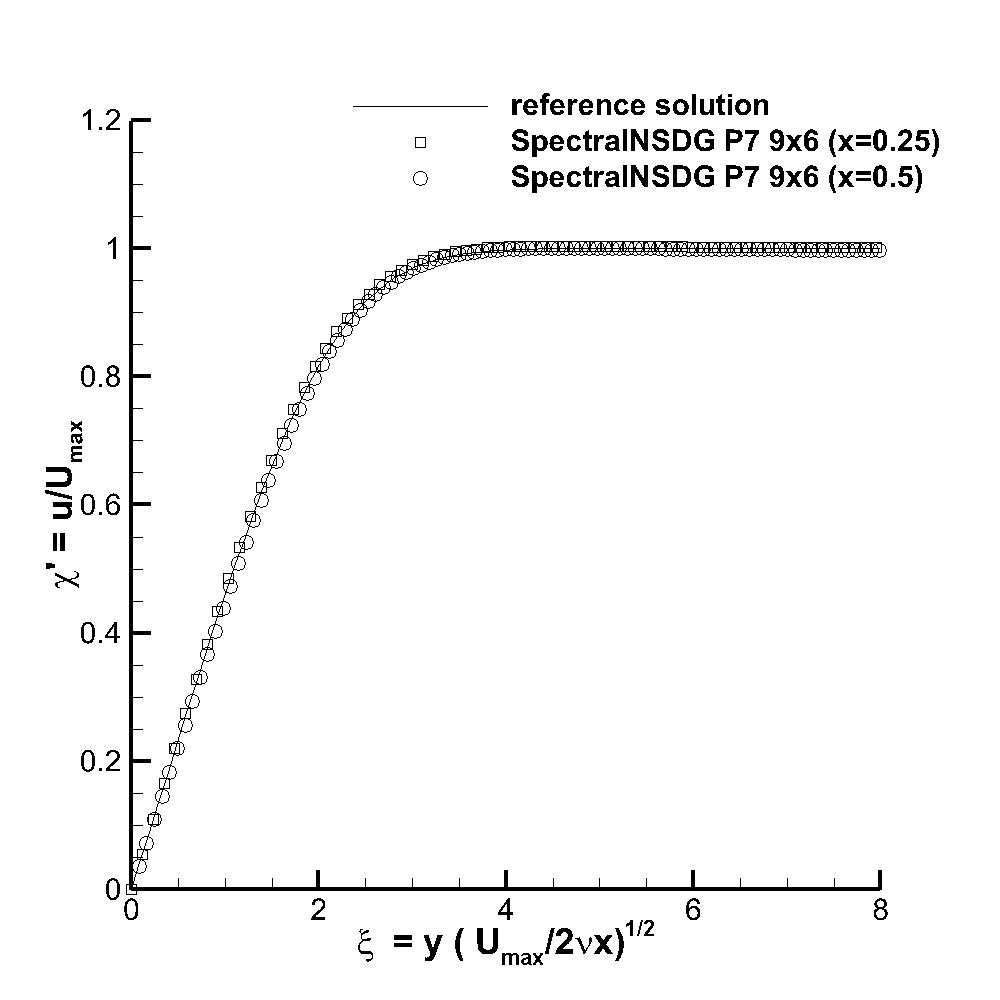}\\
			\includegraphics[width=0.9\textwidth]{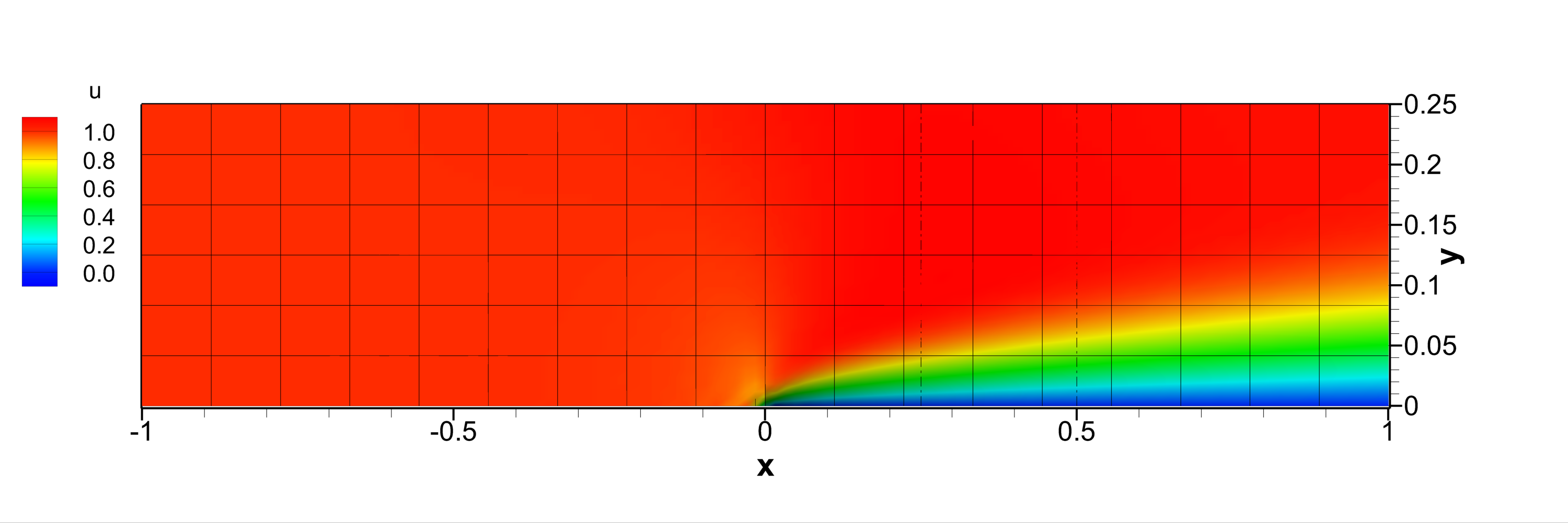}
\caption{Top: Comparison of the reference solution of Blasius with the numerical results at $t=5$ for different axial positions $x=0.25$ and $x=0.50$ obtained with a staggered semi-implicit spectral 
DG-$\p_7$ scheme on a very coarse grid of $18\times6$ elements. Bottom: numerical solution for the horizontal velocity field computed at time $t=5$; the high-order elements of the main grid are 
depicted with solid lines; the vertical cuts at $x=0.25$ and $x=0.50$ with dash-dotted lines.}\label{fig:Blasius}
\end{figure}

In this test, a steady laminar boundary layer over a flat plate is considered. According to the theory of Prandtl \cite{Prandtl1904,BLTheory}, convective terms are of the order $1$ 
in the boundary layer along the horizontal direction, whereas the vertical accelerations are of the order of the boundary layer thickness.  
The spatial domain under consideration is $\Omega = [-1,1]\times[0,0.25]$ and the chosen kinematic viscosity is $\nu=10^{-3}$.  The flat-plate boundary is imposed at $y=0$ for $x>0$. 
Constant velocity $\mathbf{v}=(1,0)$ is imposed at the left inflow boundary, constant pressure $p=0$ at the right outflow, no-slip boundary conditions along the wall and no-jump condition 
in the rest. Results are shown in Fig. \ref{fig:Blasius}, obtained with our SIDG-$\p_7$ method using $\theta=1$ and a very coarse grid of only $18\times6$ elements. A very good 
agreement between the numerical solution obtained with the semi-implicit spectral DG scheme and the Blasius reference solution can be observed. Notice that the complete boundary 
layer is well resolved inside a single element close to $x=0$.

\subsubsection{Lid-driven cavity: 2D}
\begin{figure} 
\centering %subfloat
			\includegraphics[width=0.35\textwidth]{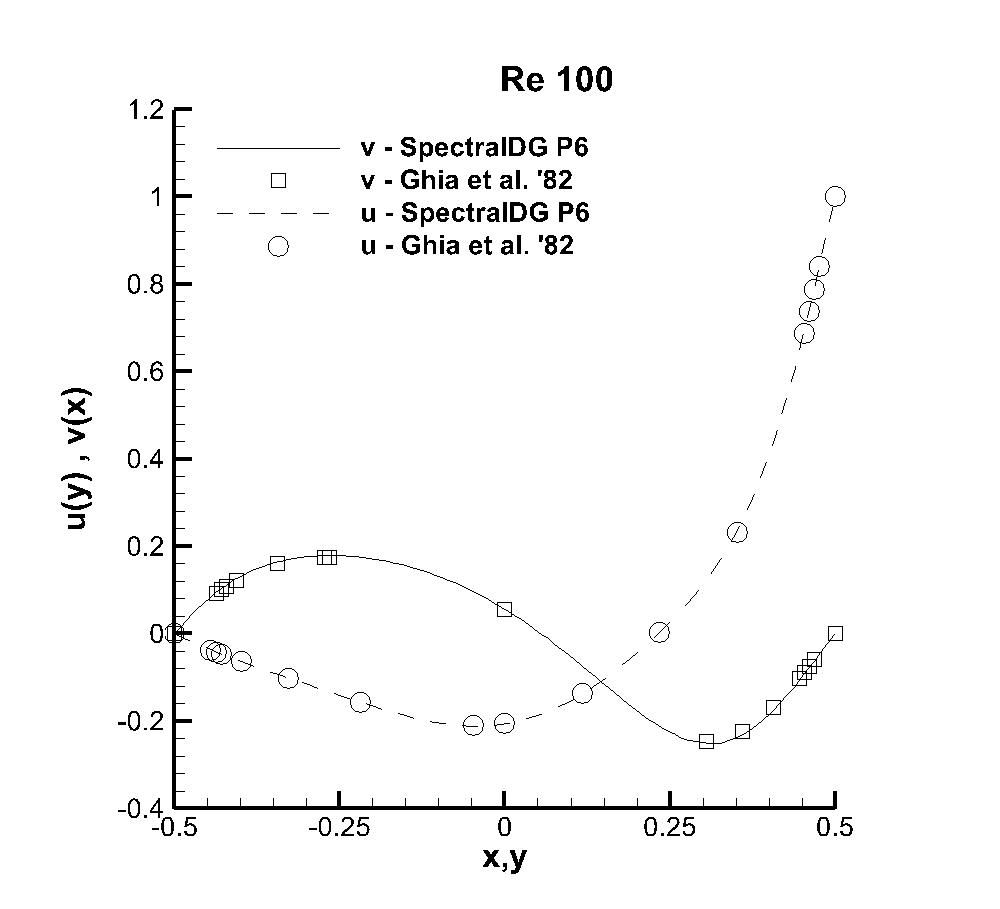}\;\includegraphics[width=0.35\textwidth]{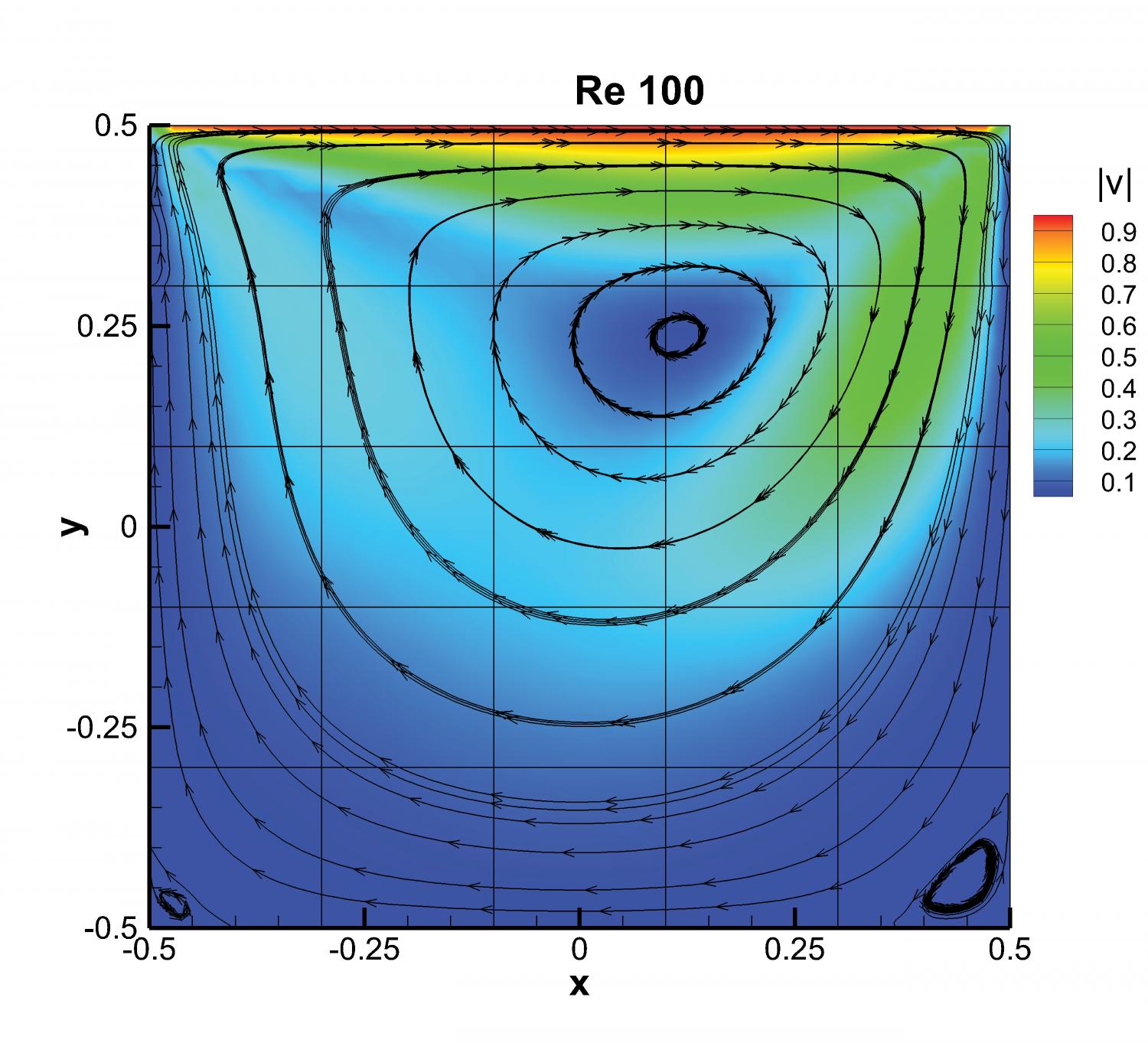}\\
			\includegraphics[width=0.35\textwidth]{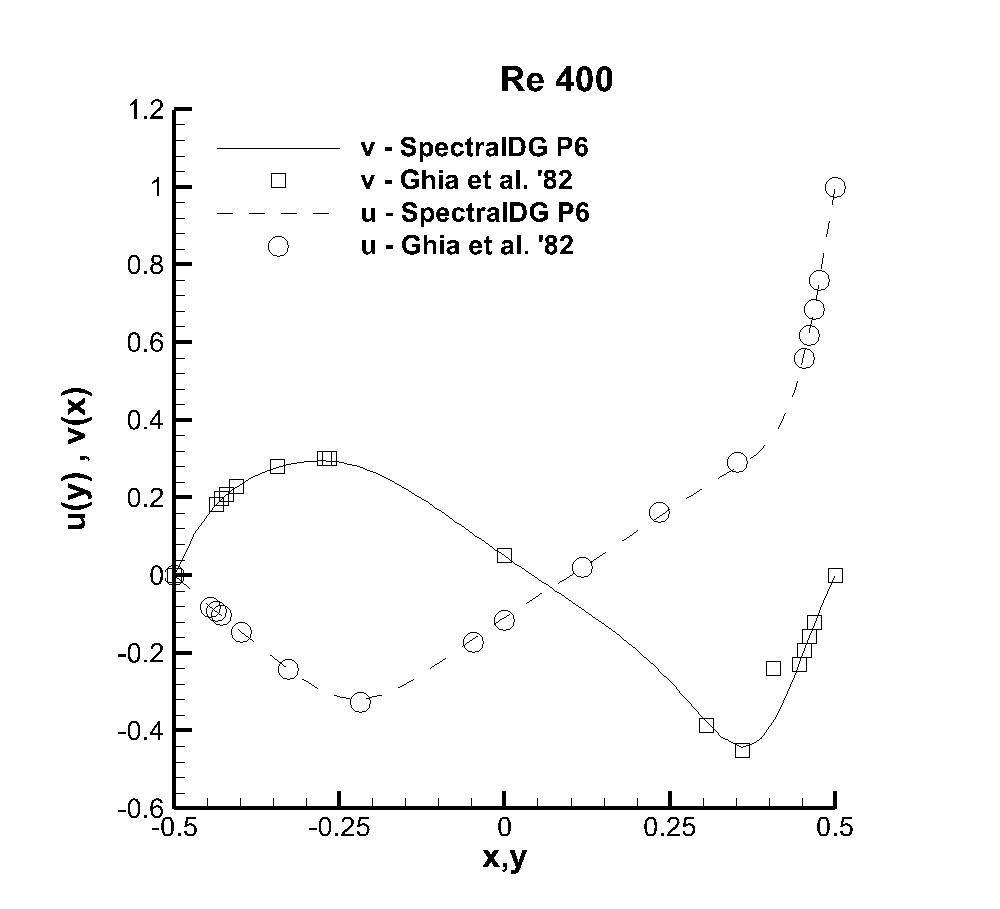}\;\includegraphics[width=0.35\textwidth]{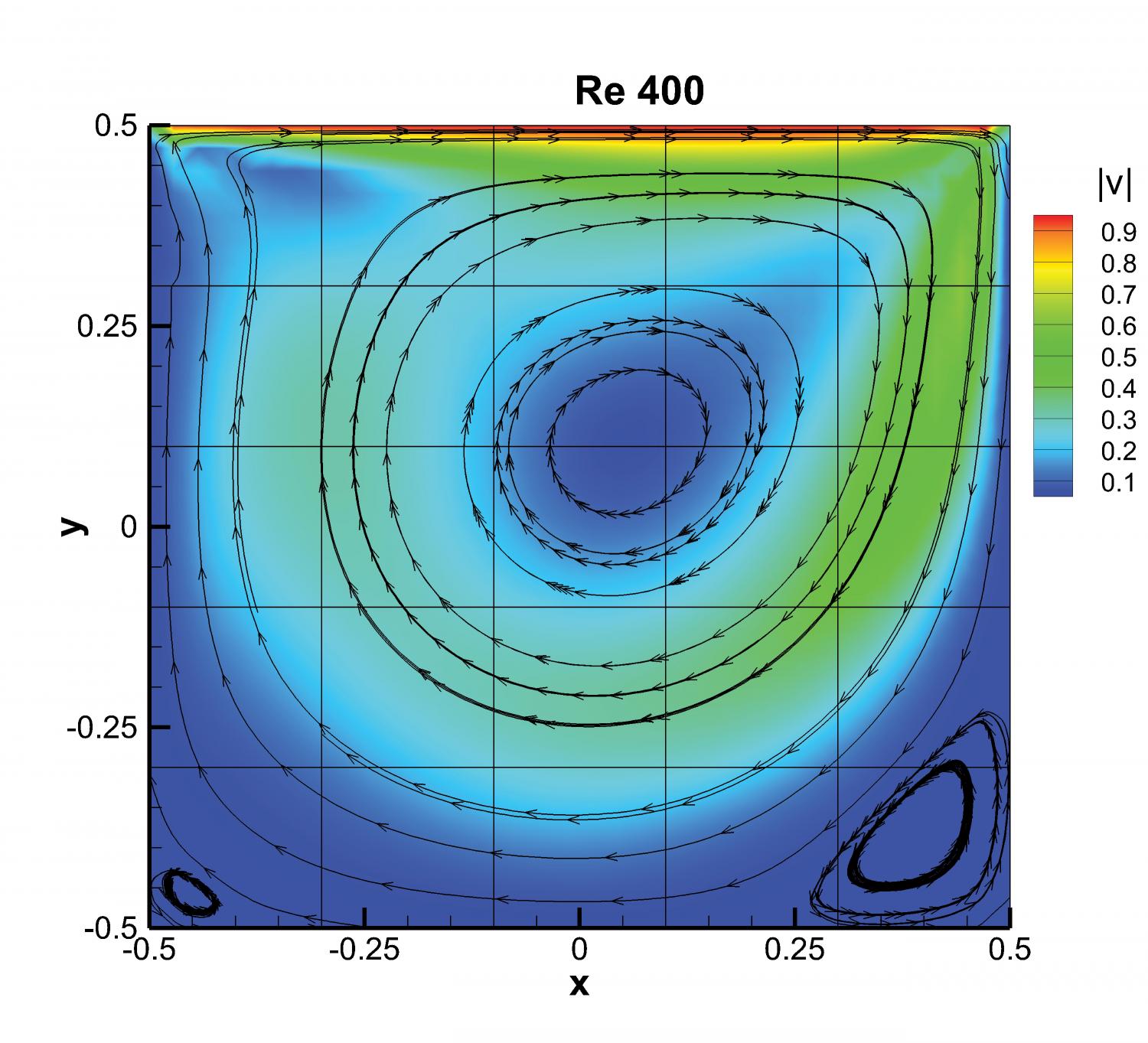}\\
			\includegraphics[width=0.35\textwidth]{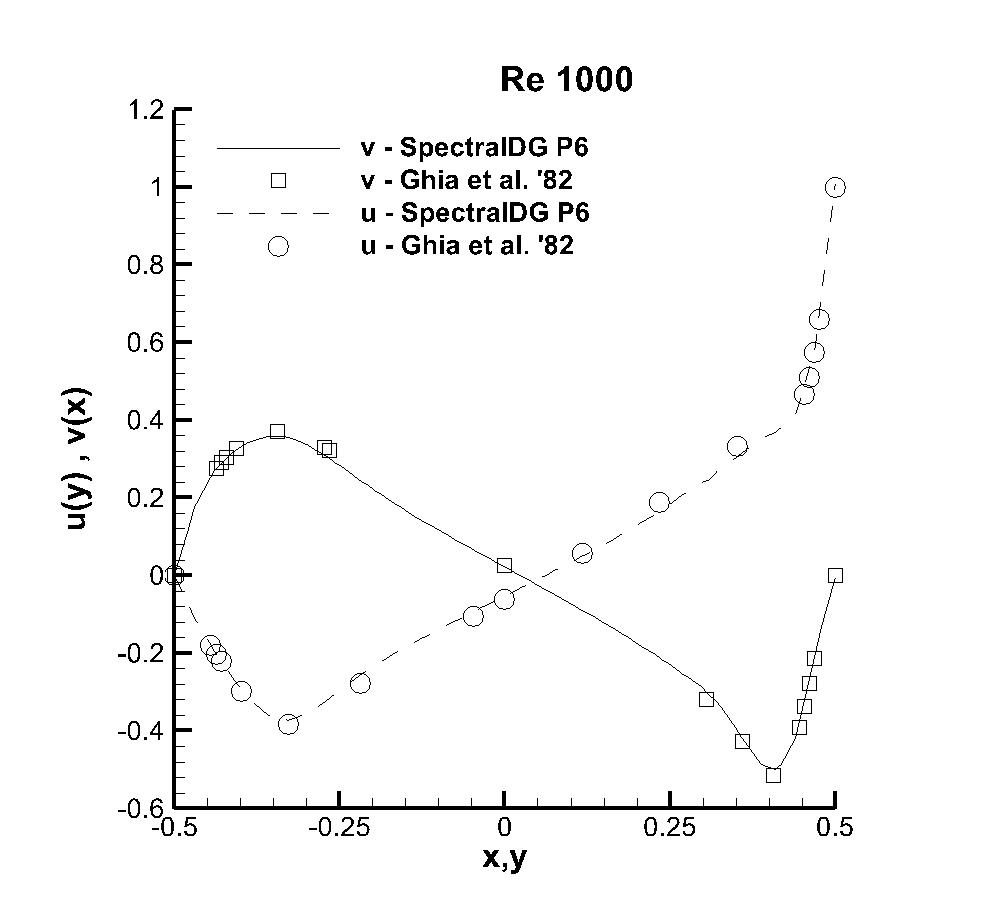}\;\includegraphics[width=0.35\textwidth]{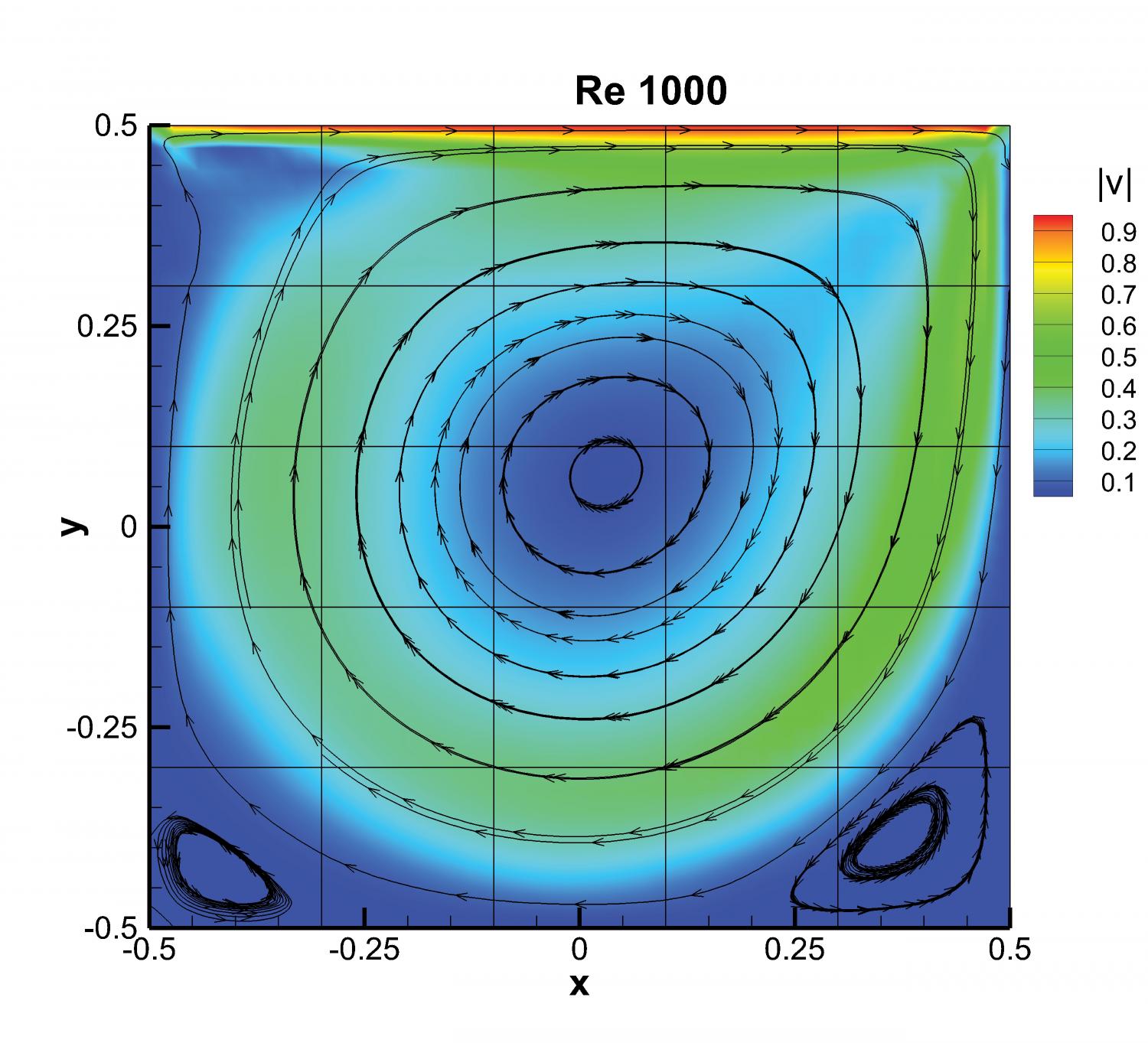}\\
			\includegraphics[width=0.35\textwidth]{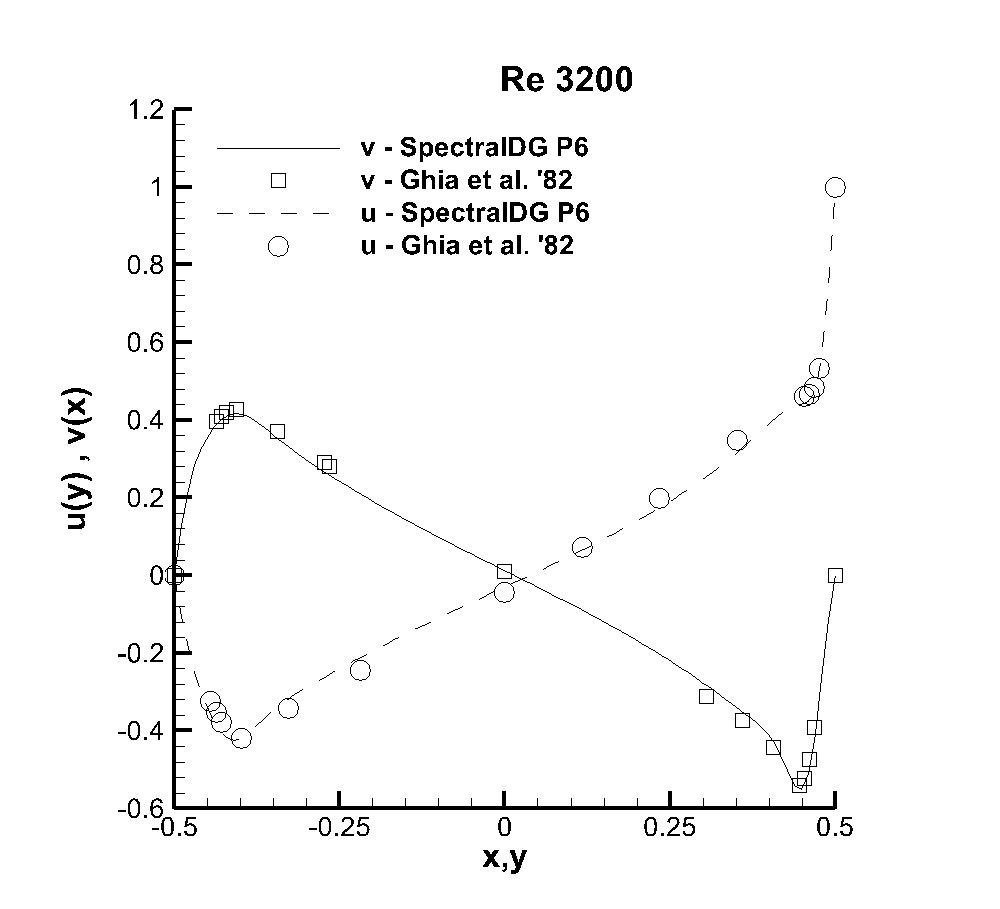}\;\includegraphics[width=0.35\textwidth]{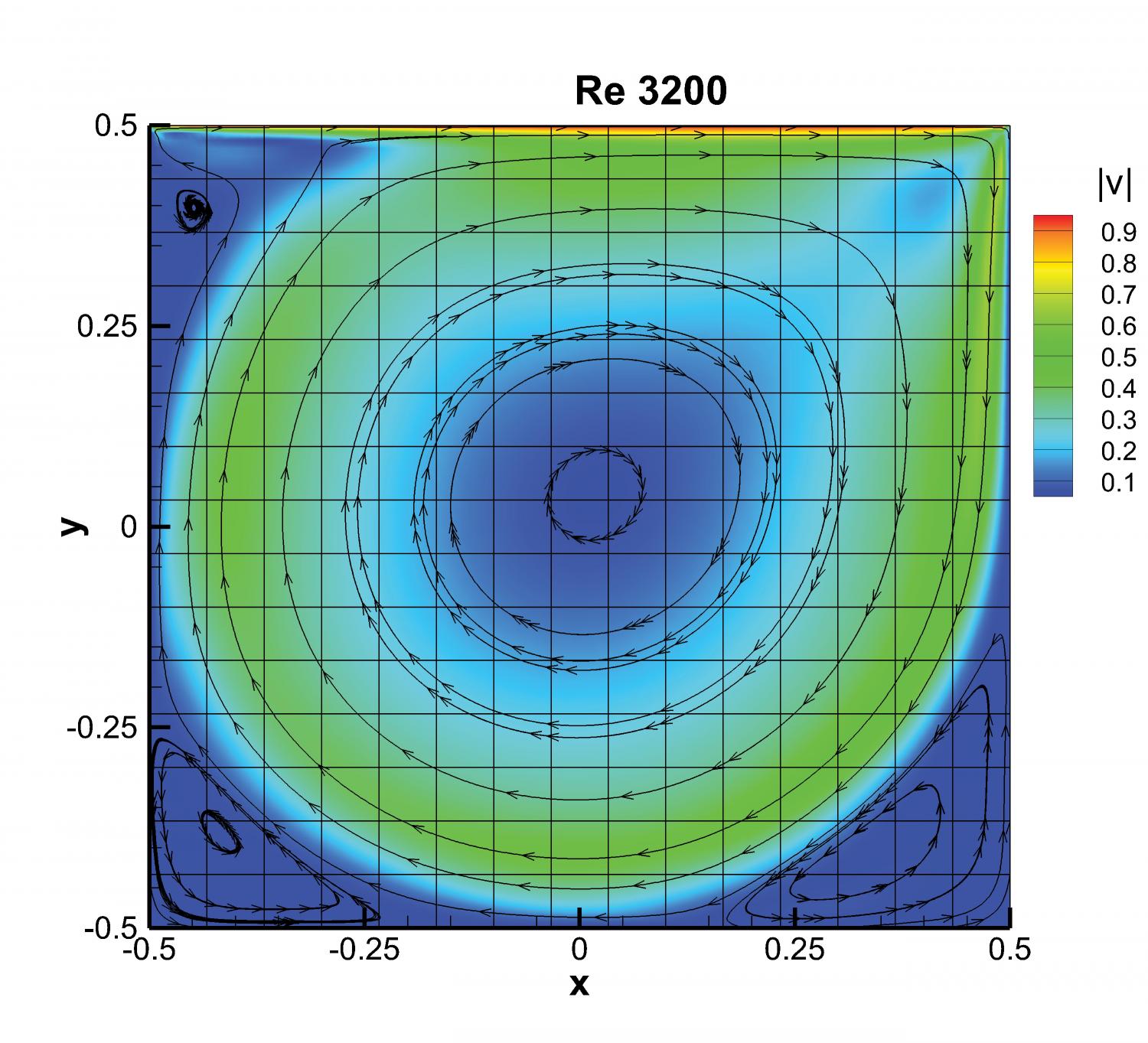}
\caption{The numerical solution obtained for the two dimensional lid-driven cavity problem compared with the numerical results of \cite{Ghia1982} at different Reynolds numbers, respectively, 
 from the top to the bottom:  Re=$100$, Re=$400$, Re=$1000$ using $5x5$ elements, and Re=$3200$ using $15x15$ elements obtained with a staggered semi-implicit spectral DG-$\p_{6}$ method.}\label{fig:LDCavity2D}
\end{figure}

An interesting standard benchmark problem for numerical methods applied to the incompressible Navier-Stokes equations is the lid-driven cavity, see \cite{Ghia1982}. 
In this test, a closed square cavity is filled with an incompressible fluid and the flow is driven by the upper wall that moves with velocity $\mathbf{v}=(1,0)$. 
The main difficulties in solving this problem arise from the singularities of the velocity gradient at the top right and at the top left corners, where the horizontal velocity component 
is a double valued function: $u=0$ at the left (or right) wall boundary and $u=1$ at the upper moving boundary. Moreover, the pressure is determined only up to a constant, because there 
are only velocity boundary conditions. 
The physical domain is $\Omega = [-0.5,0.5]\times[-0.5,0.5]$, the initial condition for velocity and pressure is set to $\mathbf{v}=(0,0)$ and $p=0$. % everywhere but $u=1$ at $y=0.5$.
Fig. \ref{fig:LDCavity2D} shows the computed results compared with the reference solution of Ghia et al. \cite{Ghia1982} next to the two-dimensional view of the velocity magnitude at 
different Reynolds numbers from Re=$100$ to Re=$3200$, obtained with the $\p_{6}$ version of our staggered semi-implicit spectral DG scheme. The implicitness factor has been chosen 
equal to $\theta=1$, since only a steady solution is sought for this test problem. 
Notice that the computed results match the reference solution very well, despite the presence of the corner singularities and the use of a very coarse mesh. 
A possibility to avoid the corner singularities in this test problem is the use the unified first order \textit{hyperbolic} formulation of viscous Newtonian fluids, 
recently proposed and used in \cite{PeshRom2014,HPRmodel}, which does not need the computation of velocity gradients in the numerical fluxes.

\subsubsection{Lid-driven cavity: 3D}
\begin{figure} 
\centering %subfloat
			\includegraphics[width=0.48\textwidth]{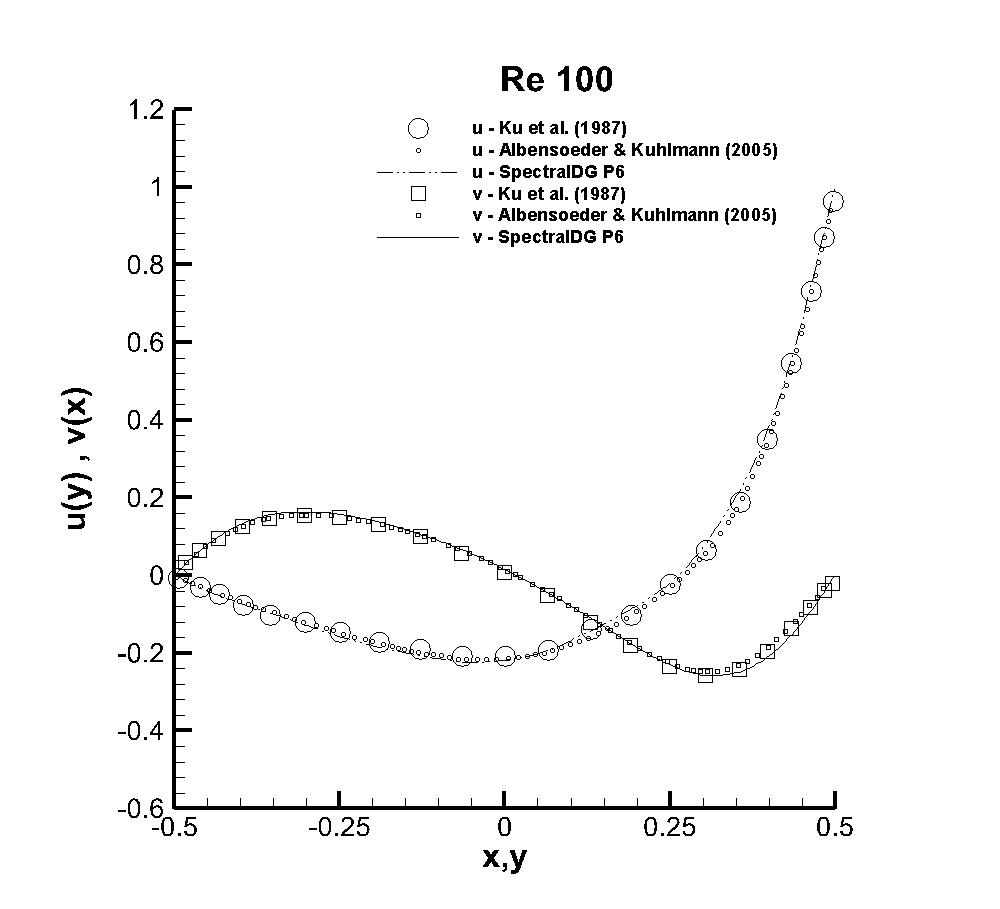}\;
			\includegraphics[width=0.48\textwidth]{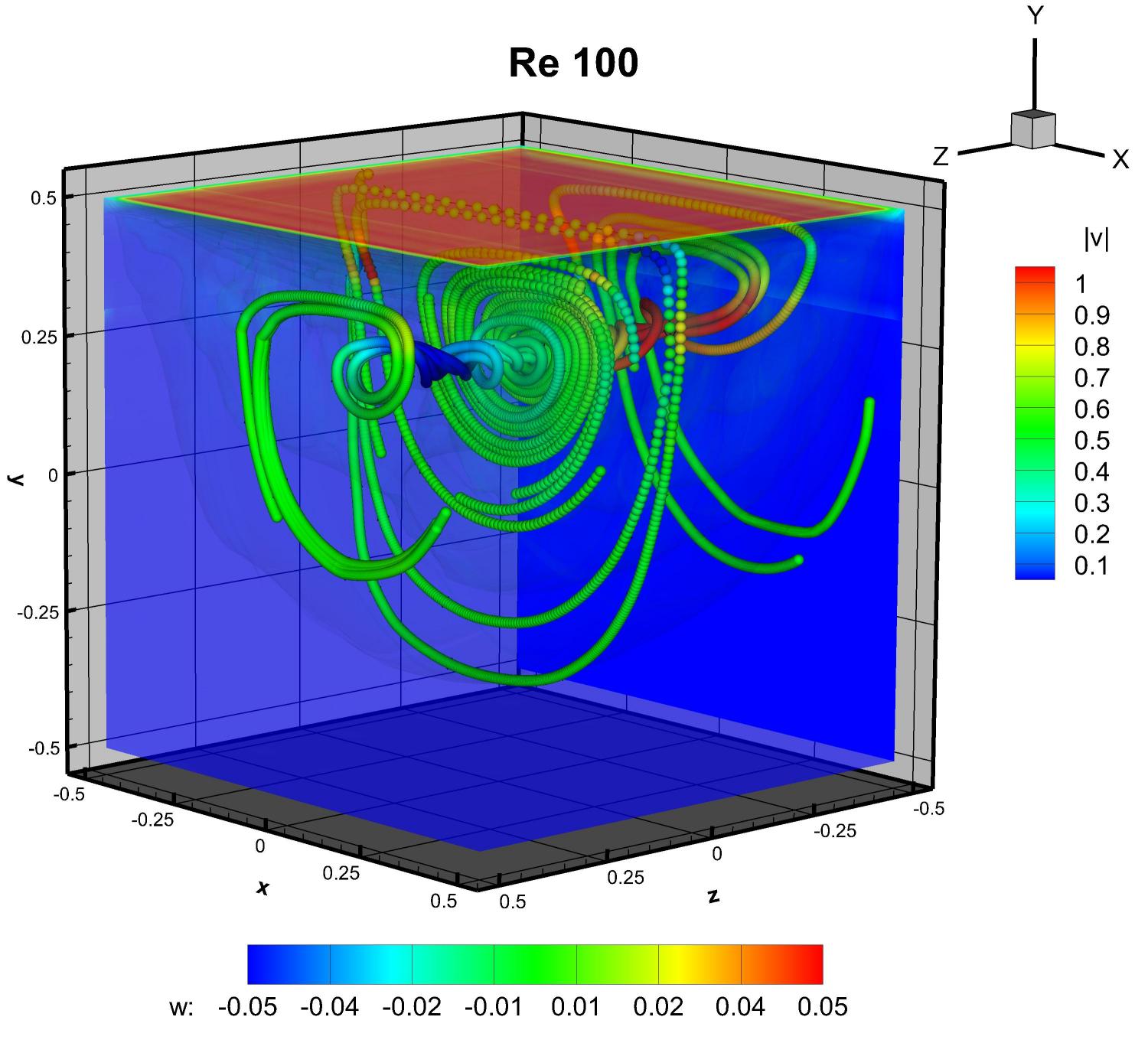}\\
			\includegraphics[width=0.48\textwidth]{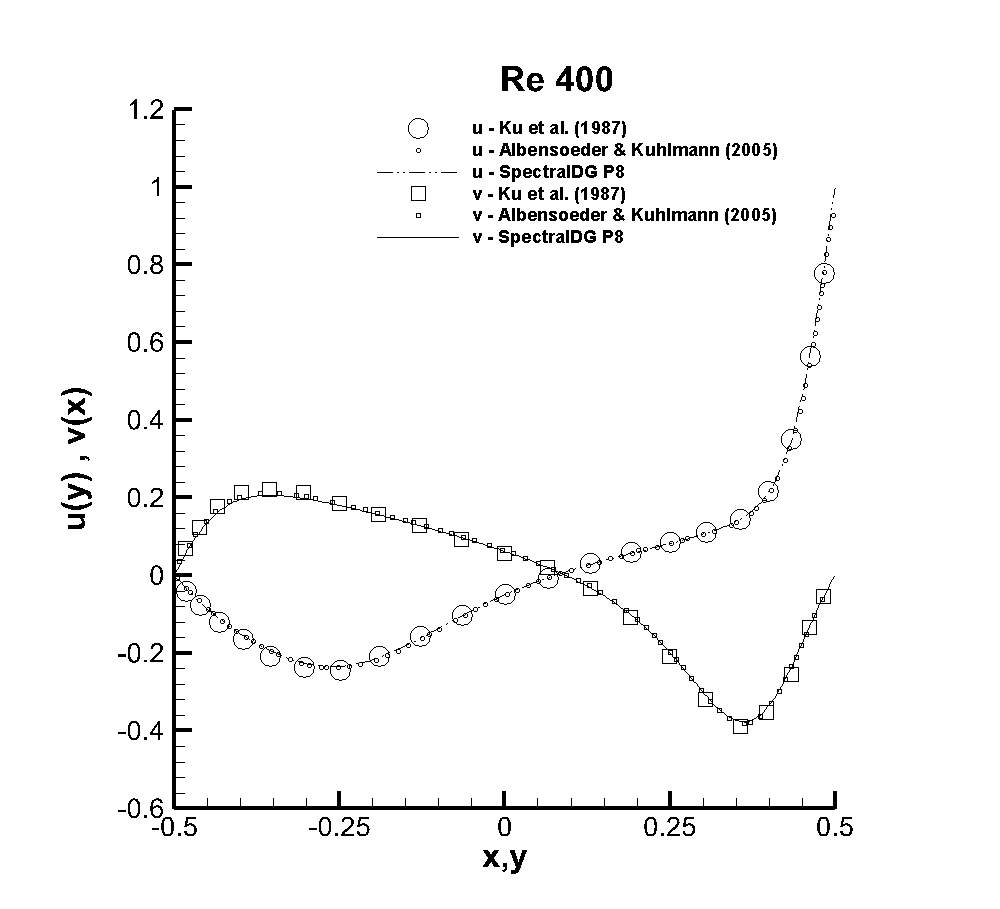}\;
			\includegraphics[width=0.48\textwidth]{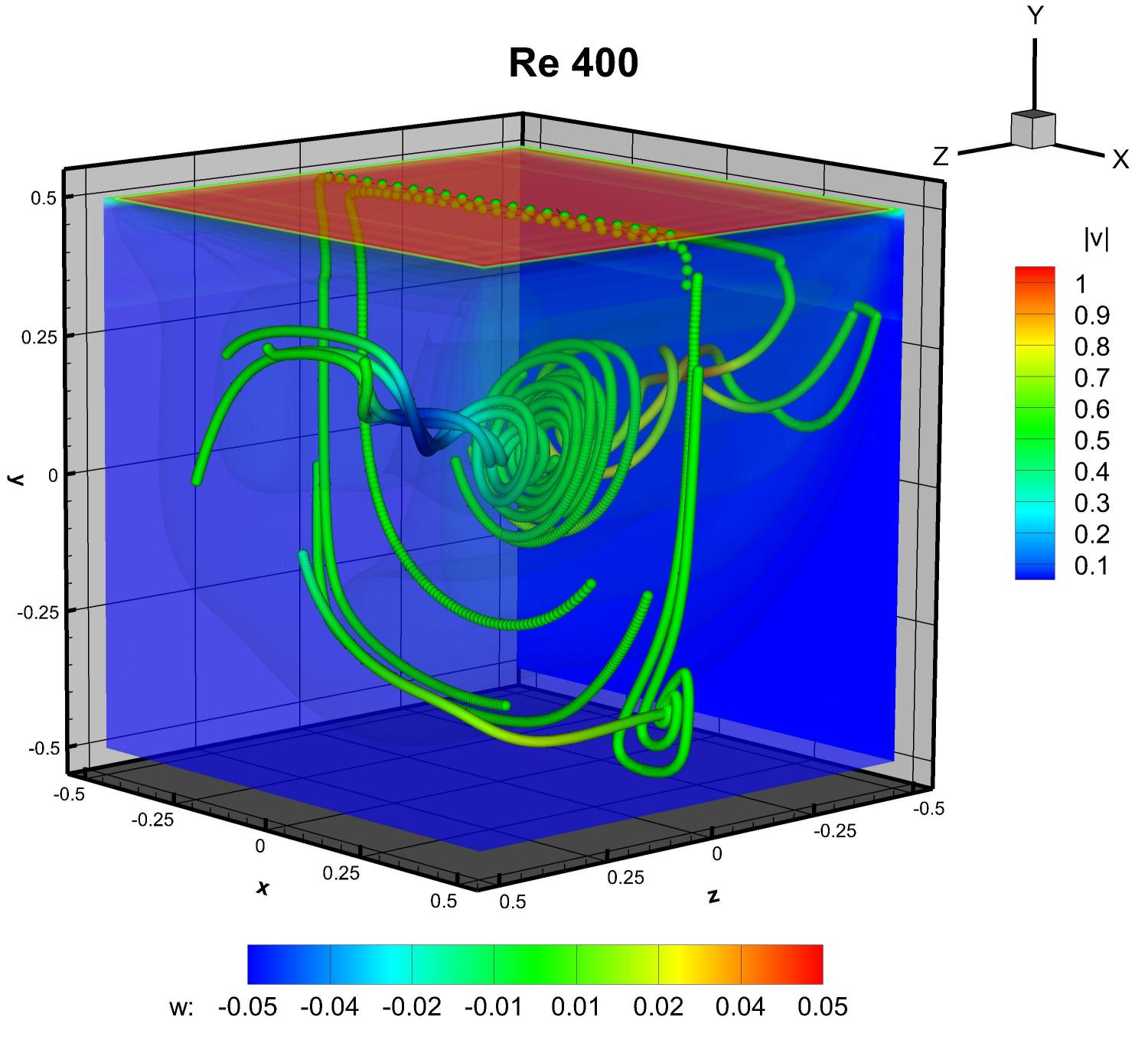}\\
			\includegraphics[width=0.32\textwidth]{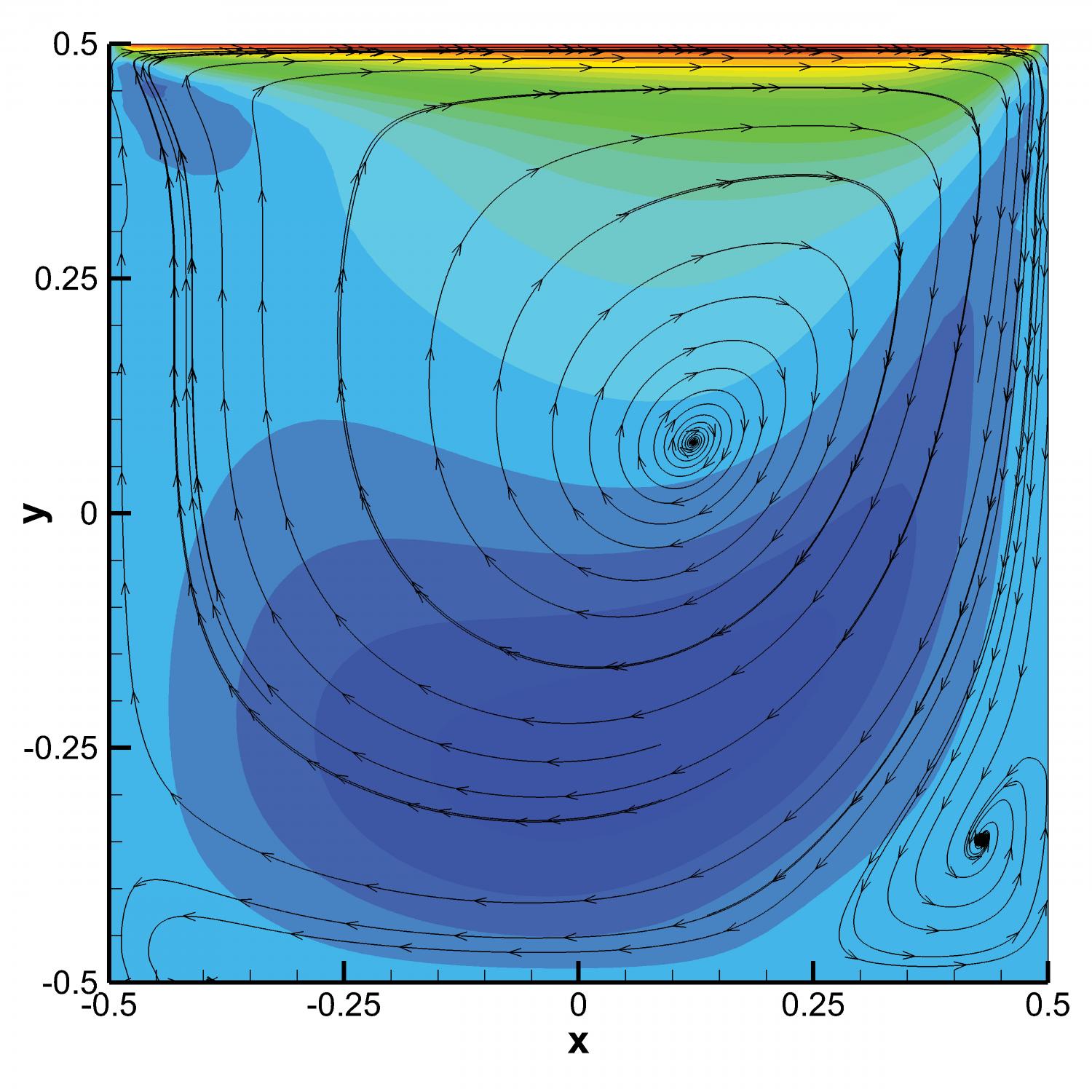}
			\includegraphics[width=0.32\textwidth]{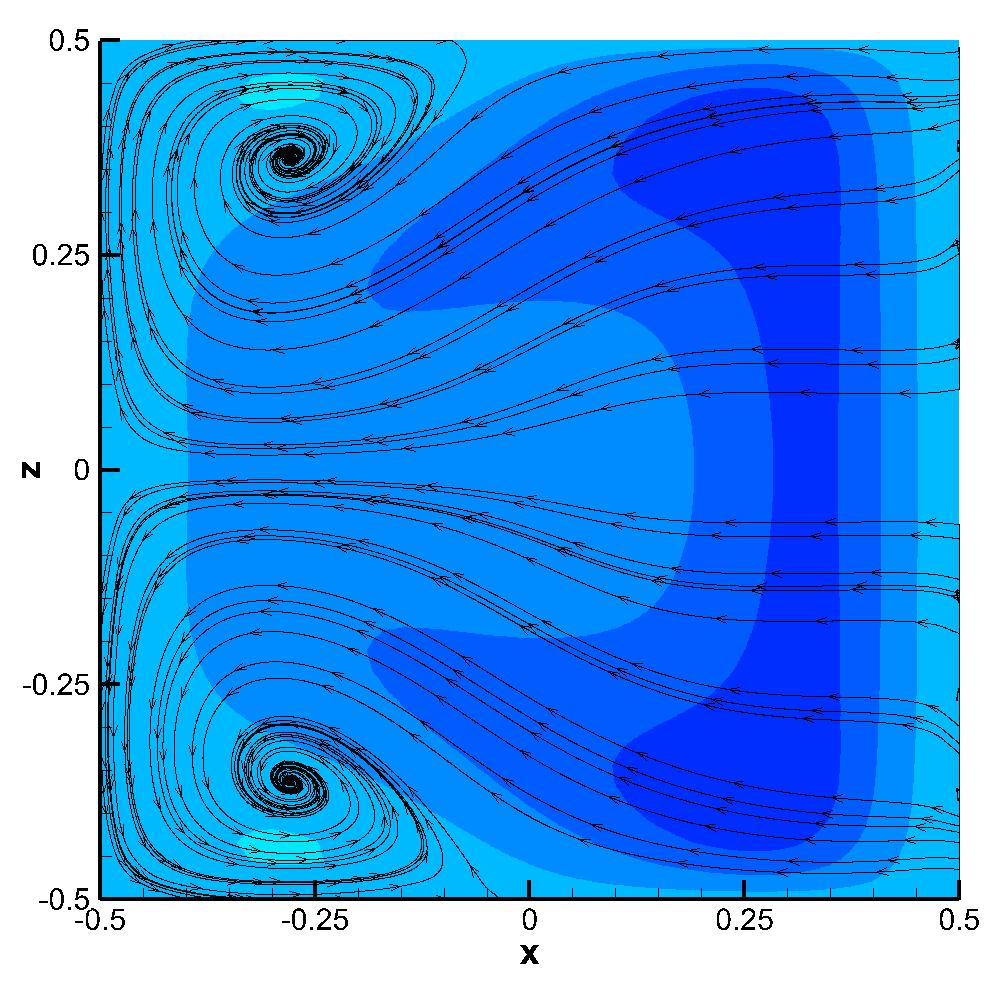}
			\includegraphics[width=0.32\textwidth]{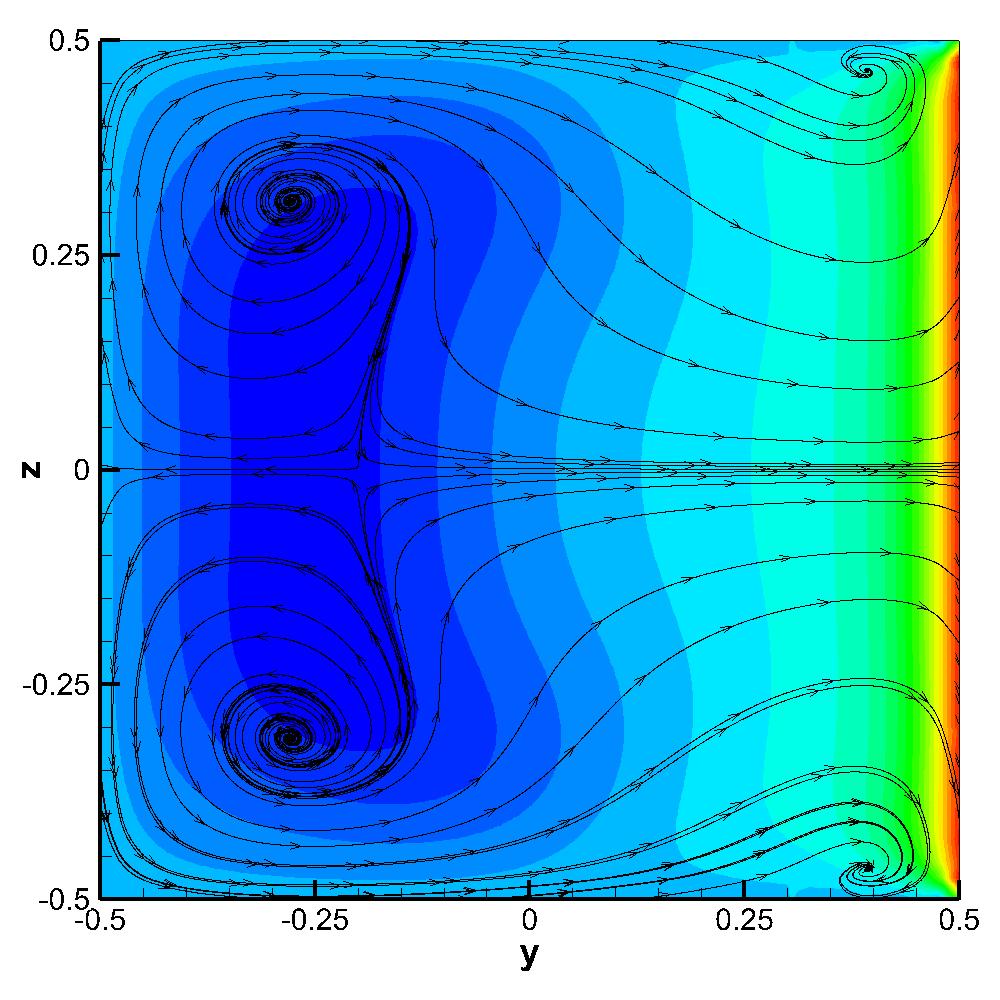}
\caption{The numerical solution obtained for the three dimensional lid-driven cavity problem compared with the numerical results of \cite{Albensoeder2005,Ku1987} at different Reynolds numbers, respectively, from the top to the center:  Re=$100$ and Re=$400$ using $5\times 5\times 5$ elements. The results have been obtained with a staggered semi-implicit spectral DG-$\p_{6}$ and DG-$\p_{8}$ 
method. The streamlines are colored with the $w$ velocity magnitude. The numerical solution for the case $Re=400$ has been interpolated along the three orthogonal planes $x-y$, $x-z$ and $y-z$ 
at the bottom: streamlines and the $u$ velocity are depicted.}\label{fig:LDCavity3D}
\end{figure}

In this section we present the three-dimensional version of the previous test case. A cubic cavity is filled by an incompressible fluid, and the upper wall boundary drives the fluid 
flow with a non-zero velocity $\mathbf{v}=(1,0,0)$. 
The presence of a third spatial dimension introduces a new degree of freedom to dynamics of the flow and the resulting flow field is different compared to the 2D case discussed before. 
The physical domain $\Omega = [-0.5, 0.5]^3$ has been divided into only $5 \times 5 \times 5$ spatial elements, with the implicitness factor $\theta = 1$ chosen for the time 
discretization. 
Fig. \ref{fig:LDCavity3D} shows the computed results compared with the reference data provided by \cite{Ku1987} and \cite{Albensoeder2005} next to the three-dimensional view of 
the flow field at Reynolds numbers Re=$100$ and Re=$400$, obtained with the $\p_{6}$ and  $\p_{8}$ version of our staggered semi-implicit spectral DG scheme.  Also for the three-dimensional 
cavity flow, our numerical results are in very good agreement with the reference data. At the bottom of Fig. \ref{fig:LDCavity3D} the numerical solution for the case $Re=400$ has been 
projected onto the three orthogonal planes ${x-y}$, ${x-z}$ and ${y-z}$. The expected secondary recirculations, which distinguish the three dimensional flow field from the two 
dimensional one, are clearly visible.

\subsubsection{Backward facing step: 2D}
\begin{figure} 
\centering %subfloat
			\includegraphics[width=0.95\textwidth]{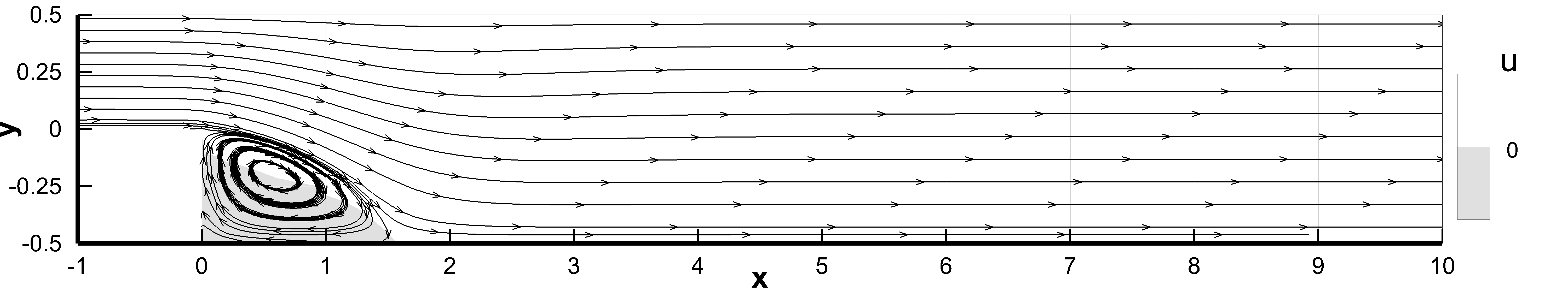}\\
			\includegraphics[width=0.95\textwidth]{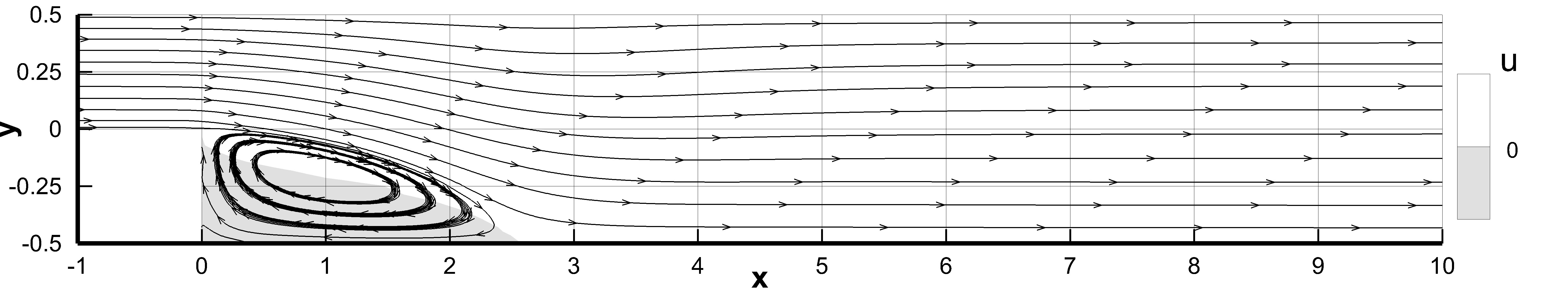}\\
			\includegraphics[width=0.95\textwidth]{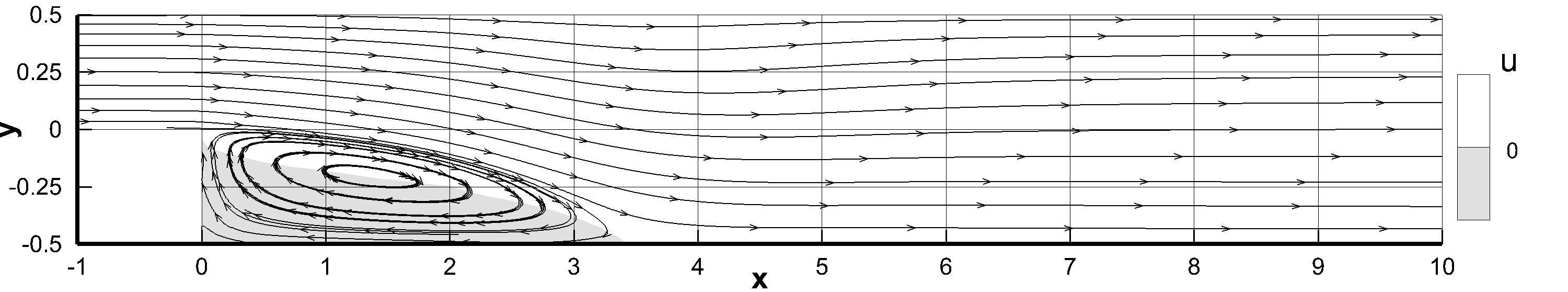}\\
			\includegraphics[width=0.95\textwidth]{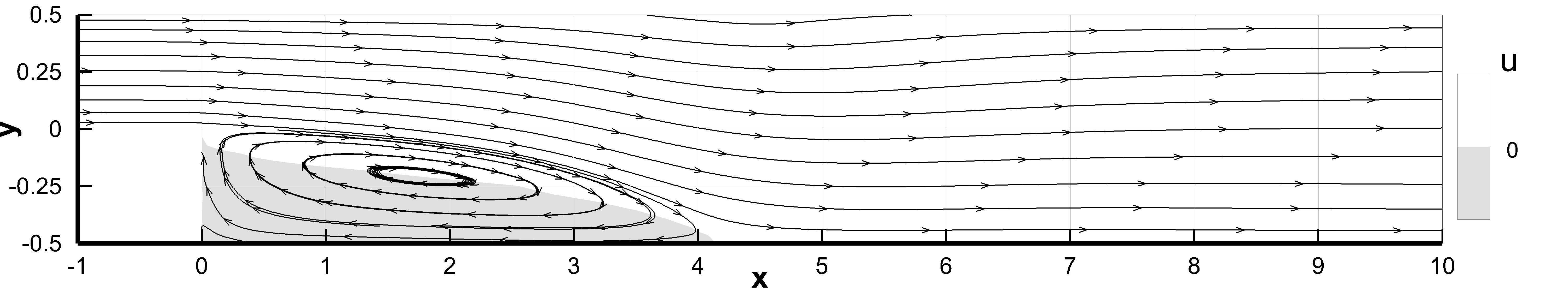}\\
			\includegraphics[width=0.95\textwidth]{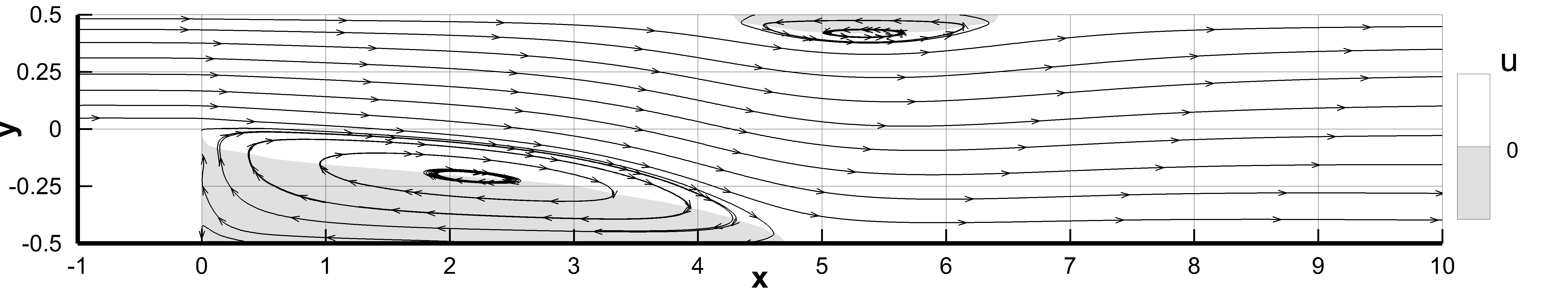}\\
			\includegraphics[width=0.95\textwidth]{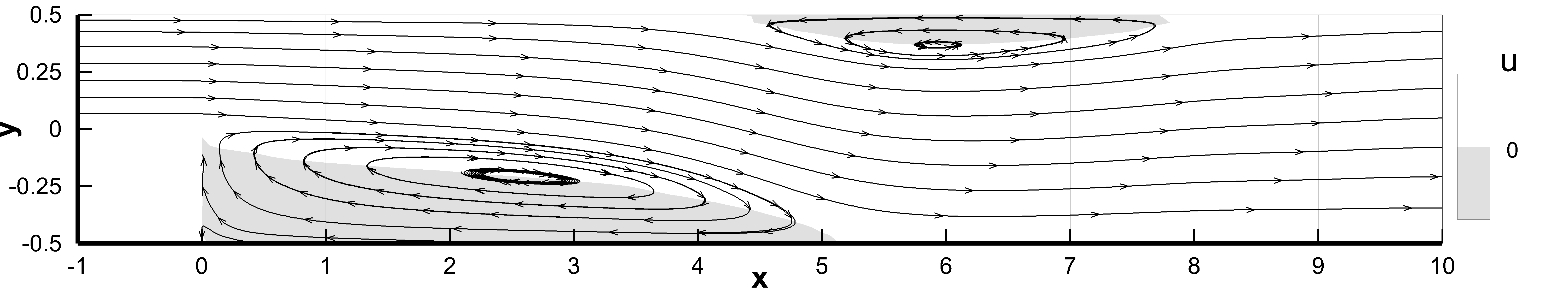}\\
			\includegraphics[width=0.95\textwidth]{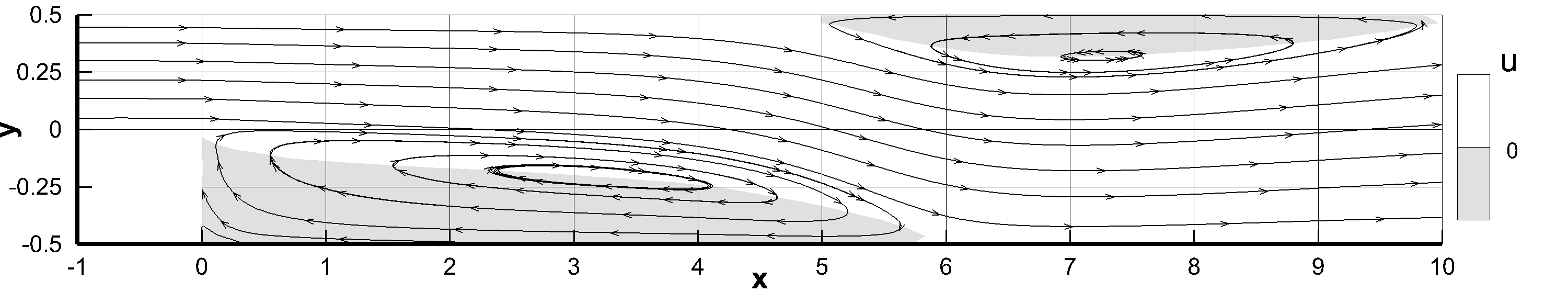}
\caption{The numerical solution obtained for the two-dimensional backward facing step problem at different Reynolds numbers, respectively,  from the top to the bottom: Re=$100$, Re=$200$, Re=$300$, Re=$400$, Re=$500$, Re=$600$, and Re=$800$ obtained with the staggered semi-implicit spectral DG-$\p_{6}$ method. Recirculations are highlighted by the sign of the axial velocity $u$.}\label{fig:BFStep2D}
\end{figure}

Another typical benchmark problem for testing the accuracy of numerical methods in computational fluid dynamics is the backward facing step problem. In this test a flow separation is induced by 
a sudden backward step inside a two dimensional duct. A main recirculation zone is generated next to the step, starting already at low Reynolds numbers. Then, by increasing the Reynolds 
number, new secondary recirculations are generated. A non-zero velocity $u=1$ is imposed at the entrance, a constant pressure $p=0$ is imposed at the outflow. In this case the axial 
spatial domain is $x\in[-10,20]$, the height of the two dimensional duct is $h_{\text{in}}=0.5$ at the entrance and $h_{\text{out}}=h_{\text{in}}+h_{\text{s}}=1.0$ at the exit, 
with an expansion ratio ER$=h_{\text{out}}/h_{\text{in}}=2$ at $x=0$, i.e. a backward facing step of height $h_s=0.5$. The spatial domain is discretized with elements of dimension 
$\Delta x = 1$, $\Delta y = 0.25$, the implicitness factor in time is taken as $\theta = 0.6$. 
Fig. \ref{fig:BFStep2D} shows the streamlines and the recirculation patterns obtained for different Reynolds number up to $Re=800$ with the $\p_{6}$ version of our staggered semi-implicit 
spectral DG method. The numerical results are compared with the two dimensional reference data provided in \cite{Erturk2008} and with the experimental measurements of 
\cite{LeeMateescu1998}. A good agreement is achieved. The plotted data in Fig. \ref{fig:BFStep2D} show some discrepancies between the two dimensional simulations and  the experimental data 
\cite{LeeMateescu1998} that become more visible at higher Reynolds number. These differences are due to three dimensional effects that are introduced by the sidewalls at higher Reynolds numbers,
as discussed in \cite{Tylli2002,Armaly1983,Mouza2005,Rani2007}.

\subsubsection{Backward facing step: 3D}
\begin{figure} 
\centering %subfloat
			\includegraphics[width=0.48\textwidth]{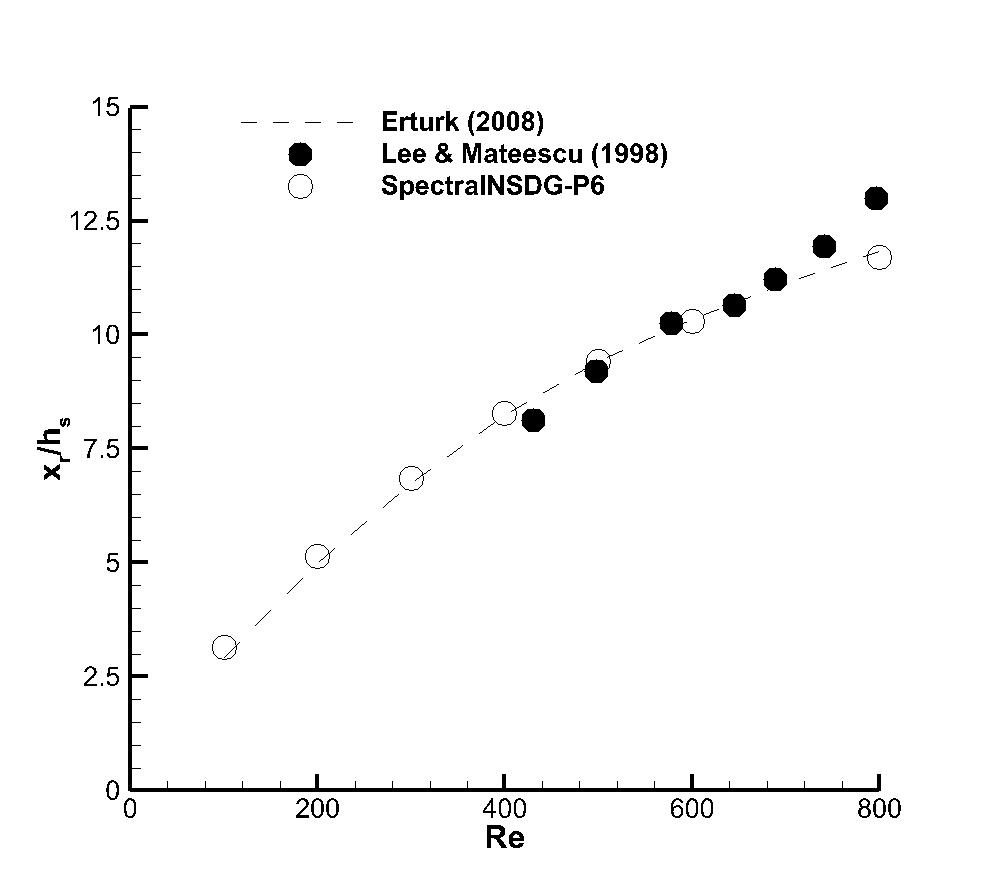}
\caption{Position of the reattachment point obtained with the staggered semi-implicit spectral DG-$\p_{6}$ method for the two dimensional backward facing step problem compared with the two dimensional 
numerical results of \cite{Erturk2008} and the experimental measurements of \cite{LeeMateescu1998} at different Reynolds numbers, in the range $Re\in(0,800)$.}\label{fig:BFStep_data} 
\end{figure}

\begin{figure} 
\centering %subfloat
			\includegraphics[width=\textwidth]{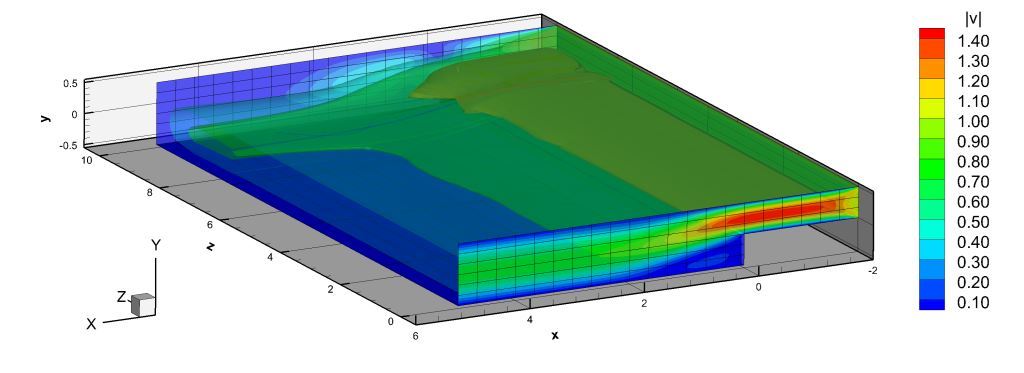}\\
			\includegraphics[width=\textwidth]{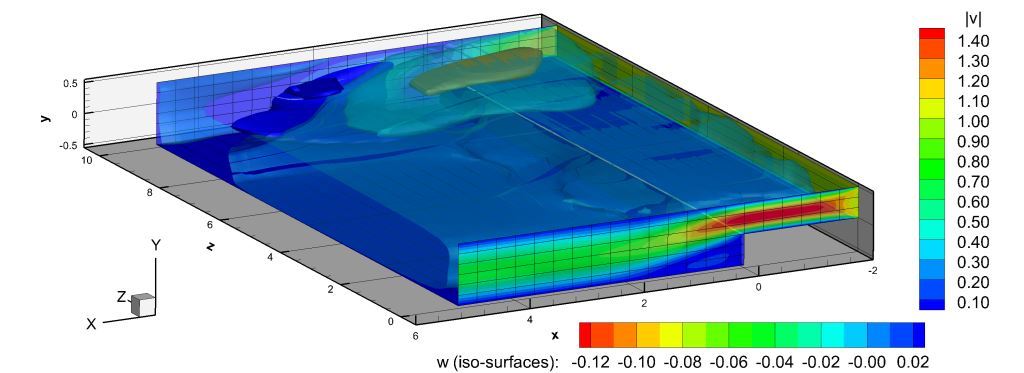} \\
			\includegraphics[width=\textwidth]{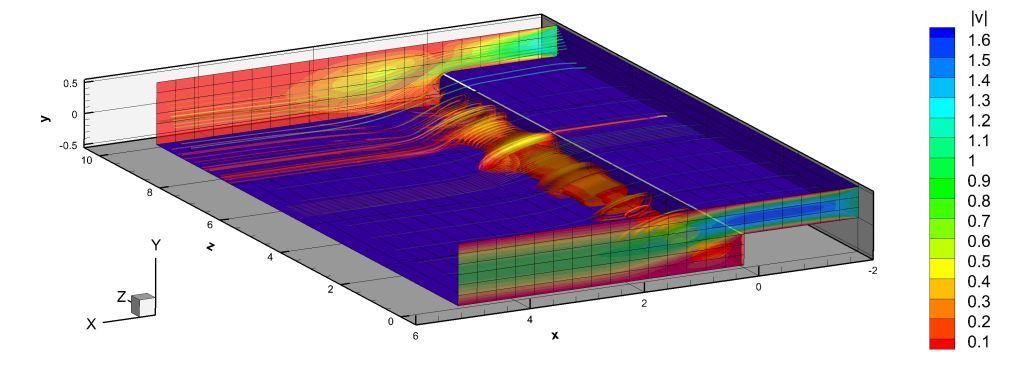}
			\caption{Numerical solution for the three dimensional backward facing step problem at time $t=25.0$ computed with the staggered semi-implicit spectral DG-$\mathbb{P}_{3}$ method for $Re=100$. 
			The iso-surfaces of the  velocity magnitude (top), the iso-surfaces of the $w$ velocity component (center) and the streamtraces of the fluid flow (bottom) are plotted for the first half of the  spatial domain $z>0$. }\label{fig:BFS3D_Re100}
\end{figure}

\begin{figure} 
\centering %subfloat
			\includegraphics[width=\textwidth]{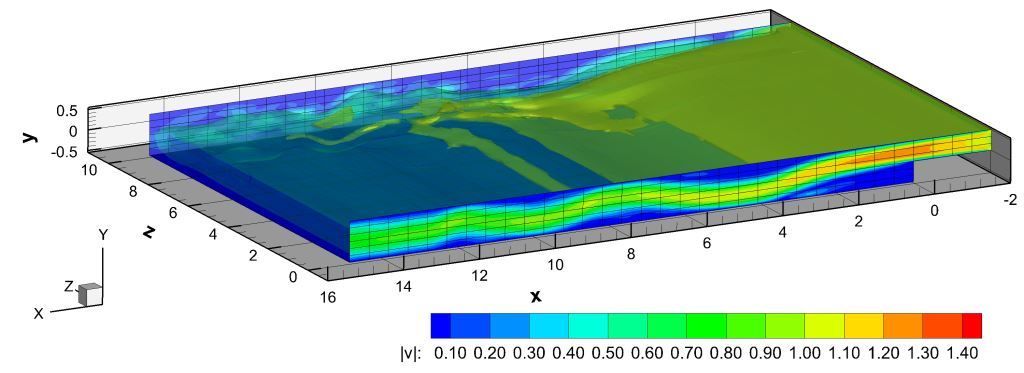}\\
			\includegraphics[width=\textwidth]{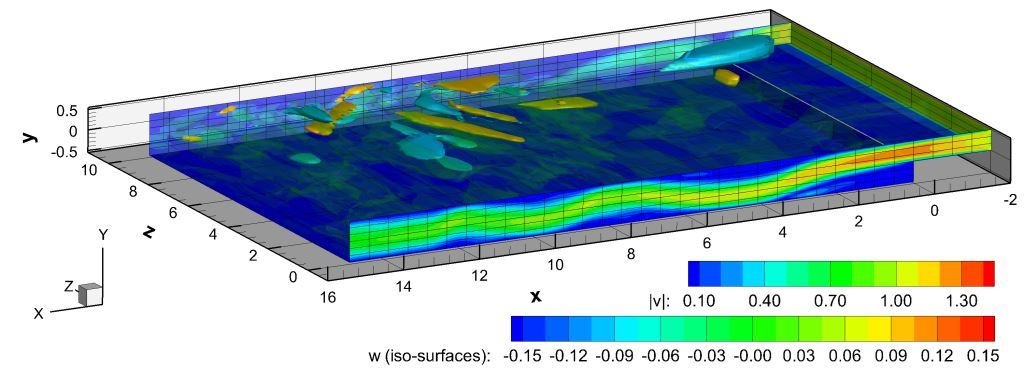} \\
			\includegraphics[width=\textwidth]{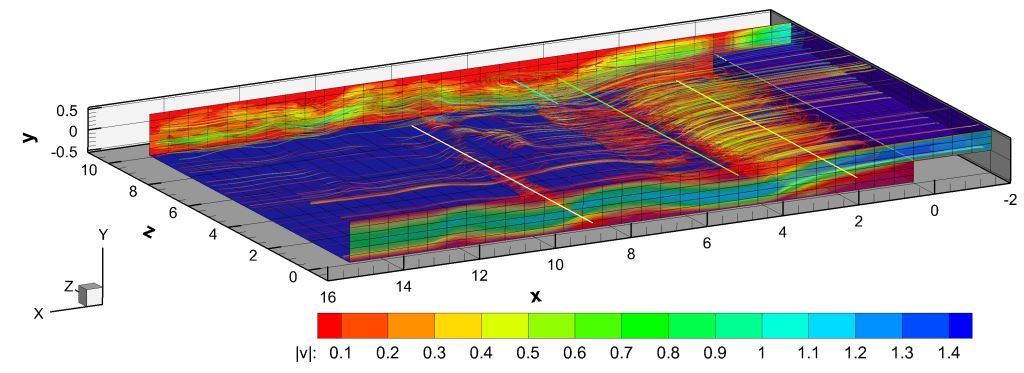}
			\caption{Numerical solution for the three dimensional backward facing step problem at time $t=25.0$ computed with the staggered semi-implicit spectral DG-$\mathbb{P}_{3}$ method for $Re=1000$. 
			The iso-surfaces of the velocity magnitude (top), the iso-surfaces of the $w$ velocity component (center) and the streamtraces of the fluid flow (bottom) are plotted for the first half of the  			
			spatial domain $z>0$. The main recirculation axes are highlighted in the figure at the bottom. }\label{fig:BFS3D_Re1000}
\end{figure}

In this section, the numerical results of the simulation of the three dimensional extension of the backward facing step problem are shown and discussed. The physical domain is described by an expansion-ratio ER$=2$, and an aspect-ratio AR$=L_z/H=40$, where $L_z$ is the width of the duct in the third spatial dimension. 
As mentioned above, the two dimensional results are show differences compared to the experimental data for higher Reynolds numbers. The main reason is that the lateral boundary layers developing 
on the side walls interact with the main recirculations of the two-dimensional flow. This interpretation is justified by the fact that at lower Reynolds numbers and higher aspect ratio, 
i.e. when the aforementioned interactions are negligible, the two-dimensional results actually match the three dimensional ones and the experimental data (see Figure \ref{fig:BFStep_data}).  
The numerical solutions for $Re=100$ (laminar regime) and Re=$1000$ (transitional regime) at time $t=25.0$ obtained with our spectral SIDG-$\p_{3}$ scheme give an overview of the 3D 
flow field, see Figures \ref{fig:BFS3D_Re100}-\ref{fig:BFS3D_Re1000}. 
The friction forces at the lateral boundary layers constrict the axial velocity profile and the main recirculation to the center of the duct. A non-zero $w$ velocity component is generated 
consequently.

\subsection{Three dimensional Taylor-Green vortex problem}
\begin{figure} 
\centering %subfloat
			\includegraphics[width=0.6\textwidth]{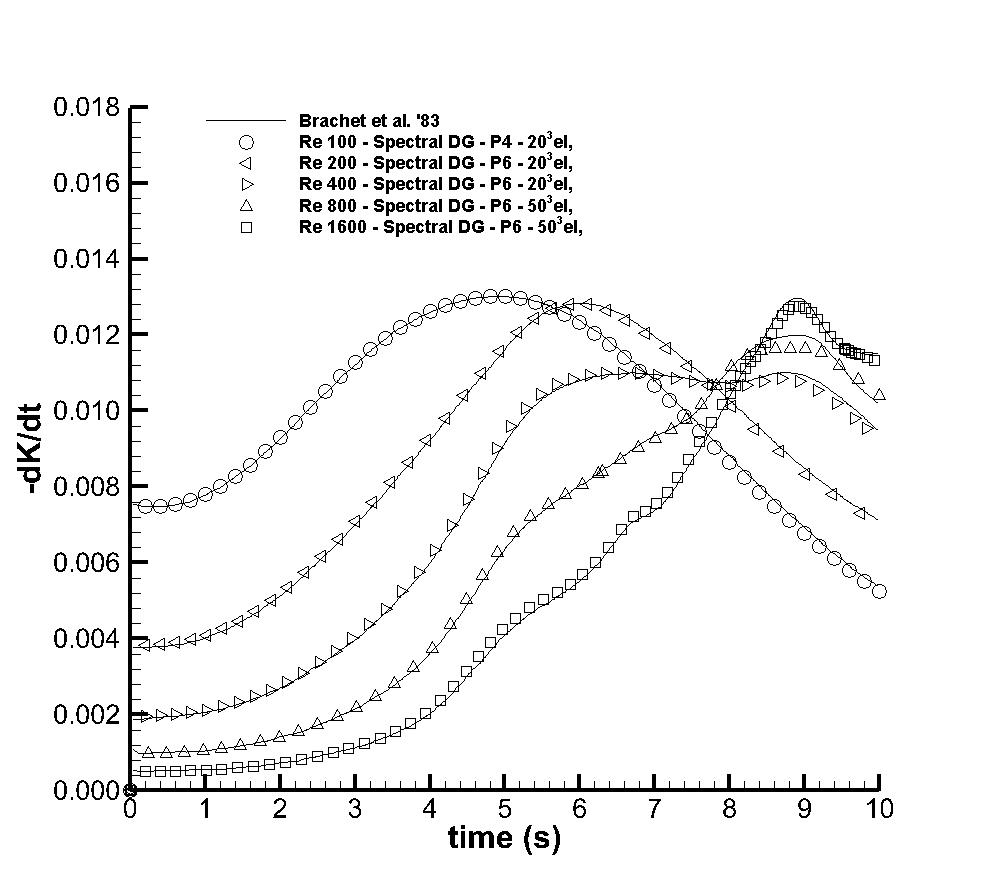}
			\caption{Time evolution of the kinetic energy dissipation rate $\epsilon(t)$ obtained with staggered semi-implicit spectral DG-$\p_{N}$ schemes 
			at different Reynolds numbers $100<Re<1600$. The DNS reference solutions of Brachet et al. \cite{Brachet1983} are plotted as continuous lines.} \label{fig:TG3D}
\end{figure}

\begin{figure} 
\centering %subfloat
			\includegraphics[width=0.32\textwidth]{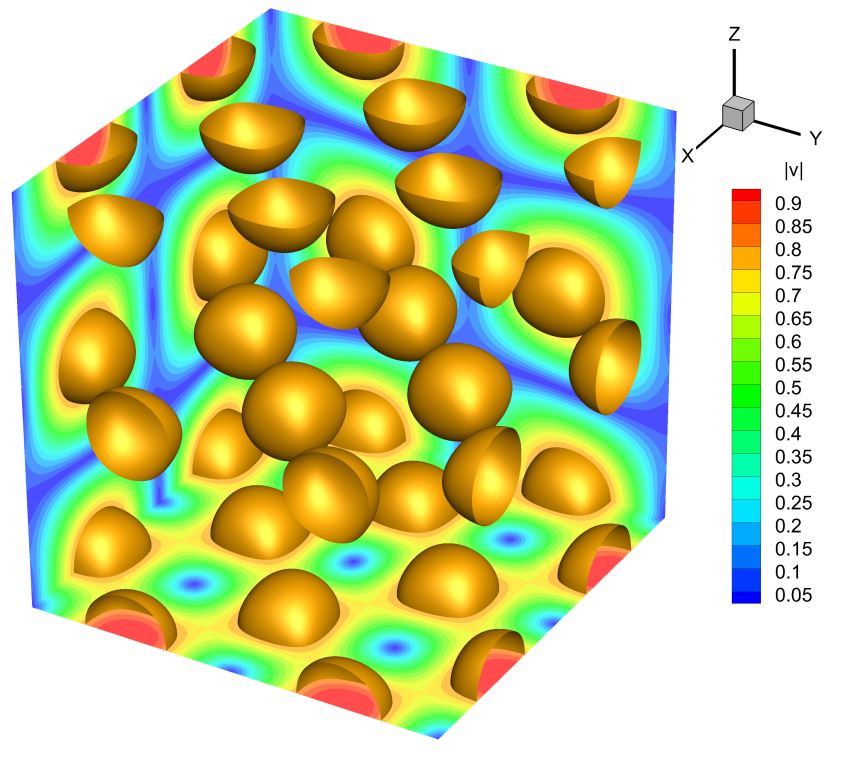}
			\includegraphics[width=0.32\textwidth]{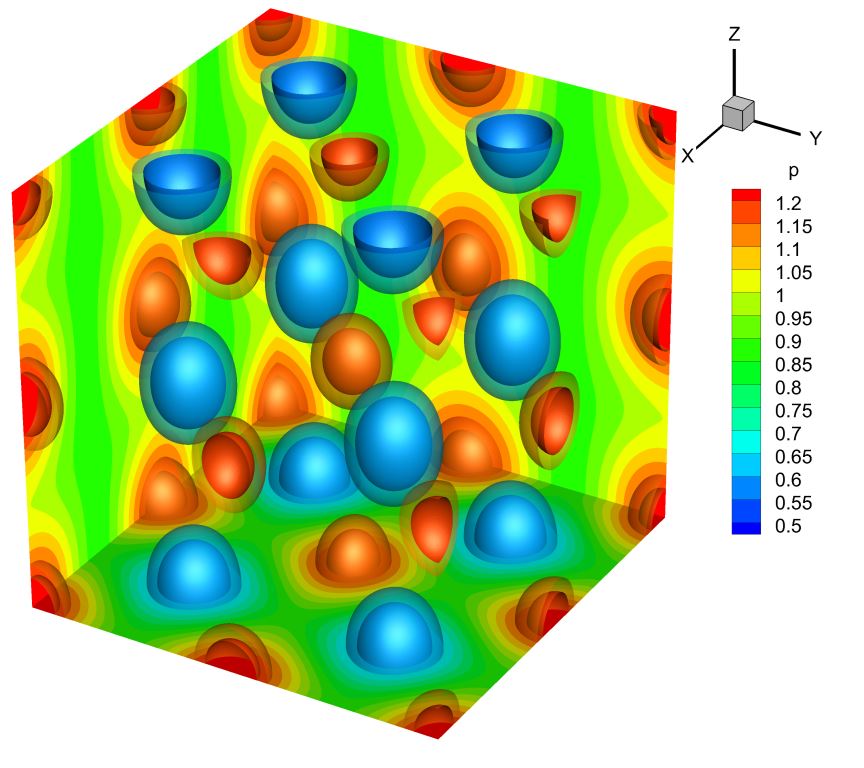}
			\includegraphics[width=0.32\textwidth]{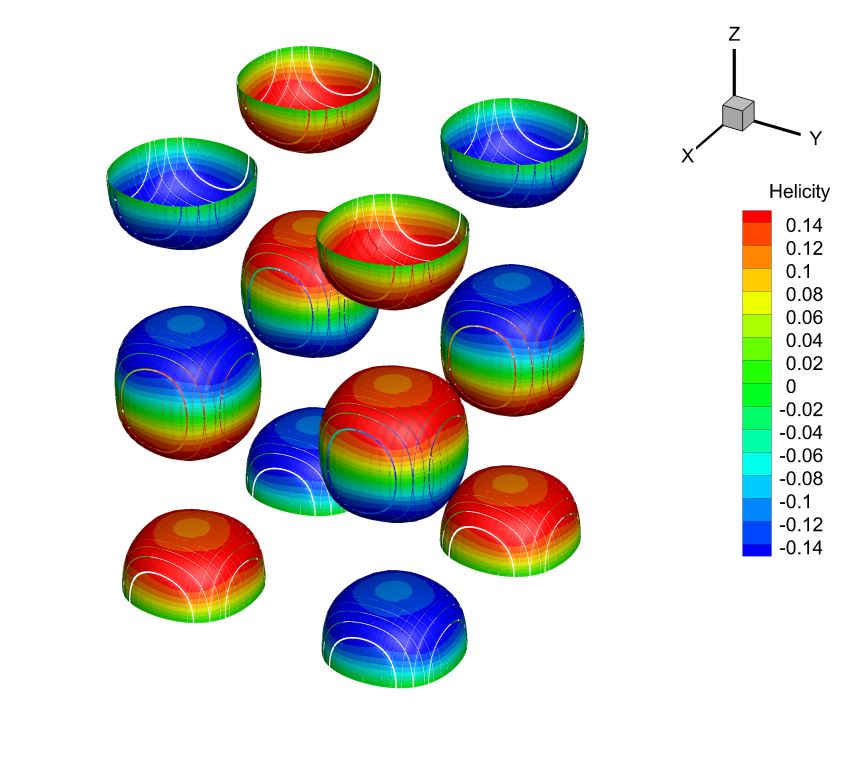}\\
			\includegraphics[width=0.32\textwidth]{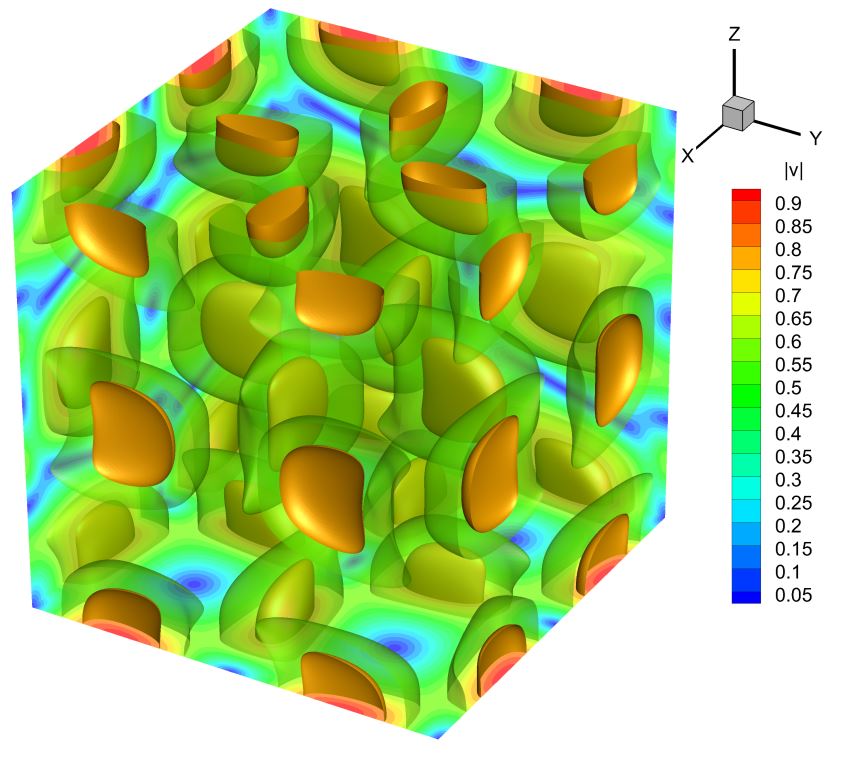}
			\includegraphics[width=0.32\textwidth]{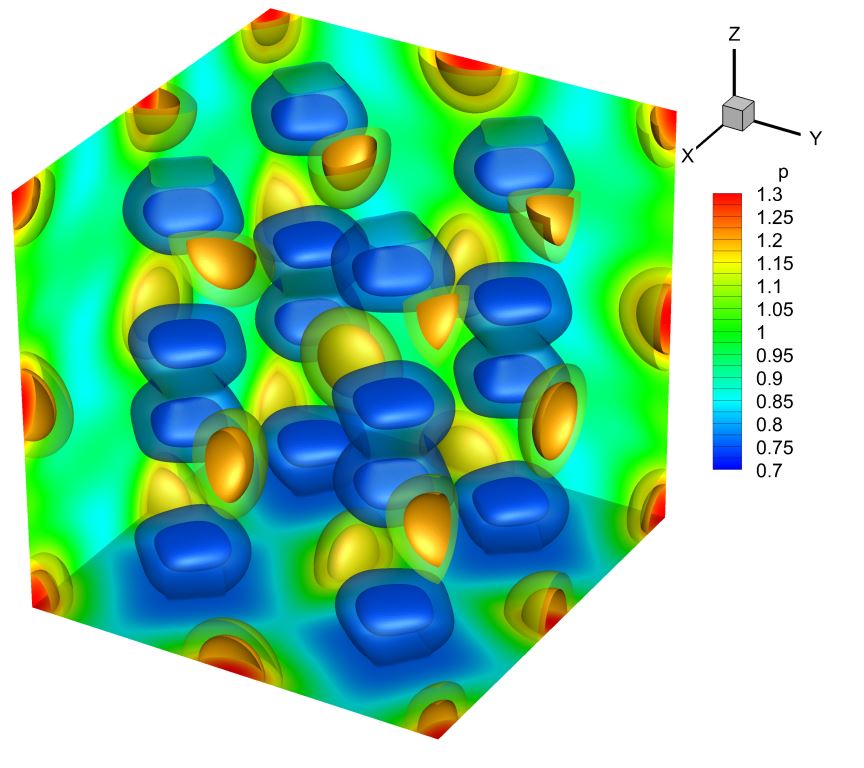}
			\includegraphics[width=0.32\textwidth]{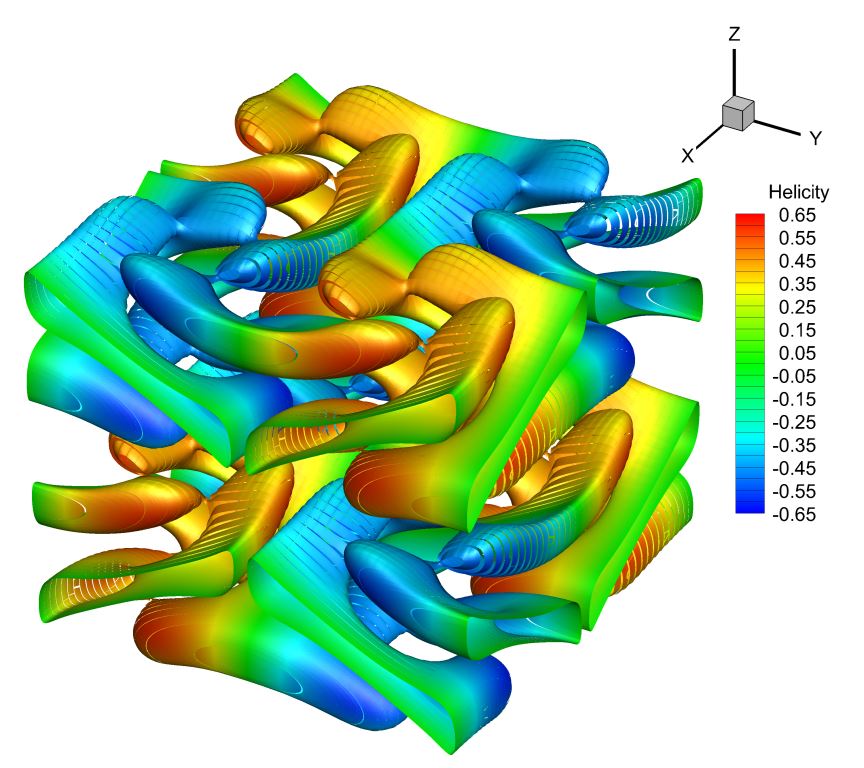}\\
			\includegraphics[width=0.32\textwidth]{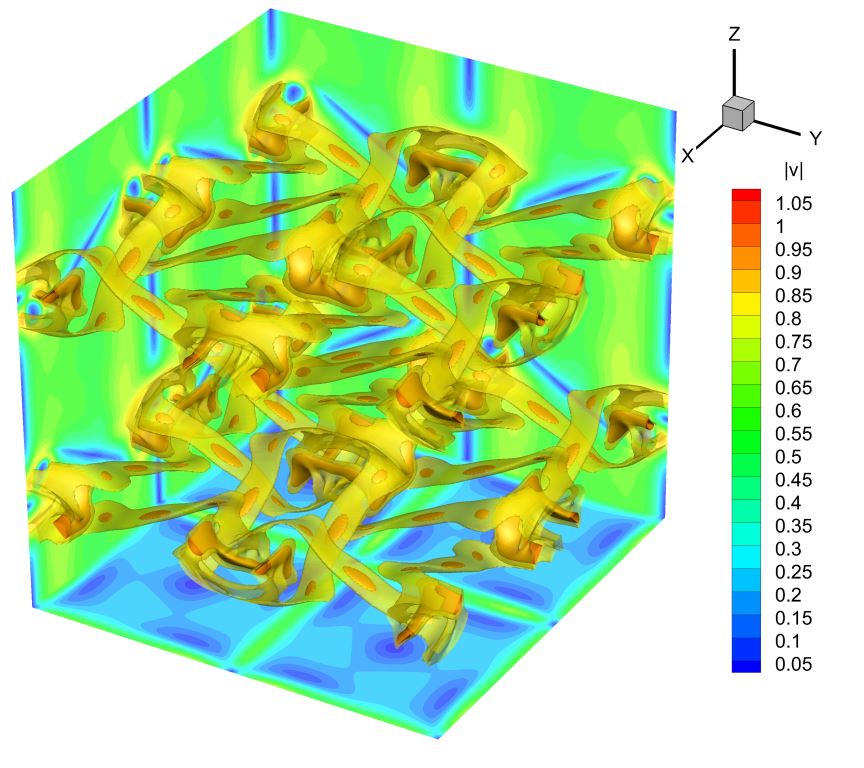}
			\includegraphics[width=0.32\textwidth]{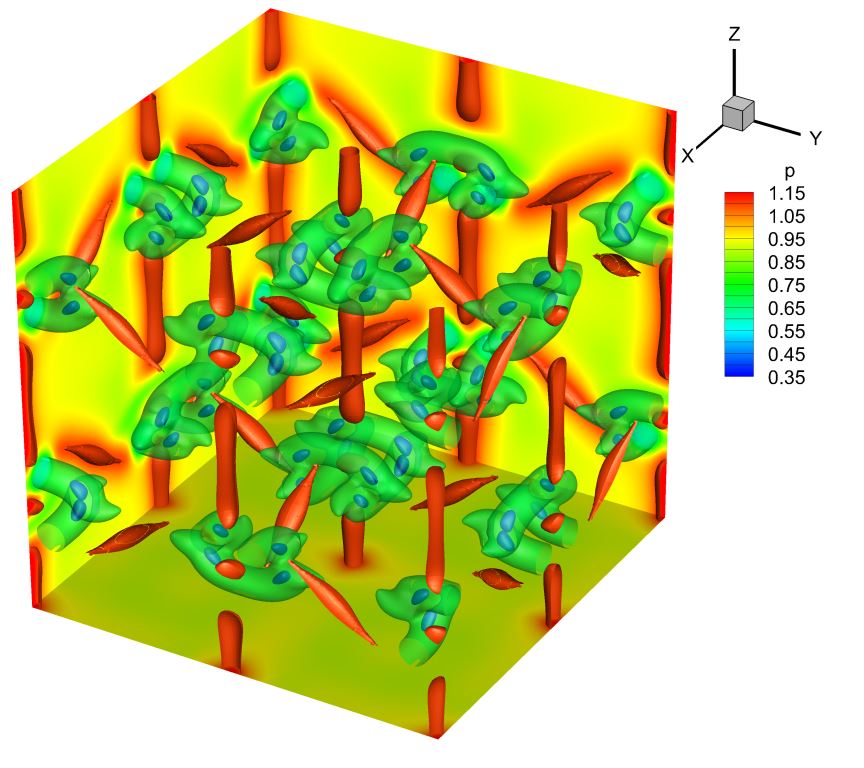}
			\includegraphics[width=0.32\textwidth]{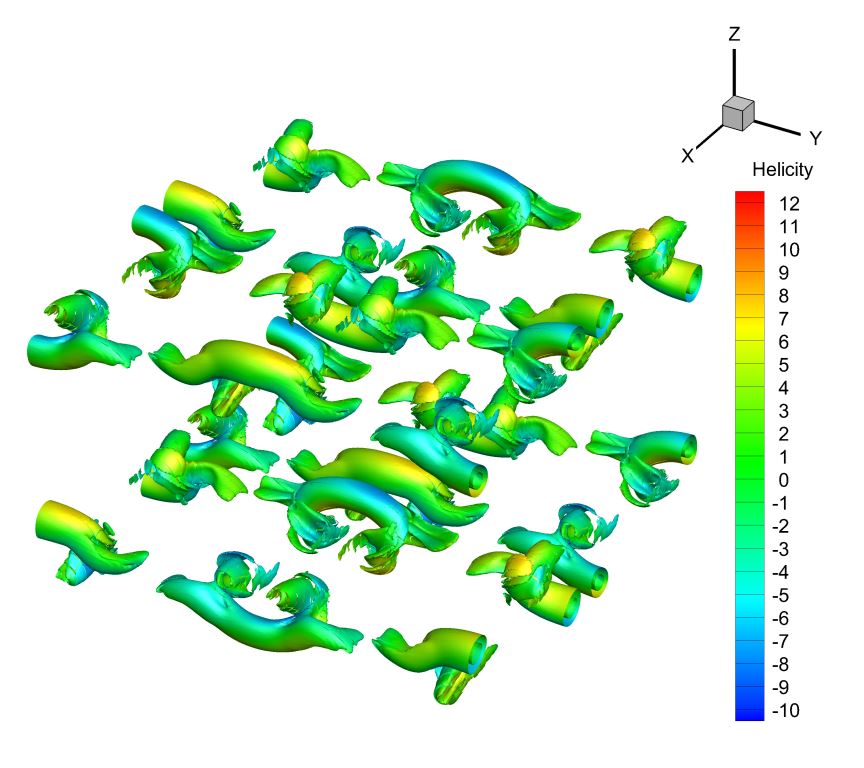}\\
			\includegraphics[width=0.32\textwidth]{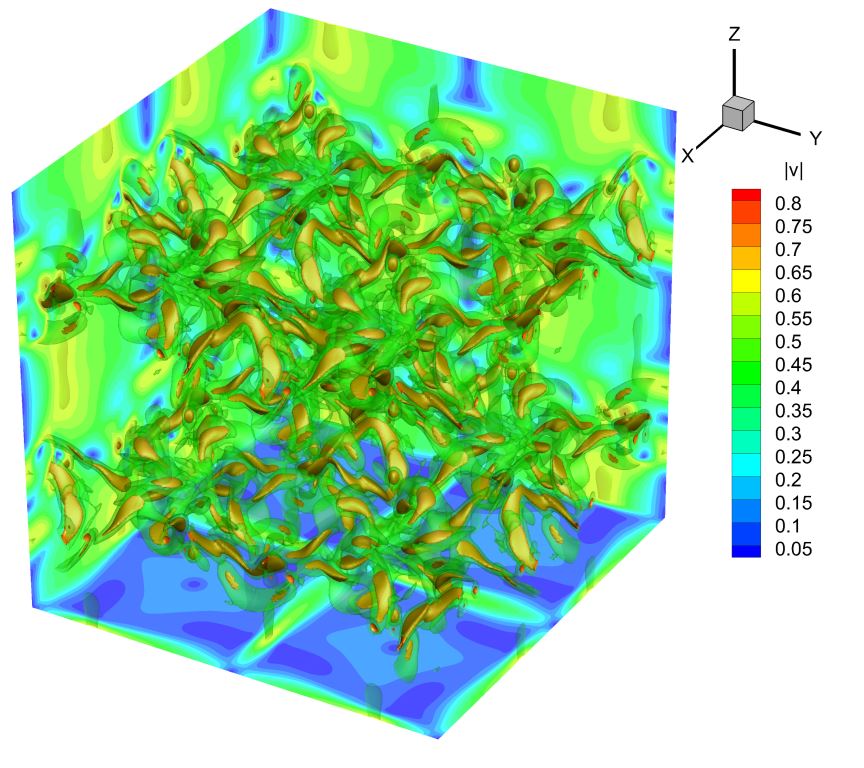}
			\includegraphics[width=0.32\textwidth]{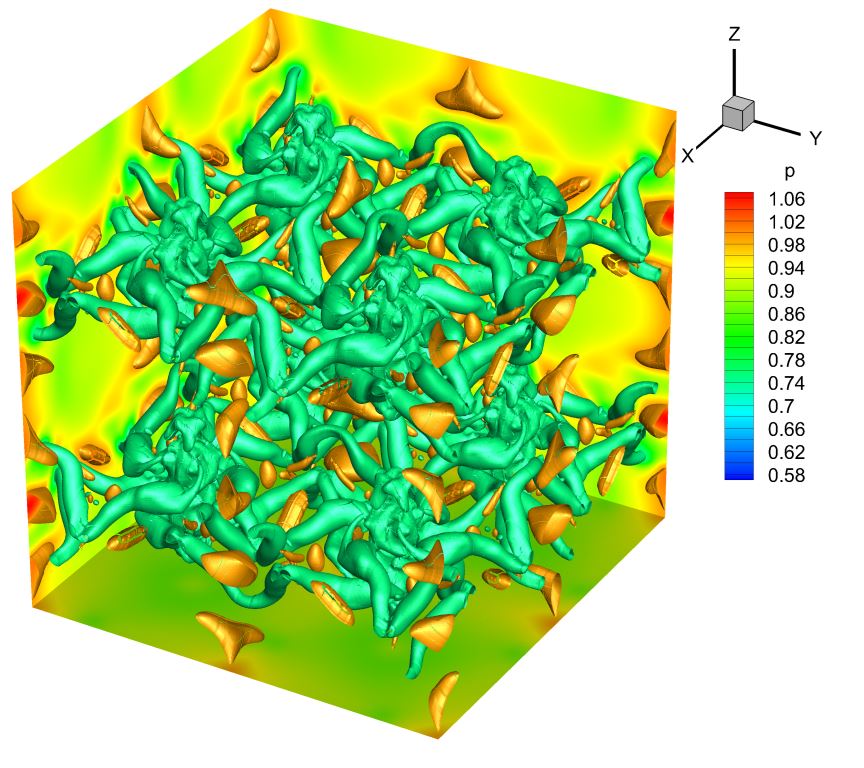}
			\includegraphics[width=0.32\textwidth]{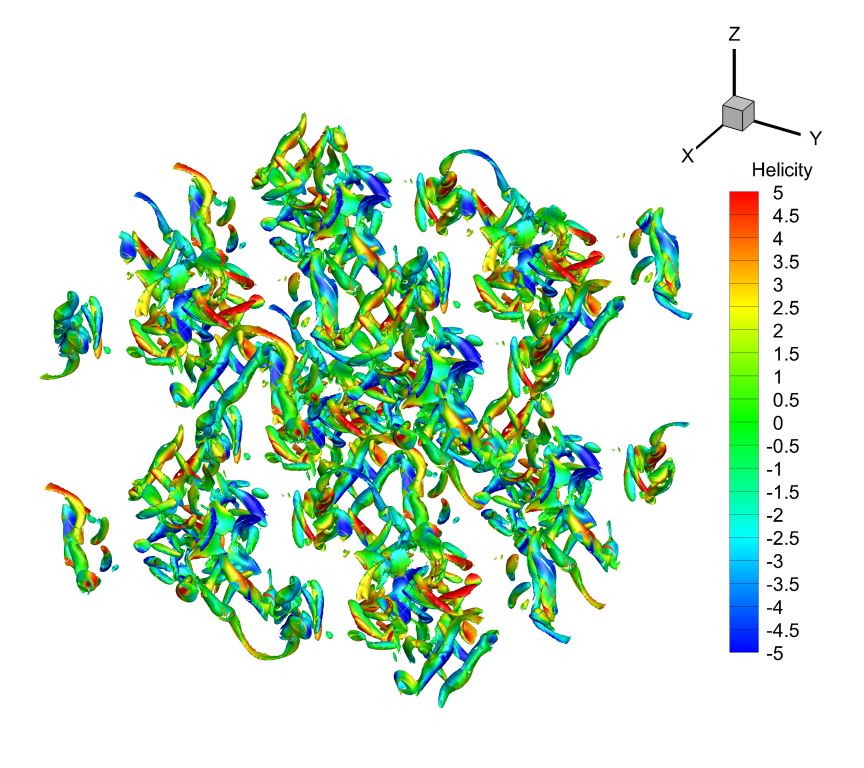}
			\caption{Numerical solution for the three dimensional Taylor-Green vortex flow at $Re=800$ computed with the staggered semi-implicit spectral DG-$\p_{6}$ method using $50^3$ elements.  
			The isosurfaces of the velocity (left), the isosurfaces of the pressure (center) and the isosurfaces of the vorticity colored by the helicity field (right) are plotted at 
			times $t=0.4$, $2.0$, $6.0$ and $10.0$ from the top to the bottom, respectively.}\label{fig:TG3D2}
\end{figure}

A classical fully three-dimensional flow that is widely used for testing the ability of a numerical method in solving the smallest scales in turbulent flows is the three dimensional Taylor-Green 
vortex problem. In this test the velocity and pressure field are initialized with 
\begin{align}
&u(x,y,z,0) = \sin(x)\cos(y)\cos(z),\\
&v(x,y,z,0) = -\cos(x)\sin(y)\cos(z),\\
&w(x,y,z,0) = 0,\\
&p(x,y,z,0) = \frac{1}{16}\left(\cos(2x)+ \cos(2y) \right)\left(\cos(2z)+2\right).
\end{align}
The resulting fluid flow is initially smooth and laminar, but the nonlinearity in the governing PDE due to the convective terms combined with a small viscosity quickly 
generates complex small-scale flow structures after finite times. 
A widely accepted reference solution for the rate of kinetic energy dissipation has been computed for this test problem by Brachet et al. in \cite{Brachet1983} through both a direct spectral 
method based on up to $256^3$ modes and a rigorous power series analysis up to order $t^{80}$ (see also \cite{Morf1980}). The computational domain is chosen as $\Omega=[0,2\pi]^3$, 
with periodic boundary conditions everywhere. The smaller the expected flow scales, the higher the necessary grid resolution. 
The time evolution of the main physical variables of the fluid flow is represented in Figure \ref{fig:TG3D2} at times  $t=0.4$, $2.0$, $6.0$ and $10.0$ for the case $Re=800$. 
%The streamlines along the three orthogonal planes ${x-y}$, ${x-z}$ and ${y-z}$ are presented in Figure \ref{fig:TG3D2b}. 
The complexity of the resulting small scale flow structures is clearly visible. In order to compare our results quantitatively with those of Brachet et al. \cite{Brachet1983}, we compute
the rate of kinetic energy dissipation 
\begin{align}
\epsilon(t) = - \frac{\partial K}{\partial t}= - \frac{1}{\left\| \Omega \right\|} \frac{\partial}{\partial t}  \int_{\Omega} \frac{1}{2} \mathbf{v}^2  d\x.
\end{align}  
Especially when the rate $\epsilon$ reaches its maximum, a high-resolution method together with a sufficiently fine grid is needed in order to resolve the flow physics properly. 
Figure \ref{fig:TG3D} shows the time evolution of the rate of the global kinetic energy dissipation $\epsilon(t)$ for different Reynolds numbers $Re=100$, Re=$200$, Re=$400$, Re=$800$ 
and Re=$1600$, obtained with our semi-implicit staggered spectral DG-$\p_{4}$ and -$\p_{6}$ schemes, along with $20^3$ and $50^3$ elements, respectively, see Figure \ref{fig:TG3D}. 
The computed results fit the DNS reference data very well, confirming that our scheme is able to resolve even the smallest flow scales properly up to $Re=1600$.  
%Figure \ref{fig:TG3D2c} shows an example of small scale structure that is generated by the flow at $Re=800$,  identified by the isosurfaces of the vorticity $\omega = 7$ and colored 
%by the helicity magnitude $\mathcal{h}=\mathbf{v}\cdot\omega$. 

\section{Spectral space-time DG schemes on staggered Cartesian grids}
 \label{sec:HO}
In this section the high-order DG formulation is extended to the time dimension by looking for discrete solutions $(u_h,v_h,w_h,p_h)$ under the form of linear combinations of \emph{piecewise space-time polynomials} of maximum degree $N$ in space and $M$ in time with respect to a reference basis $\mathpzc{B}_N$ for the spatial dependency and $\mathpzc{B}_M$ for the time dependency. In the following, the  general mathematical framework is outlined, and several numerical tests are performed in two and three space dimensions, with the aim of assessing the efficiency and the accuracy of the proposed high 
order accurate staggered spectral space-time DG scheme.

\subsection{Presentation of the numerical scheme}
 \label{sec:basicHO}
By using the same nomenclature as before, the weak formulation of the governing equations (\ref{eq:wNSmom}-\ref{eq:wNSinc}) in \emph{space-time} reads
\begin{align}
%&\int \limits_{\Omega_x} \mathit{f_u} \frac{\partial u}{\partial t} + \nabla \cdot \mathbf{F_u} + \partial_x p  = 0 \label{eq:NSmom}\\
%&\int \limits_{\Omega_y} \mathit{f_v} \frac{\partial v}{\partial t} + \nabla \cdot \mathbf{F_v} + \partial_y p = 0 \label{eq:NSmom}\\
%&\int \limits_{\Omega_z} \mathit{f_w} \frac{\partial w}{\partial t} + \nabla \cdot \mathbf{F_w} + \partial_z p = 0 \label{eq:NSmom}\\
\begin{array}{lcr}
\int \limits_{d\Omega^*_x\times T_{n+1}}  \mathit{f_u} \left( \frac{\partial u}{\partial t} + \nabla \cdot \mathbf{F}_u + \partial_x p \right) = 0, & &
\\ \int \limits_{d\Omega^*_y\times T_{n+1}} \mathit{f_v} \left( \frac{\partial v}{\partial t} + \nabla \cdot \mathbf{F}_v + \partial_y p \right) = 0, & \hspace{3cm} &
 \int \limits_{d\Omega\times T_{n+1}} \mathit{f} \, \left( \nabla \cdot \mathbf{v}\right) = 0, \\ %\\ 
\int \limits_{d\Omega^*_z\times T_{n+1}} \mathit{f_w} \left( \frac{\partial w}{\partial t} + \nabla \cdot \mathbf{F}_w + \partial_z p \right)= 0, & &
%\\       
\end{array}\label{eq:STNSeq}
\end{align}
where $T_{n+1}= \left[t_n,t_{n+1}\right]$ is the future time interval where the solution is unknown.
 Then, the definitions of the piecewise polynomials (\ref{eq:uDGpoly}-\ref{eq:pDGpoly}) are augmented by
\begin{align}
  \left\{  
\begin{array}{rcl}
 \hat{u}_{pqr}(t) &=& \sum_{s=0}^M \omega_s(t) \, \hat{u}_{pqrs}, \\ 
 \hat{v}_{pqr}(t) &=& \sum_{s=0}^M \omega_s(t) \, \hat{v}_{pqrs}, \\
 \hat{w}_{pqr}(t) &=& \sum_{s=0}^M \omega_s(t) \, \hat{w}_{pqrs}, \\  
 \hat{p}_{pqr}(t) &=& \sum_{s=0}^M \omega_s(t) \, \hat{p}_{pqrs},
\end{array}
\right.
\;\;\;\;\; 
\left.  
\begin{array}{l} 
   \text{with}\;\; t\in T_{n+1},   \\
 n=1,...,N_t,  \\
 u_{lmnp},v_{lmnp},w_{lmnp},p_{lmnp}\in \varmathbb{R},  
\end{array}
\right. \label{eq:spacetimepoly}   
\end{align}
where  $\omega_s(t)$ is a polynomial in time, generated from the basis functions $\varphi_s\in\mathpzc{B}_M$ with the rule 
\begin{align*}
&\omega(t) = \varphi(\xi), \;\;\; \textnormal{with} \; t = t_{n}+\xi \Delta t, \; 0 \leq \xi \leq 1.
\end{align*}
Then a \emph{staggered space-time DG discretization} of the incompressible Navier-Stokes equations reads 
\begin{align}
&\int \limits_{T_{i+\frac{1}{2},j,k} \times T_{n+1}}  \psi_{m_1}(x)\omega_{m_2}(y)\omega_{m_3}(z)\omega_{m_0}(t) \left( \frac{\partial u_{h}}{\partial t} + \nabla \cdot \mathbf{F}_{u} + \partial_x p_{h} \right) d\x dt = 0,\label{eq:uDGNSmomST} \\
&\int \limits_{T_{i,j+\frac{1}{2},k} \times T_{n+1}} \omega_{m_1}(x)\psi_{m_2}(y)\omega_{m_3}(z)\omega_{m_0}(t) \left( \frac{\partial v_{h}}{\partial t} + \nabla \cdot \mathbf{F}_{v} + \partial_y p_{h} \right) d\x dt = 0,\label{eq:vDGNSmomST}\\
&\int \limits_{T_{i,j,k+\frac{1}{2}} \times T_{n+1}} \omega_{m_1}(x)\omega_{m_2}(y)\psi_{m_3}(z)\omega_{m_0}(t) \left( \frac{\partial w_{h}}{\partial t} + \nabla \cdot \mathbf{F}_{w} + \partial_z p_{h} \right) d\x dt = 0, \label{eq:wDGNSmomST}\\
&\int \limits_{T_{i,j,k} \times T_{n+1}} \omega_{m_1}(x)\omega_{m_2}(y)\omega_{m_3}(z)\omega_{m_0}(t) \, \left( \nabla \cdot \mathbf{v}_{h}\right) d\x dt = 0.\label{eq:DGNSincST}
\end{align}
Again, we need to account for the jump of the pressure inside the velocity control volumes, and we perform integration by parts in space of equation \eqref{eq:DGNSincST}, which thus 
becomes 
\begin{equation}
  \int \limits_{\partial T_{i,j,k} \times T_{n+1} } \omega_{m_1}(x)\omega_{m_2}(y)\omega_{m_3}(z) \omega_{m_0}(t) \, \mathbf{v}_{h} \cdot \vec{n} \, dS dt 
	- \int \limits_{T_{i,j,k} \times T_{n+1} } \nabla \left(\omega_{m_1}(x)\omega_{m_2}(y)\omega_{m_3}(z) \omega_{m_0}(t)  \right) \cdot \mathbf{v}_{h}(\x,t) \, d\x dt = 0. \label{eq:DGNSincST2}
\end{equation} 
Eqn. \eqref{eq:DGNSincST2} is again well defined, since the velocity vector $\mathbf{v}_h$ is continuous across the element boundary $\partial T_{i,j,k} \times T_{n+1}$, thanks to the use 
of a staggered grid approach. Since $p_h$ is discontinuous inside the domains of integration of the momentum equations (\ref{eq:uDGNSmomST}-\ref{eq:wDGNSmomST}) the following \emph{jump} 
terms arise % from the integration
\begin{align}
\int \limits_{t^n}^{t^{n+1}} \int\limits_{x_i}^{x_{i+1}} \int \limits_{y_{j-\frac{1}{2}}}^{y_{j+\frac{1}{2}}} \int \limits_{z_{k-\frac{1}{2}}}^{z_{k+\frac{1}{2}}}  
\psi_{m_1}(x) \omega_{m_2}(y) \omega_{m_3}(z) \omega_{m_0}(t) \partial_x p_h(\mathbf{x},t) \, d\x dt = && \nonumber \\   
\int \limits_{t^n}^{t^{n+1}}  \int\limits_{x_i}^{x_{i+\frac{1}{2}}} \int \limits_{y_{j-\frac{1}{2}}}^{y_{j+\frac{1}{2}}} \int \limits_{z_{k-\frac{1}{2}}}^{z_{k+\frac{1}{2}}}  
\psi_{m_1}(x) \omega_{m_2}(y) \omega_{m_3}(z) \omega_{m_0}(t) \partial_x p_{i,j,k}(\mathbf{x},t) \, d\x dt + \nonumber \\ 
\int \limits_{t^n}^{t^{n+1}} \int\limits_{x_{i+\frac{1}{2}}}^{x_{i+1}} \int \limits_{y_{j-\frac{1}{2}}}^{y_{j+\frac{1}{2}}} \int \limits_{z_{k-\frac{1}{2}}}^{z_{k+\frac{1}{2}}}  
\psi_{m_1}(x) \omega_{m_2}(y) \omega_{m_3}(z) \omega_{m_0}(t) \partial_x p_{i+1,j,k}(\mathbf{x},t) \, d\x dt + && \nonumber \\
+ \int \limits_{t^n}^{t^{n+1}} \int \limits_{y_{j-\frac{1}{2}}}^{y_{j+\frac{1}{2}}} \int \limits_{z_{k-\frac{1}{2}}}^{z_{k+\frac{1}{2}}}  
\psi_{m_1}(x_{i+\frac{1}{2}}) \omega_{m_2}(y) \omega_{m_3}(z) \omega_{m_0}(t) \left(p_{i+1,j,k}(x_{i+\frac{1}{2}},y,z,t) - p_{i,j,k}(x_{i+\frac{1}{2}},y,z,t)\right) dy dz dt, && \label{eq:pressurejumpST}%\\
%\int\limits_{y_j}^{y_{j+1}} dy\,\psi_{m'}(y) \partial_y p = \int\limits_{y_j}^{y_{j+\frac{1}{2}}}dy\, \psi_{m'}(y) \partial_y p^j + \int\limits_{y_{j+\frac{1}{2}}}^{y_{j+1}}dy\, \psi_{m'}(y) \partial_y p^{j+1}   + \psi_{m'}(y_{j+\frac{1}{2}})\left(p^{j+1}(y_{j+\frac{1}{2}}) - p^j(y_{j+\frac{1}{2}})\right)\\
%\int\limits_{z_k}^{z_{k+1}}dz\, \psi_{m''}(z) \partial_z p = \int\limits_{z_k}^{z_{k+\frac{1}{2}}}dz\, \psi_{m''}(z) \partial_z p^k + \int\limits_{z_{k+\frac{1}{2}}}^{z_{k+1}}dz\, \psi_{m''}(z) \partial_z p^{k+1}   + \psi_{m''}(z_{k+\frac{1}{2}})\left(p^{k+1}(z_{k+\frac{1}{2}}) - p^k(z_{k+\frac{1}{2}})\right)\\
\end{align}
with similar expressions also in the $y$- and $z$-momentum equations, respectively.

The only real changes with respect to the previous formulation arise in the integration of the time derivatives that, after integrating by parts in time and introducing the known solution at time $t^n$
(\textit{upwinding in time}, according to the causality principle), read  
\begin{align}
\!\!\!\!\! \int \limits_{T_{i+\frac{1}{2},j,k} \times T_{n+1}} \!\!\!\!\! \psi_{m_1}(x)\omega_{m_2}(y)\omega_{m_3}(z)\omega_{m_0}(t)  \frac{\partial u_{h}}{\partial t} d\x dt % = & \Delta x\, \Delta y\, \Delta z\,  \Mass{xyz}  \cdot u_{ll'l''l'''} \nonumber \\
%
%= &  \Delta x\, \Delta y\, \Delta z\,  \Mass{xyz} \left(\int \limits_{ T_{n+1}} \omega_{m'''}(t)\omega_{l'''}'(t) \right) \cdot u_{ll'l''l'''}
= 
\int \limits_{T_{i+\frac{1}{2},j,k} } \!\! \psi_{m_1}(x)\omega_{m_2}(y)\omega_{m_3}(z) \left( \omega_{m_0}(t^{n+1}) u_{h}(\x,t^{n+1}) - \omega_{m_0}(t^{n}) u_{h}(\x,t^{n}) \right) d\x 
\nonumber \\
- \int \limits_{T_{i+\frac{1}{2},j,k} \times T_{n+1}} \!\!\!\!\! \psi_{m_1}(x)\omega_{m_2}(y)\omega_{m_3}(z) \left( \frac{\partial}{\partial t} \omega_{m_0}(t) \right) u_{h} d\x dt
=    \Mass{xyz}  \left[ \left(  \Mass{t}_1 -   \Mass{t}_V 
\right)  \cdot    \dof{U}^{n+1}_{i,j,k}  -  \Mass{t}_0   \cdot \dof{U}^{n}_{i,j,k}   % u_{ll'l''l'''}^{n} %
 \right] \times \Delta x\, \Delta y\, \Delta z\,
\end{align}
with analogous terms for $v_h$ and $w_h$, and where
\begin{align}
\Mass{t}_1 &= \left\{ M^t_{1\,pq}\right\}_{p,q=0,..,M} = \Bigg\{ \varphi_p(1)  \varphi_q(1)  \Bigg\}_{p,q=0,..,M}, \label{eq:tmatrix1}\\
\Mass{t}_0 &= \left\{ M^t_{0\,pq}\right\}_{p,q=0,..,M} = \Bigg\{ \varphi_p(0)  \varphi_q(1)  \Bigg\}_{p,q=0,..,M}, \label{eq:tmatrix2}\\
\Mass{t}_V &= \left\{ M^t_{V\,pq}\right\}_{p,q=0,..,M} = \left\{\int \limits_0^1   \varphi'_p(\xi)  \varphi_q(\xi) \,d\xi \right\}_{p,q=0,..,M}. \label{eq:tmatrix3} 
\end{align}
With the aim of simplifying the notation, one can extend the spatial $xyz$-formalism of the tensor products (\ref{eq:tensorP1}-\ref{eq:tensorP2}) to the space-time case by defining 
a generic vector of space-time degrees of freedom as 
\begin{align}
\dof{X} = \left\{\hat{x}_{mm'm''m'''}\right\}_{m,m',m''=0,..,N; \, m'''=0,..,M},
\end{align}
the space-time operators in the form of 
\begin{align}
\begin{array}{rl} 
\operator{Z}{x$_i$ x$_j$}  &  \equiv \operator{Z}{x$_i$}\cdot\operator{Z}{x$_j$},\\
\operator{Z}{x$_i$ x$_j$ x$_k$}  &\equiv \operator{Z}{x$_i$}\cdot\operator{Z}{x$_j$}\cdot\operator{Z}{x$_k$},\\
\operator{Z}{x$_i$ x$_j$ x$_k$ x$_l$}  &\equiv \operator{Z}{x$_i$}\cdot\operator{Z}{x$_j$}\cdot\operator{Z}{x$_k$}\cdot\operator{Z}{x$_l$}, \end{array}
 \;\;\; \text{with}\; x_i,x_j,x_k,x_l \in \left\{x,y,z,t\right\} 
\end{align}
which operate along a generic vector of degrees of freedom via the tensor products
\begin{align}
\operator{Z}{x}\cdot \dof{X} = Z_{ml}I_{m'l'}I_{m''l''}I_{m'''l'''}x_{l\,l'l''l'''}, \qquad 
\operator{Z}{y}\cdot \dof{X} = I_{ml}Z_{m'l'}I_{m''l''}I_{m'''l'''}x_{l\,l'l''l'''}, \\ 
\operator{Z}{z} \cdot\dof{X} = I_{ml}I_{m'l'}Z_{m''l''}I_{m'''l'''}x_{l\,l'l''l'''}, \qquad 
\operator{Z}{t} \cdot\dof{X} = I_{ml}I_{m'l'}I_{m''l''}Z_{m'''l'''}x_{l\,l'l''l'''},  
\end{align}
where $\mathbf{Z}$ is a real square matrix, $I$ is the identity operator and the Einstein convention of summation over repeated indexes is assumed. In this notation the mass matrix that corresponds 
to the time coordinate can be written as $\Mass{t}$, according to the definition of the mass matrix in equations (\ref{eq:matrices}). 
Then, from equations (\ref{eq:uDGNSmomST}-\ref{eq:DGNSincST}) the following system is obtained: 
\begin{align}
\Mass{xyz}  \left[ \left(  \Mass{t}_1 -   \Mass{t}_V 
\right) \cdot\dof{U}^{n+1}_{i+\frac{1}{2},j,k}- \Mass{t}_0 \cdot \dof{Fu}^{n}_{i+\frac{1}{2},j,k}  \right] + \frac{\Delta t}{\Delta x}\Mass{yzt}\cdot \left( \tild{R}{x} \cdot\dof{P}^{n+1}_{i+1,j,k} - \tild{L}{x} \cdot\dof{P}^{n+1}_{i,j,k} \right)&=0,\label{eq:uSDG1_ST} \\
\Mass{xyz}   \left[ \left(  \Mass{t}_1 -   \Mass{t}_V 
\right) \cdot \dof{V}^{n+1}_{i,j+\frac{1}{2},k}- \Mass{t}_0 \cdot \dof{Fv}^{n}_{i,j+\frac{1}{2},k}  \right] + \frac{\Delta t}{\Delta y}\Mass{zxt} \cdot\left( \tild{R}{y} \cdot\dof{P}^{n+1}_{i,j+1,k} - \tild{L}{y}\cdot \dof{P}^{n+1}_{i,j,k} \right)&=0,\label{eq:vSDG1_ST} \\
\Mass{xyz}   \left[ \left(  \Mass{t}_1 -   \Mass{t}_V 
\right) \cdot \dof{W}^{n+1}_{i,j,k+\frac{1}{2}}- \Mass{t}_0 \cdot \dof{Fw}^{n}_{i,j,k+\frac{1}{2}} \right] + \frac{\Delta t}{\Delta z} \Mass{xyt}\cdot\left( \tild{R}{z}\cdot \dof{P}^{n+1}_{i,j,k+1} - \tild{L}{z} \cdot\dof{P}^{n+1}_{i,j,k} \right)&=0,\label{eq:wSDG1_ST} %\\
\end{align}
\begin{equation} 
 \frac{\Mass{yzt} \left(\barr{R}{x} \cdot\dof{U}^{n+1}_{i+\frac{1}{2},j,k} \!-\! \barr{L}{x} \cdot\dof{U}^{n+1}_{i-\frac{1}{2},j,k}\right)}{\Delta x} + 
 \frac{\Mass{zxt} \left(\barr{R}{y} \cdot\dof{V}^{n+1}_{i,j+\frac{1}{2},k} \!-\! \barr{L}{y}\cdot \dof{V}^{n+1}_{i,j-\frac{1}{2},k}\right)}{\Delta y} + 
 \frac{\Mass{xyt} \left(\barr{R}{z} \cdot\dof{W}^{n+1}_{i,j,k+\frac{1}{2}} \!-\! \barr{L}{z} \cdot\dof{W}^{n+1}_{i,j,k-\frac{1}{2}}\right)}{\Delta z} = 0, \label{eq:pSDG1_ST}
\end{equation}
which is analogous to the system of equations (\ref{eq:uSDG1}-\ref{eq:pSDG1}), where 
now the advective-diffusive terms are computed according to
\begin{align}
 \Mass{t}_0  \cdot \dof{Fu}^{n}_{i,j,k} &= \Mass{t}_0  \cdot \dof{U}^n_{i,j,k} - \frac{\Delta t}{\Delta x \Delta y \Delta z}\left(\Mass{xyz}\right)^{-1}\Mass{t} \cdot\left(\, \int \limits_{\partial T_{i,j,k}  } \omega \mathbf{F}_u\cdot \vec{n} dS  - \int \limits_{T_{i,j,k} } \nabla \omega \cdot \mathbf{F}_u \, d\x  \right) \label{eq:Fu_ST},\\
\Mass{t}_0  \cdot \dof{Fv}^{n}_{i,j,k} &= \Mass{t}_0  \cdot \dof{V}^n_{i,j,k} - \frac{\Delta t}{\Delta x \Delta y \Delta z}\left(\Mass{xyz}\right)^{-1}\Mass{t} \cdot\left(\, \int \limits_{ \partial T_{i,j,k}  } \omega \mathbf{F}_v\cdot \vec{n} dS  - \int \limits_{T_{i,j,k}  } \nabla \omega \cdot \mathbf{F}_v \, d\x  \right) \label{eq:Fv_ST}, \\
\Mass{t}_0  \cdot \dof{Fw}^{n}_{i,j,k} &= \Mass{t}_0  \cdot \dof{W}^n_{i,j,k} - \frac{\Delta t}{\Delta x \Delta y \Delta z}\left(\Mass{xyz}\right)^{-1}\Mass{t} \cdot\left(\, \int \limits_{ \partial T_{i,j,k}  } \omega \mathbf{F}_w\cdot \vec{n} dS  - \int \limits_{T_{i,j,k}  } \nabla \omega \cdot \mathbf{F}_w \, d\x  \right)  \label{eq:Fw_ST}.
\end{align}
The adopted numerical strategy for the implicit diffusion is actually the higher order time extension of the aforementioned implicit approach and it  will be described later in this section. Following the philosophy of section \ref{sec:LO}, after multiplying equations (\ref{eq:uSDG1_ST}-\ref{eq:wSDG1_ST}) by the inverse of the matrix $\Mass{xyz}  \left(  \Mass{t}_1 -   \Mass{t}_V 
\right)$, the following direct definitions of the degrees of freedom of the velocity components are obtained
\begin{align}
\dof{U}^{n+1}_{i+\frac{1}{2},j,k}& =   \left(  \Mass{t}_1 -   \Mass{t}_V 
\right)^{-1}   \Mass{t}_0 \cdot \dof{Fu}^{n}_{i+\frac{1}{2},j,k}   \reds{-} \frac{\Delta t}{\Delta x}\left(  \Mass{t}_1 -   \Mass{t}_V 
\right)^{-1} \Mass{t} \left(\Mass{x}\right)^{-1} \cdot \left( \tild{R}{x} \cdot\dof{P}^{n+1}_{i+1,j,k} - \tild{L}{x} \cdot\dof{P}^{n+1}_{i,j,k} \right),\label{eq:uSDG1_ST2} \\
 \dof{V}^{n+1}_{i,j+\frac{1}{2},k}&=  \left(  \Mass{t}_1 -   \Mass{t}_V 
\right)^{-1}  \Mass{t}_0 \cdot \dof{Fv}^{n}_{i,j+\frac{1}{2},k}  \reds{-}   \frac{\Delta t}{\Delta y}  \left(  \Mass{t}_1 -   \Mass{t}_V 
\right)^{-1}\Mass{t} \left(\Mass{y}\right)^{-1} \cdot\left( \tild{R}{y} \cdot\dof{P}^{n+1}_{i,j+1,k} - \tild{L}{y}\cdot \dof{P}^{n+1}_{i,j,k} \right),\label{eq:vSDG1_ST2} \\
 \dof{W}^{n+1}_{i,j,k+\frac{1}{2}}&=   \left(  \Mass{t}_1 -   \Mass{t}_V 
\right)^{-1}  \Mass{t}_0 \cdot \dof{Fw}^{n}_{i,j,k+\frac{1}{2}} \reds{-}   \frac{\Delta t}{\Delta z}   \left(  \Mass{t}_1 -   \Mass{t}_V 
\right)^{-1} \Mass{t} \left(\Mass{z}\right)^{-1} \cdot\left( \tild{R}{z}\cdot \dof{P}^{n+1}_{i,j,k+1} - \tild{L}{z} \cdot\dof{P}^{n+1}_{i,j,k} \right),\label{eq:wSDG1_ST2} 
\end{align}

Then, after  substitution of the resulting equations in the discrete incompressibility condition (\ref{eq:pSDG1_ST}), one obtains
\begin{align}
\left[ \Mass{t} \left(  \Mass{t}_1 -   \Mass{t}_V 
\right)^{-1} \Mass{t}\;\; \mathbb{H}^{i,j,k} \right] \cdot \mathpzc{P}^{n+1} = \left[ \Mass{t} \left(  \Mass{t}_1 -   \Mass{t}_V 
\right)^{-1} \Mass{t}_0\right] \cdot  \mathpzc{b}^n_{i,j,k} \label{eq:pSyst_ST} %\\
%\text{for}\;\;\;i=2,...,N_x-1;\;\;\;j=2,...,N_y-1;\;\;\,k=2,...,N_z-1 \nonumber
\end{align}
that is the higher order time-accurate version of the pressure equation, analogous to (\ref{eq:pSyst}). The right hand side $\mathpzc{b}^n_{i,j,k}$ collects all the known terms, 
i.e. the advective and diffusive terms $\{\dof{Fu}^{n}\}$,  $\{\dof{Fv}^{n}\}$ and  $\{\dof{Fw}^{n}\}$. This system is not symmetric because of the non-symmetric time-matrices 
(\ref{eq:tmatrix1}) and (\ref{eq:tmatrix3}).  After multiplication by the inverse of $\Mass{t} \left(  \Mass{t}_1 -   \Mass{t}_V \right)^{-1} $, the non-symmetric contribution 
of the time-matrix can be \textit{removed}, and the same well suited coefficient matrix $\mathbb{H}$ of section \ref{sec:LO} is obtained, i.e. 
\begin{align}
 &\Mass{t} \mathbb{H}^{i,j,k} \cdot \mathpzc{P}^{n+1} = \left[ \Mass{t}_0\right] \cdot  \mathpzc{b}^n_{i,j,k} \label{eq:pSyst_ST2} %\\
%&\text{for}\;\;\;i=2,...,N_x-1;\;\;\;j=2,...,N_y-1;\;\;\,k=2,...,N_z-1 \nonumber
\end{align}
and consequently, the resulting system is \textbf{symmetric} and \textbf{strictly positive definite} (for appropriate pressure boundary conditions). Hence, it can be solved very efficiently by 
means of a classical conjugate gradient method. Once the system for the higher order accurate space-time expansion coefficients of the pressure $\mathpzc{P}^{n+1}$ has been solved, 
the velocity can be readily updated according to equations (\ref{eq:uSDG1_ST2}-\ref{eq:wSDG1_ST2}). Note, however, that although the presented space-time DG framework is formally 
high order accurate in time, the final numerical scheme is strongly influenced by the time-splitting between advection, diffusion and incompressibility condition, which constrains 
the final method to be only first order accurate in time. In section \ref{subsec:Picard} a very simple numerical procedure based on the Picard iteration is outlined in order to 
circumvent the order limitation induced by the time-splitting and to enable the final solution to preserve the original high-order time accuracy of the presented spectral staggered 
space-time DG discretization.

\subsection{Implicit diffusion}
Following the same procedure outlined in Section \ref{sec:LO}, the high-order time accurate version of the implicit scheme for diffusion (\ref{eq:uSDG12}) reads
\begin{align}
&\Mass{xyz} \cdot\left[\left(  \Mass{t}_1 -   \Mass{t}_V 
\right) \cdot  \dof{U}^{n+1}_{i,j,k}- \Mass{t}_0 \cdot \dof{Fu}^{n}_{i,j,k}  \right] - \frac{\Delta t}{\Delta x} \Mass{yzt}\cdot\left(\barr{R}{x} \cdot\dof{U}^{\textsl{(x)}}_{i+\frac{1}{2},j,k} - \barr{L}{x} \cdot\dof{U}^{\textsl{(x)}}_{i-\frac{1}{2},j,k}\right) +\nonumber\\
&- \frac{\Delta t}{\Delta y}\Mass{zxt}\cdot\left(\barr{R}{y} \cdot\dof{U}^{\textsl{(y)}}_{i,j+\frac{1}{2},k} - \barr{L}{y}\cdot \dof{U}^{\textsl{(y)}}_{i,j-\frac{1}{2},k}\right)    - \frac{\Delta t}{\Delta z} \Mass{xyt}\cdot\left(\barr{R}{z} \cdot\dof{U}^{\textsl{(z)}}_{i,j,k+\frac{1}{2}} - \barr{L}{z} \cdot\dof{U}^{\textsl{(z)}}_{i,j,k-\frac{1}{2}}\right) = 0,\label{eq:uSDG12_ST}
\end{align}
Then, after substituting the definitions of the velocity derivatives (\ref{eq:Ux}-\ref{eq:Uz}) the high-order accurate space-time DG version of (\ref{eq:implicitviscosity}) can be written as 
\begin{align}
\left[\left( \mathbb{M}^t_1- \mathbb{M}^t_V \right)\mathbb{M}^{xyz} + \nu \mathbb{M}^t\mathbb{H}^{xyz}\right] \cdot \mathpzc{U} = & \mathbb{M}^t_0 \mathbb{M}^{xyz}   \cdot \mathpzc{Fu}^n, \nonumber\\ %\;\;\; %\nonumber\\
\left[ \left( \mathbb{M}^t_1- \mathbb{M}^t_V \right) \mathbb{M}^{xyz} + \nu  \mathbb{M}^t \mathbb{H}^{xyz}\right] \cdot \mathpzc{V} = & \mathbb{M}^t_0\mathbb{M}^{xyz}   \cdot \mathpzc{Fv}^n,\nonumber\\ %\;\;\; %\mathpzc{b}_v \label{eq:implicitviscosity} \\
\left[\left( \mathbb{M}^t_1- \mathbb{M}^t_V \right) \mathbb{M}^{xyz} + \nu  \mathbb{M}^t \mathbb{H}^{xyz}\right] \cdot \mathpzc{W} = & \mathbb{M}^t_0\mathbb{M}^{xyz}   \cdot \mathpzc{Fw}^n,\label{eq:implicitviscosity_ST} %\mathpzc{b}_w\nonumber
\end{align}
that is \emph{non-symmetric} because of the time-matrices 
\begin{align}
\left( \mathbb{M}^t_1- \mathbb{M}^t_V \right)=\left\{\left(  \Mass{t}_1 -   \Mass{t}_V \right)_{ijk}\right\},\hspace{1cm} \mathbb{M}^t_0 =\left\{ \left. \Mass{t}_0 \right._{ijk}\right\}
\end{align}
and can be efficiently solved by means of a classical GMRES method \cite{saad:gmres}. Notice that the non-symmetric component of system (\ref{eq:implicitviscosity_ST}) can be shifted 
to the viscous terms, i.e. the second term on the left-hand-side, by multiplying the equations with the inverse of $\left( \mathbb{M}^t_1- \mathbb{M}^t_V \right)$ from the left. 
In that case, for small viscosities, the system can be seen as a non-symmetric perturbation of the inviscid case. 

%%%%%%%%%%%%%%%%%%%%%%%%%

\subsection{Space-time pressure correction algorithm}
\label{subsec:Picard}
In Section \ref{sec:LO} the final staggered semi-implicit DG scheme (\ref{eq:LOvelocityU}-\ref{eq:LOdiffusionW}) consists of two main blocks that are solved sequentially by the use of a 
fractional time-step approach. If only high order of accuracy in space is needed, such a splitting is possible. The first fractional block is described by the discrete advection-diffusion 
equations (\ref{eq:LOdiffusionU}-\ref{eq:LOdiffusionW}), which itself contains a first fractional step for the purely explicit advection and a second fractional step for the implicit 
discretization of the diffusive terms. Then, the second fractional block contains the solution of the discrete pressure Poisson equation that results from substituting the discrete
momentum equations into the discrete incompressibility condition, i.e. combining (\ref{eq:LOvelocityU}-\ref{eq:LOvelocityW}) with (\ref{eq:LOpressure}). The important fact is that 
the chosen fractional time discretization is only first order accurate. In principle, higher order schemes for fractional time-stepping or other more sophisticated techniques could 
be adopted in defiance of simplicity or generality \cite{Karniadakis91,Zang86,KimMoin85,Marcus84}. In this work a simple Picard method has been implemented. In this manner, the first 
order time-splitting approach of system (\ref{eq:LOvelocityU}-\ref{eq:LOdiffusionW}) can then be generalized to arbitrary high order of accuracy in time at the aid of the Picard 
procedure. \reds{We emphasize that at the moment we have no rigorous mathematical proof for the fact that the Picard iterations actually increase the order of accuracy by 
one per iteration. We only have numerical evidence which support this claim in the context of high order ADER schemes, see \cite{Dumbser2008}, as well as the numerical convergence 
tables shown later in this paper for a set of test cases. }    
The final version of the spectral staggered space-time DG scheme, which is written in terms of a \textit{space-time pressure correction algorithm}, reads: 
for $k=0,\ldots,M$ do  
%%
%\begin{align}
%\left[\left( \mathbb{M}^t_1- \mathbb{M}^t_V \right)\mathbb{M}^{xyz} + \nu \mathbb{M}^t \mathbb{H}^{xyz}\right] \cdot \mathpzc{U}^{n+\frac{1}{2},(k+1)}_{\#} = &\mathbb{M}^{xyz}  \cdot \left[ \mathbb{M}^t_0 \cdot \mathpzc{Fu}^{n+1,(k)}_{\#}  %+ \right.\nonumber \\
%%& \left. 
%- \frac{\Delta t}{\Delta x} \mathbb{M}^t \cdot \left(\mathbb{M}^{-1}  \mathcal{D}\right)^x \cdot \mathpzc{P}^{n+1,(k)}_{{\#}}  \right], \label{eq:LOdiffusionU_HO}\\
%\left[\left( \mathbb{M}^t_1- \mathbb{M}^t_V \right)\mathbb{M}^{xyz} + \nu \mathbb{M}^t \mathbb{H}^{xyz}\right] \cdot \mathpzc{V}^{n+\frac{1}{2},(k+1)}_{\#} = & \mathbb{M}^{xyz}  \cdot \left[ \mathbb{M}^t_0 \cdot \mathpzc{Fv}^{n+1,(k)}_{\#} %+ \right.\nonumber \\
%%& \left. 
%- \frac{\Delta t}{\Delta y}   \mathbb{M}^t \cdot \left(\mathbb{M}^{-1}  \mathcal{D}\right)^y \cdot \mathpzc{P}^{n+1,(k)}_{{\#}} \right], \label{eq:LOdiffusionV_HO_HO}\\
%\left[\left( \mathbb{M}^t_1- \mathbb{M}^t_V \right)\mathbb{M}^{xyz} + \nu \mathbb{M}^t \mathbb{H}^{xyz}\right] \cdot \mathpzc{W}^{n+\frac{1}{2},(k+1)}_{\#} = & \mathbb{M}^{xyz}  \cdot \left[ \mathbb{M}^t_0 \cdot\mathpzc{Fw}^{n+1,(k)}_{\#} %+ \right.\nonumber \\
%%& \left.
%- \frac{\Delta t}{\Delta z}  \mathbb{M}^t \cdot \left(\mathbb{M}^{-1}  \mathcal{D}\right)^z \cdot \mathpzc{P}^{n+1,(k)}_{{\#}}  \right], \label{eq:LOdiffusionW_HO}
%\end{align}
\begin{align}
\left[\left( \mathbb{M}^t_1- \mathbb{M}^t_V \right)\mathbb{M}^{xyz} + \nu \mathbb{M}^t \mathbb{H}^{xyz}\right] \cdot \mathpzc{U}_{\#} ^{\reds{ n+1,(k+\frac{1}{2})}}= &\mathbb{M}^{xyz}  \cdot \left[ \mathbb{M}^t_0 \cdot \mathpzc{Fu}^{n+1,(k)}_{\#}  %+ \right.\nonumber \\
%& \left. 
- \reds{ \pi_{:*\rightarrow \#}}\left\{\frac{\Delta t}{\Delta x} \mathbb{M}^t \cdot \left(\mathbb{M}^{-1}  \mathcal{D}\right)^x \cdot \mathpzc{P}^{n+1,(k)}_{{\#}} \right\} \right], \label{eq:LOdiffusionU_HO}\\
\left[\left( \mathbb{M}^t_1- \mathbb{M}^t_V \right)\mathbb{M}^{xyz} + \nu \mathbb{M}^t \mathbb{H}^{xyz}\right] \cdot \mathpzc{V}_{\#} ^{\reds{ n+1,(k+\frac{1}{2})}}  = & \mathbb{M}^{xyz}  \cdot \left[ \mathbb{M}^t_0 \cdot \mathpzc{Fv}^{n+1,(k)}_{\#} %+ \right.\nonumber \\
%& \left. 
-  \reds{ \pi_{:*\rightarrow \#}}\left\{\frac{\Delta t}{\Delta y}   \mathbb{M}^t \cdot \left(\mathbb{M}^{-1}  \mathcal{D}\right)^y \cdot \mathpzc{P}^{n+1,(k)}_{{\#}} \right\}\right], \label{eq:LOdiffusionV_HO_HO}\\
\left[\left( \mathbb{M}^t_1- \mathbb{M}^t_V \right)\mathbb{M}^{xyz} + \nu \mathbb{M}^t \mathbb{H}^{xyz}\right] \cdot \mathpzc{W}_{\#} ^{\reds{ n+1,(k+\frac{1}{2})}}  = & \mathbb{M}^{xyz}  \cdot \left[ \mathbb{M}^t_0 \cdot\mathpzc{Fw}^{n+1,(k)}_{\#} %+ \right.\nonumber \\
%& \left.
-  \reds{ \pi_{:*\rightarrow \#} }\left\{\frac{\Delta t}{\Delta z}  \mathbb{M}^t \cdot \left(\mathbb{M}^{-1}  \mathcal{D}\right)^z \cdot \mathpzc{P}^{n+1,(k)}_{{\#}} \right\} \right], \label{eq:LOdiffusionW_HO}
\end{align}
%
%\begin{gather}
%\mathbb{M}^t\mathbb{H}^{xyz}\cdot \left( \mathpzc{P}^{n+1,(k+1)}_{\#} - \mathpzc{P}^{n+1,(k)}_{\#} \right) = \mathbb{M}^t_0 \mathpzc{b}^{n+\frac{1}{2},(k)}_{\#}, \label{eq:LOpressure_HO}\\
%\text{with} \;\;\; \mathpzc{b}^{n+\frac{1}{2},(k)}_{\#}\equiv\mathbb{M}^{yz} \mathcal{D}^{T\,x}  \cdot \mathpzc{U}^{n+\frac{1}{2},(k+1)}_{*}    +
%\mathbb{M}^{zx} \mathcal{D}^{T\,y} \cdot \mathpzc{V}^{n+\frac{1}{2},(k+1)}_{*} + \mathbb{M}^{xy} \mathcal{D}^{T\,z} \cdot \mathpzc{W}^{n+\frac{1}{2},(k+1)}_{*}, \nonumber
%\end{gather}
%
\begin{gather}
\mathbb{M}^t\mathbb{H}^{xyz}\cdot \left( \mathpzc{P}^{n+1,(k+1)}_{\#} - \mathpzc{P}^{n+1,(k)}_{\#} \right) =  \reds{ \mathpzc{b}_{\#}}^{\reds{n+1,(k+\frac{1}{2})}}, \label{eq:LOpressure_HO}\\
\text{with} \;\;\; \mathpzc{b}_{\#}^{\reds{ n+1,(k+\frac{1}{2})}}\equiv\mathbb{M}^{yz} \mathcal{D}^{T\,x}  \cdot \mathpzc{U}_{*}^{\reds{ n+1,(k+\frac{1}{2})}}     +
\mathbb{M}^{zx} \mathcal{D}^{T\,y} \cdot \mathpzc{V}_{*}^{\reds{ n+1,(k+\frac{1}{2})}} + \mathbb{M}^{xy} \mathcal{D}^{T\,z} \cdot \mathpzc{W}_{*}^{\reds{ n+1,(k+\frac{1}{2})}} , \nonumber
\end{gather}
\begin{align}
\mathpzc{U}^{n+1,(k+1)}_{*}& =  \left( \mathbb{M}^t_1- \mathbb{M}^t_V \right)^{-1} \cdot \left[ \mathbb{M}^t_0 \cdot\mathpzc{U}_{*}^{\reds{ n+1,(k+\frac{1}{2})}}   - \frac{\Delta t}{\Delta x}  \mathbb{M}^t \cdot \left(\mathbb{M}^{-1}  \mathcal{D}\right)^x \cdot \left( \mathpzc{P}^{n+1,(k+1)}_{{\#}} - \mathpzc{P}^{n+1,(k)}_{{\#}}  \right)\right], \label{eq:LOvelocityU_HO}\\
\mathpzc{V}^{n+1,(k+1)}_{*} &= \left( \mathbb{M}^t_1- \mathbb{M}^t_V \right)^{-1}\cdot\left[ \mathbb{M}^t_0\cdot \mathpzc{V}_{*}^{\reds{ n+1,(k+\frac{1}{2})}} - \frac{\Delta t}{\Delta y}  \mathbb{M}^t\cdot   \left(\mathbb{M}^{-1}  \mathcal{D}\right)^y \cdot  \left(\mathpzc{P}^{n+1,(k+1)}_{{\#}} - \mathpzc{P}^{n+1,(k)}_{{\#}}   \right)\right],   \label{eq:LOvelocityV_HO} \\
\mathpzc{W}^{n+1,(k+1)}_{*}& = \left( \mathbb{M}^t_1- \mathbb{M}^t_V \right)^{-1}\cdot\left[ \mathbb{M}^t_0 \cdot\mathpzc{W}_{*}^{\reds{ n+1,(k+\frac{1}{2})}}   - \frac{\Delta t}{\Delta z}  \mathbb{M}^t\cdot \left(\mathbb{M}^{-1}  \mathcal{D}\right)^z \cdot  \left( \mathpzc{P}^{n+1,(k+1)}_{{\#}} - \mathpzc{P}^{n+1,(k)}_{{\#}}   \right)\right],    \label{eq:LOvelocityW_HO} 
\end{align}
%
%\begin{gather}
%\text{for} \;\;\, k=0,\ldots,M \;\;\; \nonumber
%\end{gather}
% 
where $M$ is the maximum degree of the time-polynomials and $k$ is the Picard iteration number. Note that the Picard process allows to gain one order of accuracy in time 
per Picard iteration \reds{when applied to an ODE}, see \cite{Layton2008,minion2003,minion2001}. 
%\begin{align}
%\mathpzc{Fu}^{n,(0)} = \mathpzc{Fu}^{n}, \;\;\;\mathpzc{Fv}^{n,(0)} = \mathpzc{Fv}^{n}, \;\;\; \mathpzc{Fw}^{n,(0)} = \mathpzc{Fw}^{n},
$\mathpzc{Fu}=\left\{\dof{Fu}_{i,j,k}\right\}$,  $\mathpzc{Fv}=\left\{ \dof{Fv}_{i,j,k}\right\}$,  $\mathpzc{Fw}=\left\{ \dof{Fw}_{i,j,k}\right\}$ 
%\end{align}
are the advective terms computed according to \eqref{eq:Fu_ST}-\eqref{eq:Fw_ST}, without taking into account the diffusive flux. 
 %to the velocity components $\mathpzc{U}^{n+1,(k)}$, $\mathpzc{V}^{n+1,(k)}$ and $\mathpzc{W}^{n+1,(k)}$, i.e.
%\begin{align}
 %\Mass{t}_0  \cdot \dof{Fu}^{n+1,(k)}_{i,j,k} &= \Mass{t}_0  \cdot \dof{U}^n_{i,j,k} - \frac{\Delta t}{\Delta x \Delta y \Delta z}\iMass{xyz}\Mass{t} \cdot\left(\, \int \limits_{\partial T_{ijk}} \omega \mathbf{F}_u^{n+1,(k)} \cdot \vec{n} dS- \int \limits_{T_{ijk}} \nabla \omega \cdot \mathbf{F}_u^{n+1,(k)}  \, dxdydz\right)\label{eq:Fu_ST1},\\
%%
 %\Mass{t}_0  \cdot \dof{Fv}^{n+1,(k)}_{i,j,k} &= \Mass{t}_0  \cdot \dof{V}^n_{i,j,k} - \frac{\Delta t}{\Delta x \Delta y \Delta z}\iMass{xyz}\Mass{t} \cdot\left(\, \int \limits_{ \partial T_{ijk}} \omega \mathbf{F}_v^{n+1,(k)} \cdot \vec{n} dS - \int \limits_{T_{ijk}} \nabla \omega \cdot \mathbf{F}_v^{n+1,(k)}  \, dxdydz\right)^{n+1,(k)} \label{eq:Fv_ST1}, \\
%%
%\Mass{t}_0  \cdot \dof{Fw}^{n+1,(k)}_{i,j,k} &= \Mass{t}_0  \cdot \dof{W}^n_{i,j,k} - \frac{\Delta t}{\Delta x \Delta y \Delta z}\iMass{xyz}\Mass{t} \cdot\left(\, \int \limits_{ \partial T_{ijk}} \omega \mathbf{F}_w^{n+1,(k)} \cdot \vec{n} dS- \int \limits_{T_{ijk}} \nabla \omega \cdot \mathbf{F}_w^{n+1,(k)}  \, dxdydz\right)^{n+1,(k)}  \label{eq:Fw_ST1}.
%\end{align}  with %of the previous time step $t\in T^n$.
We furthermore set 
\begin{align}
\mathpzc{Fu}^{n+1,(0)} = \mathpzc{Fu}^{n}, \;\;\;\mathpzc{Fv}^{n+1,(0)} = \mathpzc{Fv}^{n}, \;\;\; \mathpzc{Fw}^{n+1,(0)} = \mathpzc{Fw}^{n},
%$\mathpzc{Fu}^{n+1,(0)} = \mathpzc{Fu}^{n}$,  $\mathpzc{Fv}^{n+1,(0)} = \mathpzc{Fv}^{n+1}$,  $\mathpzc{Fw}^{n,(0)} = \mathpzc{Fw}^{n}$,
\end{align}
and $\mathpzc{P}^{n+1,(k)}$ is the $k$-th iterate for the discrete pressure, for which we use the trivial initial guess 
\begin{align}
\mathpzc{P}^{n+1,(0)}_{{\#}} = 0.
%$\mathpzc{Fu}^{n+1,(0)} = \mathpzc{Fu}^{n}$,  $\mathpzc{Fv}^{n+1,(0)} = \mathpzc{Fv}^{n+1}$,  $\mathpzc{Fw}^{n,(0)} = \mathpzc{Fw}^{n}$,
\end{align}
%are the advective terms (\ref{eq:Fu_ST}-\ref{eq:Fw_ST}) computed accordingly to the velocity components $\mathpzc{U}^{n,(M)}$, $\mathpzc{V}^{n,(M)}$ and $\mathpzc{W}^{n,(M)}$ of the previous time step $t\in T^n$.
Thanks to the Picard procedure the desired properties of the presented spectral space-time DG method are re-established, so that the final algorithm  (\ref{eq:LOvelocityU_HO})-(\ref{eq:LOdiffusionW_HO}) is arbitrary high-order accurate both in space and time. Finally, it is important to stress that the proposed iterative solution of the non-trivial system of equations (\ref{eq:LOvelocityU_HO})-(\ref{eq:LOdiffusionW_HO}) is feasible in practice,  thanks to the fact that the coefficient matrix $\mathbb{H}$ that enters into the discrete Poisson equation (i.e. the incompressibility condition) and the discrete diffusion equation is well conditioned and can be solved in a very efficient way via modern matrix-free Krylov subspace methods, even \textit{without the use of any preconditioner}.  %In fact,  the resolution of the linear systems is very efficient 
Finally, note that when the degree of the time-polynomials $M$ is set to be zero, then 
\begin{align}
\Mass{t}_1 \equiv 1, \;\;\; \Mass{t}_0 \equiv 1, \;\;\; \Mass{t}_V \equiv 0, \nonumber
\end{align}
and the method collapses to the previous spectral staggered semi-implicit DG scheme with a classical first order backward Euler discretization in time. 
Moreover, if at the same time the spatial and the temporal polynomial approximation degrees are chosen to be zero ($M=N=0$), then the following equalities arise from (\ref{eq:matrices}) 
\begin{align}
\Mass{} \equiv 1, \;\;\; \tild{R}{}  \equiv 1, \;\;\; \tild{L}{}   \equiv 1,  \;\;\; \barr{R}{}   = 1, \;\;\; \barr{L}{} \equiv 1, \nonumber
\end{align}
and the method collapses to a classical staggered semi-implicit finite-difference finite-volume method for the incompressible Navier-Stokes equations, where the pressure field is defined 
at the barycenters of the main grid and the velocity components are defined at the middle points of the cell interfaces, i.e. the classical family of efficient semi-implicit 
methods on staggered grids of Casulli et al. \cite{CasulliDumbserToro,BoscheriDumbser,CasulliWalters,Casulli1999,CasulliCheng,Casulli1990,Blood3D2014,DumbserIbenIoriatti,Casulli2009,CasulliZanolli2012,CasulliZanolli2002,CasulliCattani,CasulliStelling2011} is obtained.

%\begin{itemize}
%\item iterative nonlinear methods: \cite{BrugnanoCasulli2009,Casulli2009,CasulliZanolli2012,CasulliZanolli2002}. 
%
%\item semi-implicit staggered grids \cite{CasulliDumbserToro,BoscheriDumbser,CasulliWalters,Casulli1999,CasulliStelling1998,Casulli1995,CasulliCheng,Casulli1990,Blood3D2014,DumbserIbenIoriatti}
%\item stability efficiency \cite{CasulliCattani}
%\item purely applicative articles SI staggered \cite{CasulliStelling1996}
%\item subgrid \cite{CasulliStelling2011}
%\item semi-impl higher order DG \cite{DumbserCasulli2013,TavelliDumbser2014,TavelliDumbser2014b,TavelliDumbser2015}
%\end{itemize}

\subsection{Numerical validation}
In this section the capabilities of our new spectral space-time DG method are tested against several numerical benchmark problems in two and three space dimensions for which either an 
analytical or other numerical reference solutions exist. In particular, three different numerical convergence tables are produced, with the aim of assuring that the presented method 
is really arbitrary high-order accurate in both space \textit{and} time. Note that achieving high order time accuracy for the incompressible Navier-Stokes equations is far from being
straightforward. % because this may be true in theory, but in principle not in practice.

%Here we present the convergence table \ref{tab:Womersley} for the Womersley oscillatory flow between two flat plates for assessing that the final numerical solution is really '\emph{arbitrary'} higher order accurate in space and time.

\subsubsection{Oscillatory viscous flow between two flat plates}

\begin{figure} 
\centering %subfloat
			\includegraphics[width=0.5\textwidth]{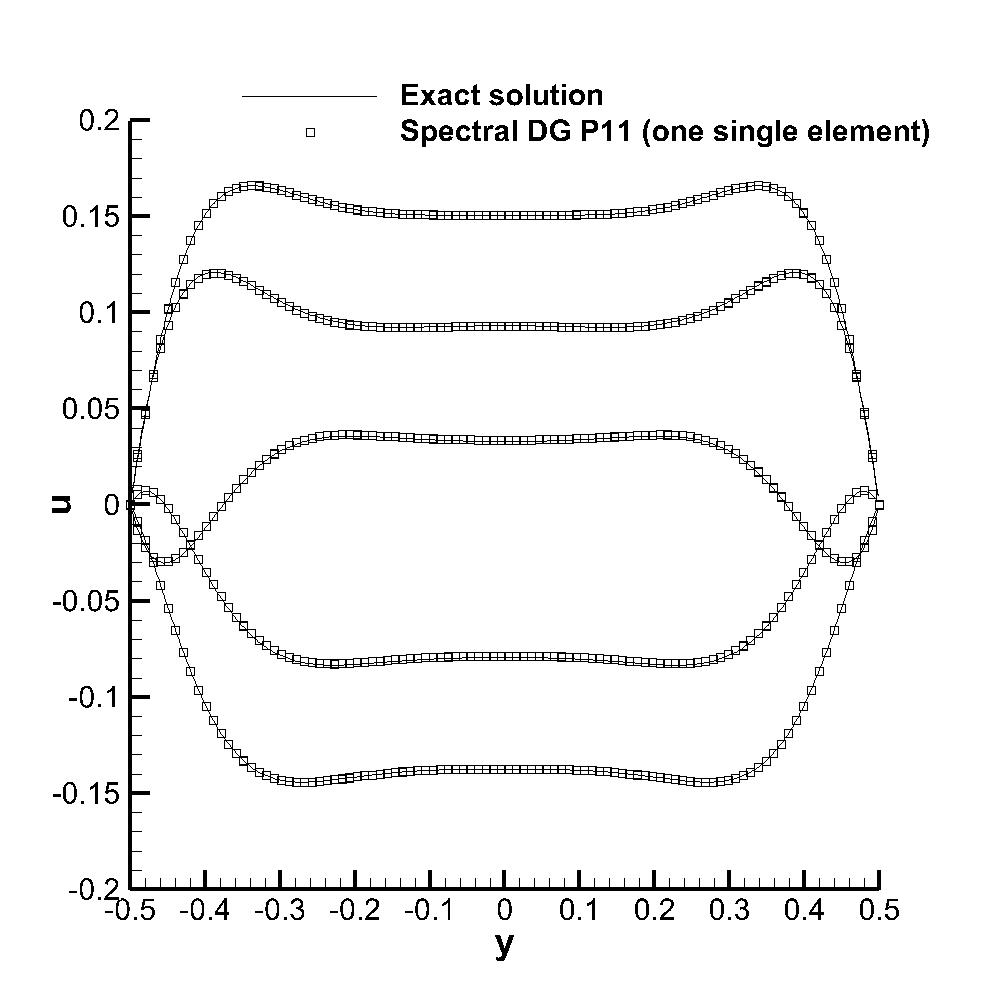}
\caption{The numerical solution interpolated along 100 equidistant spatial points obtained for the unsteady Womersley problem compared with the exact solution \cite{Womersley,Loudon1998} 
at different times for $\nu = 2\cdot 10^{-2}$:  $t=1.8$, $1.6$, $2.0$, $1.4$, $2.2$,  respectively,  from the bottom to the top. A staggered spectral space-time DG-$\mathbb{P}_{11}$ 
method has been run using only \textit{one single} space-time element.}\label{fig:Womersley_ST}
\end{figure}

\begin{figure} 
\centering %subfloat
			\includegraphics[width=\textwidth]{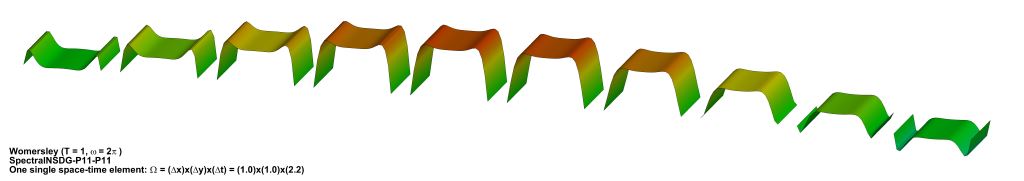}\\
			\includegraphics[width=\textwidth]{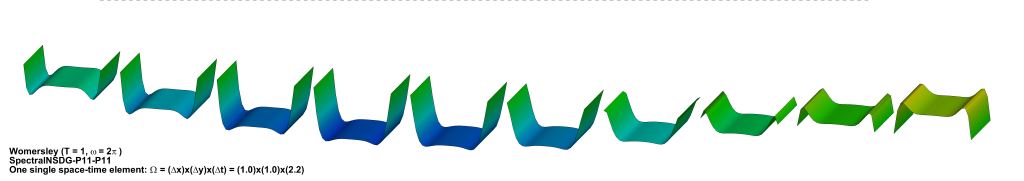}\\
			\includegraphics[width=\textwidth]{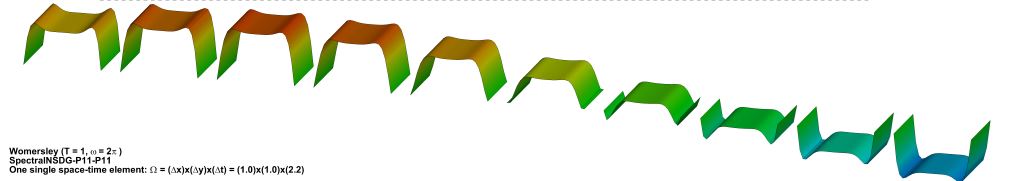}\\
			\includegraphics[width=\textwidth]{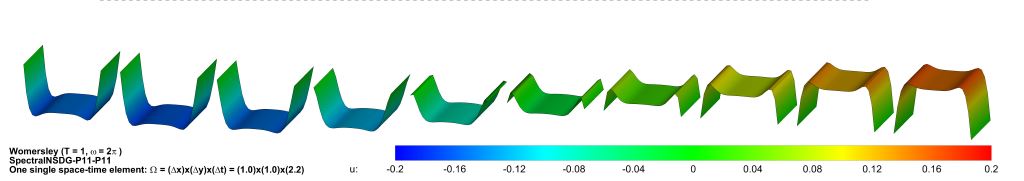}
\caption{The numerical solution obtained for the two-dimensional oscillatory flow between two flat plates with the staggered spectral space-time DG-$\mathbb{P}_{11}$ scheme. 
The computational domain in space and time $\Omega=\Delta x \times \Delta y \times \Delta t = 1.0 \times 1.0 \times 2.2$ has been discretized by using only \textbf{\emph{one single space time element}}, and the plotted numerical solution for the velocity field has been interpolated along $40$ time slices with $t\in [0,2.2]$, respectively, from top left to bottom right.}\label{fig:Wom_ST_clips}
\end{figure}

In this test, the fluid flow between two parallel flat plates is driven by a time harmonic pressure gradient.  According to \cite{Landau-Lifshitz6,Kurzweg85,Loudon1998}, by neglecting 
the nonlinear convective terms, the resulting axial velocity profile is only a function of time and the distance from the plates. The flow furthermore depends only on one single 
dimensionless parameter, known as the Womersley number $\alpha_W = R\sqrt{\omega/\nu}$, see \cite{Womersley}, where $R$ is the half distance between the two plates, 
$\omega$ is the frequency of the oscillations and $\nu$ is the kinematic viscosity. In particular the fluid velocity and pressure are given by 
\begin{align}
&u\left(x,y\right) = \frac{A}{i\,\omega}\left[1 - \frac{\cosh \left( \alpha_W \sqrt{i} \left.y\middle/R\right.\right)}{\cosh\left( \alpha_W \sqrt{i}\right)}  \right],\nonumber\\
&\frac{\partial p}{\partial x} = \frac{p(x_R) - p(x_L)}{L} = - A\,e^{i\omega t},\nonumber
\end{align}
where $i=\sqrt{-1}$ is the imaginary unit, $L=x_R-x_L$ is the total length of the duct and the amplitude has been chosen equal to $A=1$. 
The exact solution has been chosen as initial condition at $t=0$, then pressure conditions are 
imposed on the left and right boundaries, while no-slip boundary conditions have been imposed at the upper and lower walls. The other parameters
of this test problem were chosen as $L=1$, $R=0.5$ and $\omega=1$. 
Figs. \ref{fig:Womersley_ST} and \ref{fig:Wom_ST_clips} show
 the numerical results obtained for $\nu = 2 \cdot 10^{-2}$ with our spectral staggered space-time DG scheme using only \emph{\textbf{one single }} 
$\p_{11}$ \emph{space-time} element (M=N=11), completing the entire simulation within the time interval $t\in[0,2.2]$ in one single time step. 
The results are compared with the exact analytical solution at different intermediate output times.  
In particular, for this test problem, two periods of oscillation are resolved within a single time-step, and the complete velocity profile is resolved within a 
single spatial cell.  
%Despite the low order of time accuracy (SIDG-$\p_{11}$ is only first order),  
From the obtained results one can conclude that the proposed staggered spectral space-time DG scheme is indeed very accurate in both space and time, since it is able 
to resolve all flow features within \textit{one single} space-time element. %The chosen time step, without considering CFL condition, is $dt=10^{-3}$.

%In order to give a demonstration of the numerical potential of the presented method, this problem has been resolved with the
% $\mathbb{P}_{11}$  version of our space-time spectral DG by forcing the code to complete the simulation $t\in[0,2.2]$
%with only \textbf{a single space-time element}.
%Figure $\ref{fig:Womersley_ST}$ show the numerical solution for the axial velocity field interpolated at different 
%time slices along $100$ spatial points next to the exact solution. Figure \ref{fig:Wom_ST_clips} shows the same numerical solution interpolated along $40$ different
%time slices and it gives an idea of the quality of the information that are stored inside a single space-time polynomial which keeps inside exactly $(N+1)\cdot(M+1)=144$ degrees of freedom. In particular, two periods of oscillation are resolved within a single time-step, and the complete velocity profile is resolved within a single spatial cell.
%As a result, it is confirmed that the computed solution perfectly matches the analytical one with the correct order of accuracy for this test problem.

Furthermore, Table \ref{tab:Womersley} contains the results of a numerical convergence study that we have performed with this smooth unsteady two-dimensional flow 
problem, for which an exact solution is available. The order of accuracy has been verified up to order $7$ in space and time by evaluating the $L_2$ and $L_{\infty}$ errors 
\begin{align}
\epsilon_{L_2} = \sqrt{ \int_{\Omega} \left(u_h - u \right)^2 }, \;\;\; \text{and} \;\;\; \epsilon_{L_{\infty}} = \max\limits_{\Omega}  \left|u_h - u \right|,\nonumber 
\end{align}
at different discretization numbers for the polynomial degrees $N=M=1,\ldots,6$. From the obtained results we conclude that the designed order of accuracy of the scheme
has been reached in both space and time. \reds{For the polynomial degree $N=M=1$, only sub-optimal convergence rates have been verified experimentally, and will be 
subject of future research.}

\subsubsection{2D Taylor-Green vortex}
\begin{figure} 
\centering %subfloat
			\includegraphics[width=0.59\textwidth]{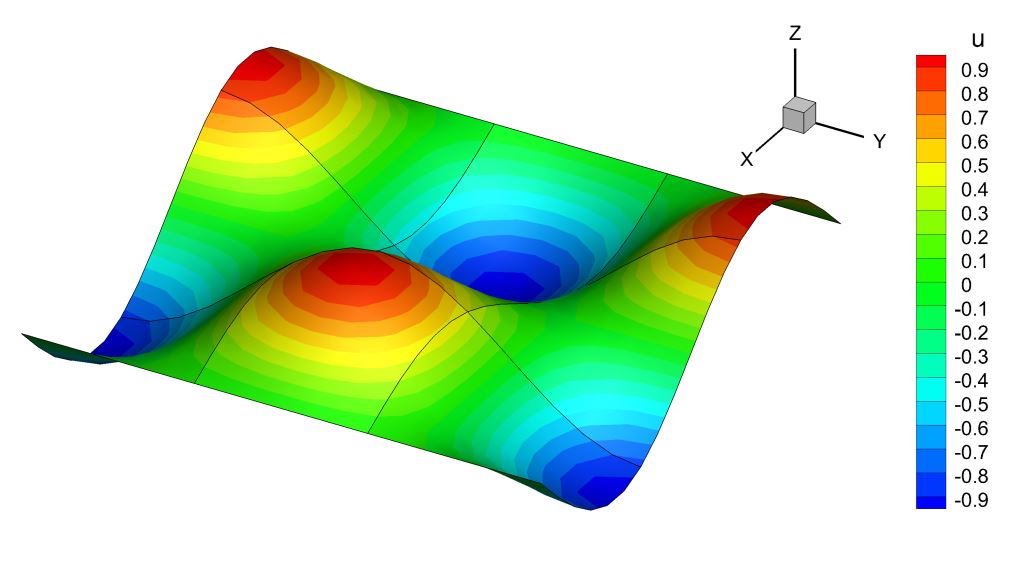}
			\includegraphics[width=0.39\textwidth]{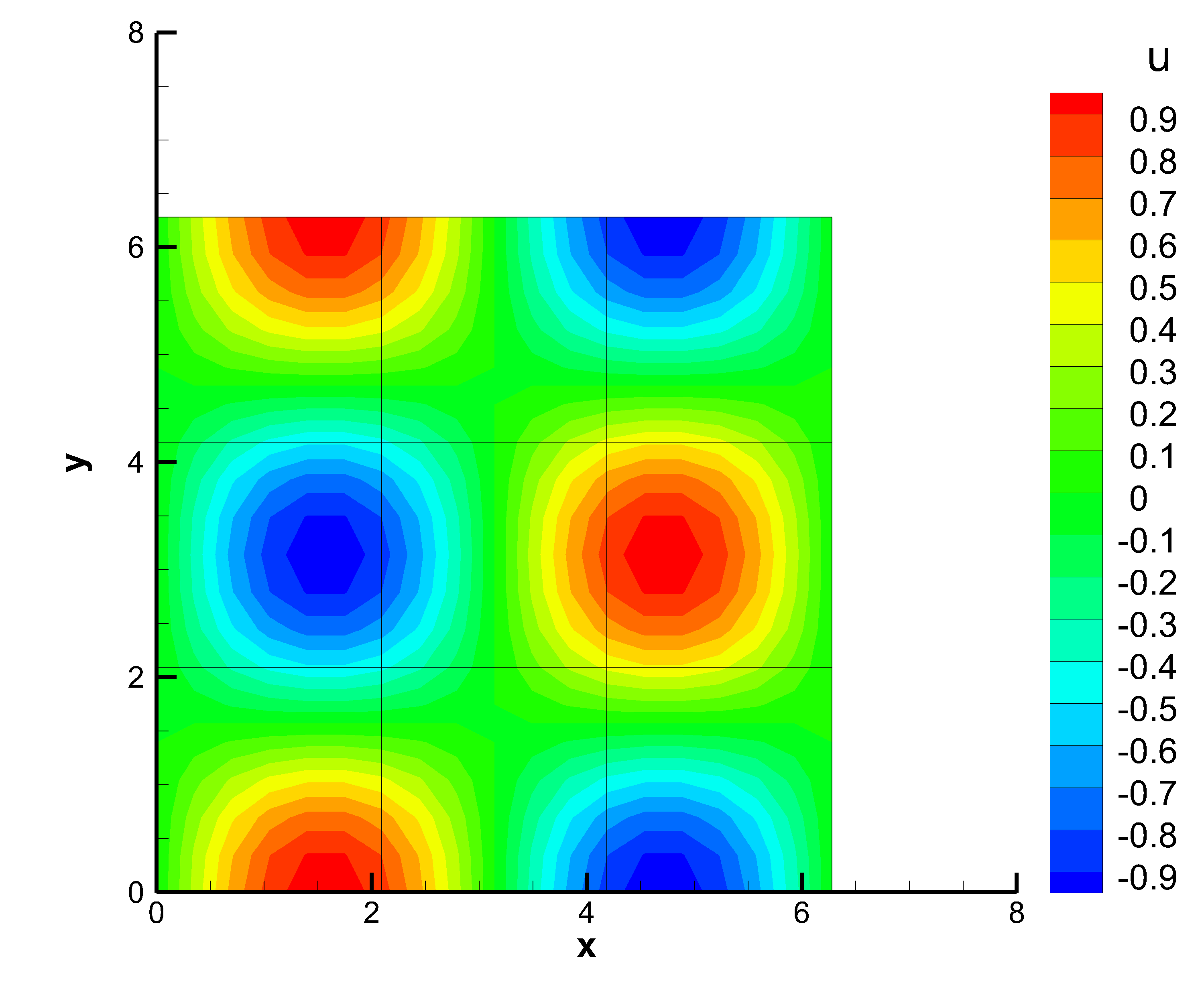}\\
			%\caption{ciao}
			\includegraphics[width=0.59\textwidth]{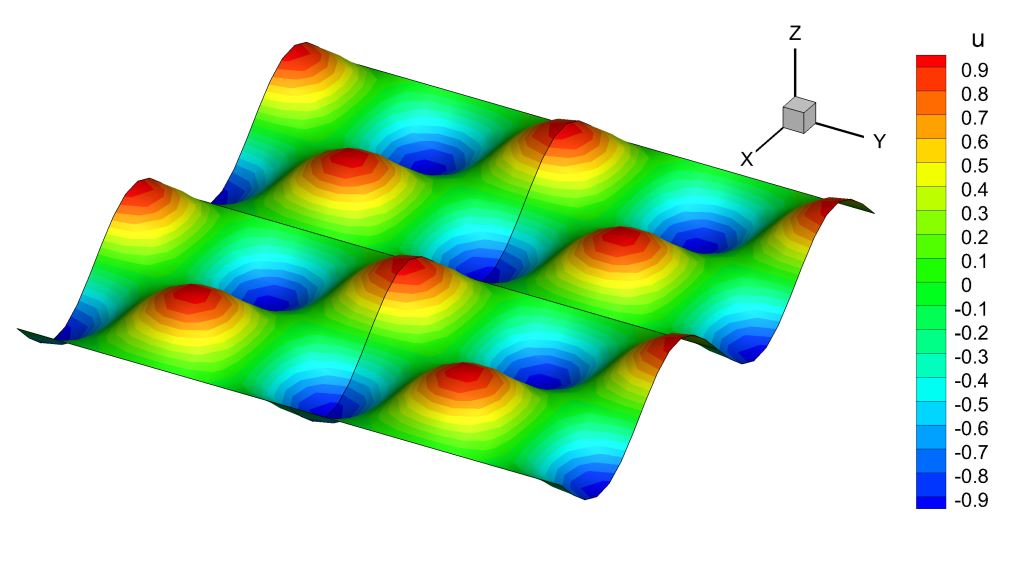} 
			\includegraphics[width=0.39\textwidth]{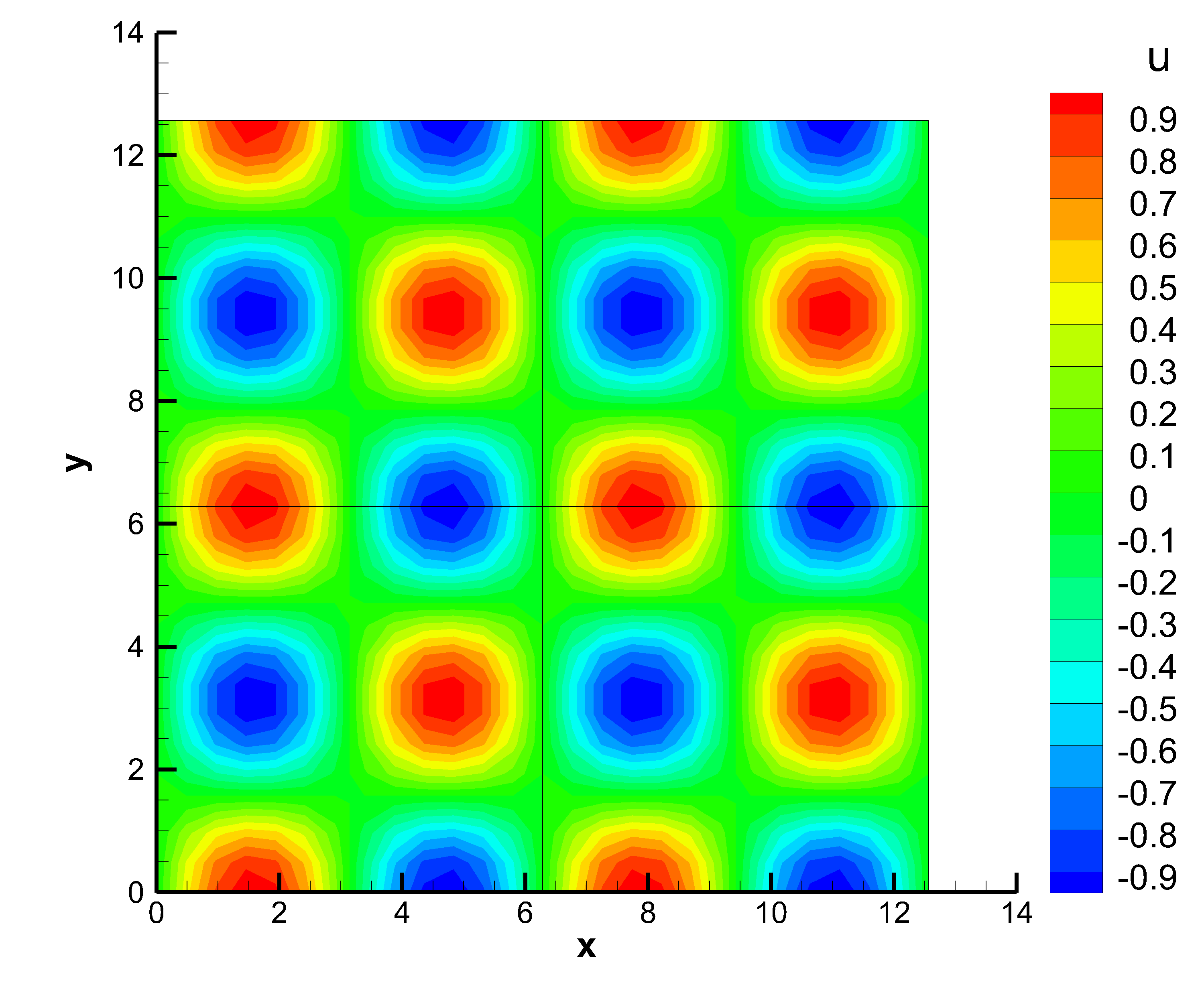}
			\caption{Numerical solution of the $u$ velocity component for the two dimensional Taylor-Green vortex problem computed with the staggered spectral space-time DG-$\mathbb{P}_{5}$ method 
			using $3^2$ elements with $L=2\pi$ (top) and the staggered spectral space-time DG-$\mathbb{P}_{12}$ scheme using $2^2$ elements with $L=4\pi$ (bottom).}\label{fig:TGV2D}
%\caption{The numerical solution obtained for the two-dimensional oscillatory flow between two flat plates with our space-time spectral SIDG-$\mathbb{P}_{11}$. The space-time domain $\Omega=\Delta x \times \Delta y \times \Delta t = 1.0 \times 1.0 \times 2.2$ has been discretized into \textbf{\emph{one single space time element}}, and the plotted numerical solution for the velocity field has been interpolated along $40$ time slices with $t\in [0,2.2]$, respectively  from the left to right, from the top to the bottom.}\label{fig:TGV2D}
\end{figure}

\begin{figure} 
\centering %subfloat
			\includegraphics[width=0.7\textwidth]{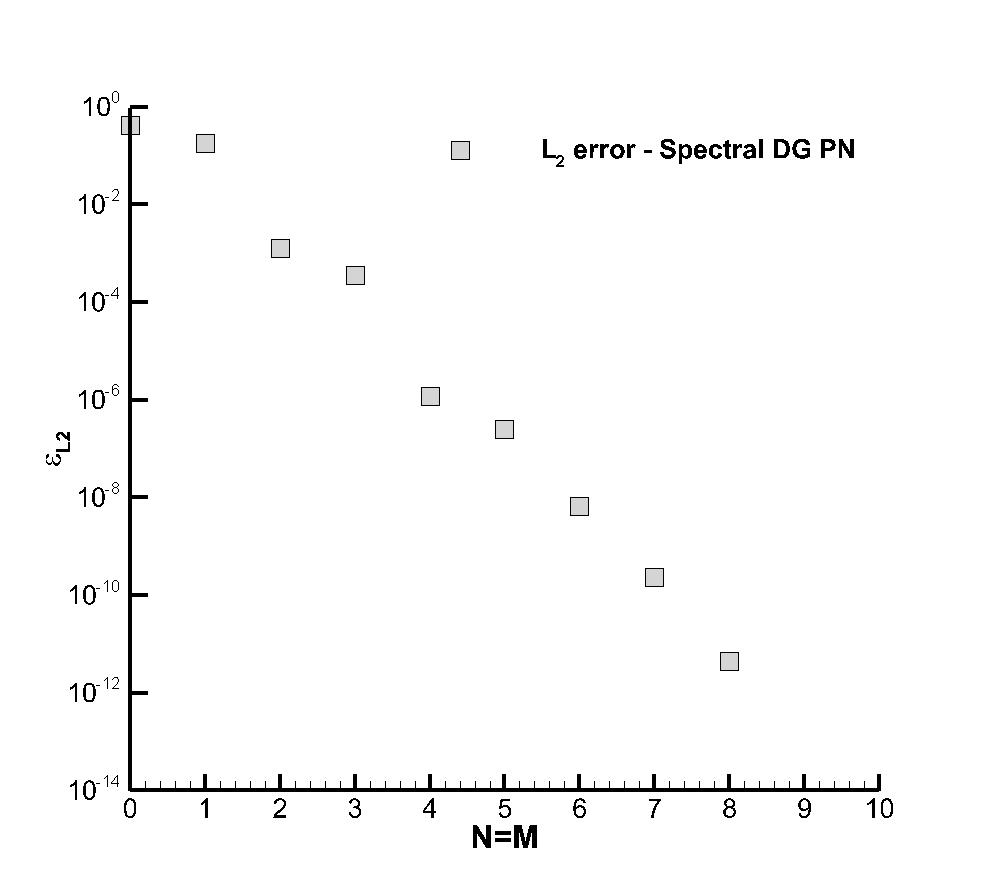}
			\caption{Numerical $L_2$ error $\epsilon_{L_2}$ of the $u$ velocity component for the two dimensional Taylor-Green vortex problem computed with staggered spectral space-time DG-$\mathbb{P}_{N}$ schemes as a function of the polynomial degree $N=M$ on a fixed grid of $12^2$ elements.}\label{fig:spectral}
%\caption{The numerical solution obtained for the two-dimensional oscillatory flow between two flat plates with our space-time spectral SIDG-$\mathbb{P}_{11}$. The space-time domain $\Omega=\Delta x \times \Delta y \times \Delta t = 1.0 \times 1.0 \times 2.2$ has been discretized into \textbf{\emph{one single space time element}}, and the plotted numerical solution for the velocity field has been interpolated along $40$ time slices with $t\in [0,2.2]$, respectively  from the left to right, from the top to the bottom.}\label{fig:TGV2D}
\end{figure}

The two dimensional Taylor-Green vortex problem is widely used for testing the accuracy of numerical schemes, because
it offers another smooth unsteady analytical solution of the incompressible Navier-Stokes equations with periodic
boundary conditions. The exact solution of this problem is given by
\begin{align}
&u(x,y,t) = \sin(x) \cos(y) e^{-2\nu t}, \qquad  v(x,y,t) = - \cos(x) \sin(y) e^{-2\nu t}, \nonumber \\
&p(x,y,t) = \frac{1}{4} \left( \cos(2x) + \cos (2y)\right) e^{-4\nu t}.\nonumber 
\end{align}
The computational domain is $\Omega=[0,L]^2$ with periodic boundary conditions. The initial sinusoidal velocity field is smoothed in time
by the viscous dissipative forces. The convergence study for this test is summarized in Table \ref{tab:TGV2D}. The accuracy of our staggered 
spectral space-time DG scheme is verified for polynomial degrees $N=M=1,\ldots,8$.  
Figure \ref{fig:TGV2D} shows the numerical solution obtained by setting $L=2\pi$ for the staggered spectral space-time DG-$\mathbb{P}_{5}$ scheme, using a very coarse mesh composed of only 
$3^2$ spatial elements. Furthermore, we repeat this test with $L=4\pi$ using a staggered spectral space-time DG-$\mathbb{P}_{12}$ scheme using only $2^2$ spatial elements. 
Moreover, Figure \ref{fig:spectral} shows the behavior of the error $\epsilon_{L_2}$ as a function of the polynomial degree ($N=M$) for a fixed mesh: 
\emph{the exponential decay} of the error, i.e. the \textbf{spectral convergence} obtained with our scheme by increasing the polynomial approximation degree in space and time, 
is explicitly verified. The results confirm the designed accuracy in space and time and show how the presented numerical method works properly even when using very high order 
approximation polynomials and very coarse meshes. \reds{Also in this two dimensional test, for the polynomial degree $N=M=1$ a non-optimal convergence has been experimetally verified.}

\subsubsection{3D Arnold-Beltrami-Childress flow}
\begin{figure} 
\centering %subfloat
			\includegraphics[width=0.4\textwidth]{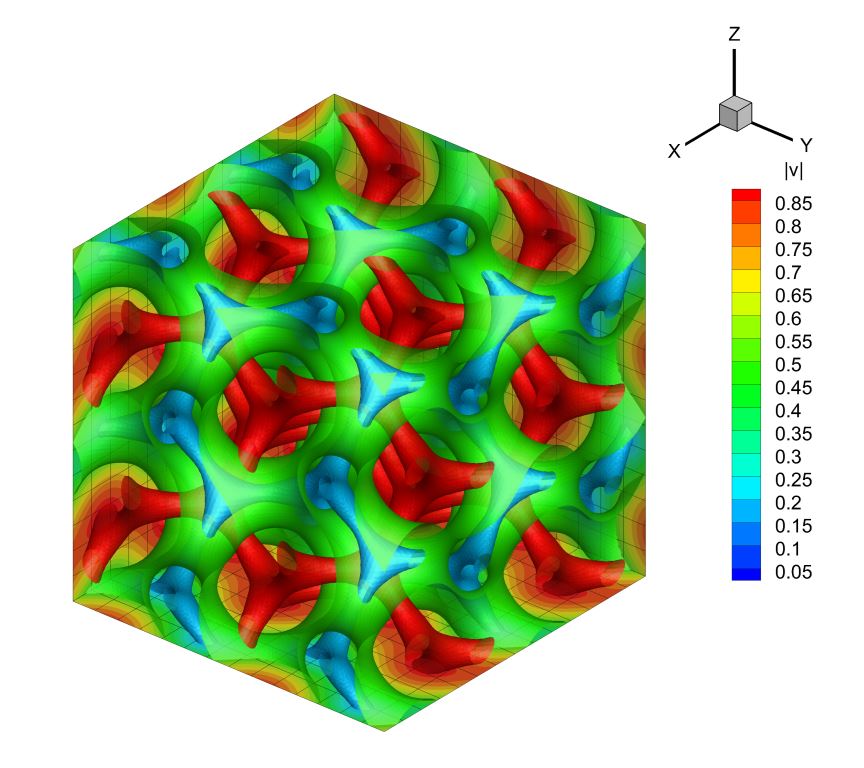}
			\includegraphics[width=0.4\textwidth]{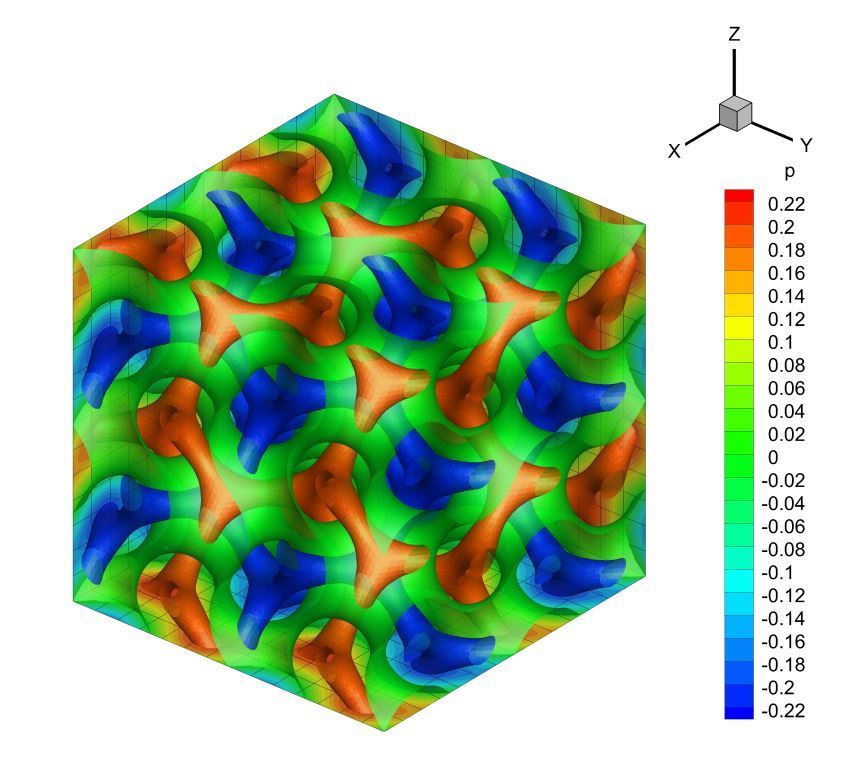}\\
			\includegraphics[width=0.6\textwidth]{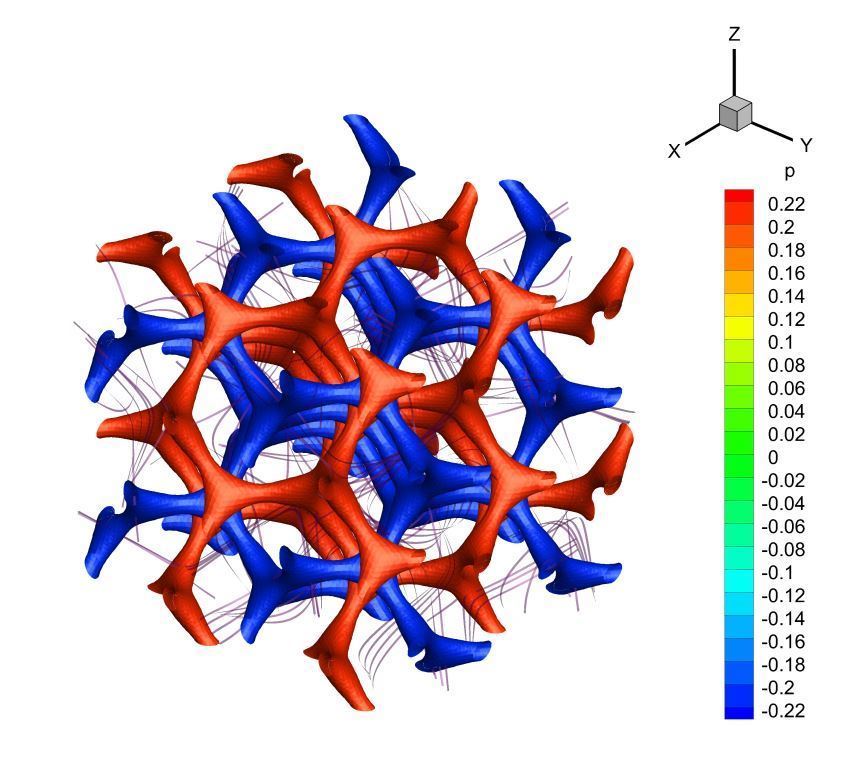}
			\caption{Numerical solution for the three dimensional Arnold-Beltrami-Childress (ABC) flow at time $t=1.0$ computed with the staggered spectral space-time DG-$\mathbb{P}_{5}$ scheme using 
			only $7^3$ elements in space. 
			The periodic solution has been replicated along the three dimensional cube of edge $L_c=4\pi$ for giving a better view of the field variables. At the top of the figure the 
			velocity is plotted on the left and the pressure is depicted on the right; at the bottom the 3D stream-traces together with the pressure isosurfaces are plotted.}\label{fig:ABC} 
\end{figure}

In order to test the accuracy of our staggered spectral space-time DG scheme also against an unsteady three dimensional benchmark problem,  
the Arnold-Beltrami-Childress (ABC) flow, proposed by Arnold \cite{Arnold65} and Childress in \cite{Childress70}, is considered. 
For this smooth unsteady test problem, the exact solution reads 
\begin{align}
&u(x,y,z,t) = \left[\sin(z) + \cos(y)\right] e^{-\nu t}, \nonumber \\
&v(x,y,z,t) = \left[\sin(x) +  \cos(z)\right] e^{-\nu t}, \label{eq:ABC} \\
&w(x,y,z,t) = \left[\sin(y) +  \sin(x)\right] e^{-\nu t}. \nonumber 
\end{align}
The computational domain is the cube $\Omega = [0,2 \pi]^3$, with periodic boundary conditions everywhere. 
Given the initial condition (\ref{eq:ABC}) at time $t=0$, the corresponding analytical solution decays exponentially in time according to the chosen kinematic viscosity.  
Also for this three dimensional time-dependent test problem, the designed high order of accuracy of our staggered spectral space-time DG scheme 
has been confirmed up to order $9$ by a numerical convergence study that is summarized in Table \ref{tab:ABC}. Similar to the two dimensional 
Taylor-Green vortex, in the 3D ABC flow the advective terms, the pressure forces and the incompressibility condition are highly coupled. 
The numerical solution for $\nu=0.1$ at time $t=10$ is depicted in Figure \ref{fig:ABC}.  
%These features make the ABC flow a very good test in validating a numerical method suitable for turbulent simulations.

 %Supported by the above numerical results, the presented new spectral-DG method actually showed the desirable resulution properties for the simulation of turbulent flows.
%\begin{align}
%%K = \left.\int_{\omega} \frac{1}{2} \mathbf{v}^2 \middle/  \int_{\Omega} 1  \right.
%%K = \frac{\int\limits_{\Omega} \frac{1}{2} \mathbf{v}^2 }{ \left\| \Omega \right\|  }
%K = \left.\int\limits_{\Omega} \frac{1}{2} \mathbf{v}^2 \middle/  \left\| \Omega \right\| \right.
%\end{align}

%--------------------------

\section{Conclusion} 
\label{sec:conclusion}

In this paper the new family of staggered spectral semi-implicit DG methods, recently proposed by Dumbser and Casulli in \cite{DumbserCasulli2013} for the shallow water equations on 
staggered Cartesian grids, has been extended to the incompressible Navier-Stokes equations in two and three space dimensions and to arbitrary high order of accuracy in time, adopting a novel staggered spectral space-time 
DG formalism. A similar formulation has been recently presented in \cite{TavelliDumbser2014,TavelliDumbser2014b,TavelliDumbser2015} for unstructured staggered meshes, but there the chosen 
staggered grid was slightly different, and the use of unstructured meshes did \textit{not} allow to produce a spectral DG scheme based on simple tensor products of one-dimensional operators. 
Of course, unstructured meshes as those used in \cite{TavelliDumbser2014,TavelliDumbser2014b,TavelliDumbser2015} allow to fit very complicate geometries and complex physical boundaries, 
however, by choosing staggered Cartesian grids, some interesting advantages follow, in particular: 
\begin{enumerate}
\item 
Cartesian grids allow the use of \emph{tensor-products} of the basis and test functions; this means that the weak formulation of the governing equations can be written as a very handy 
combination of one dimensional integrals over the canonical reference element $\xi\in[0,1]$; 
\item this fact significantly minimizes the computational costs and difficulties for evaluating 
integrals, because the defined matrices are the same for all the elements in the Cartesian framework; 
\item by using basis functions that are built from the Lagrange interpolation polynomials passing through the Gauss-Legendre quadrature points, the basis functions are \textit{orthogonal} and thus the resulting 
mass matrices are \textit{diagonal}; this fact reduces significantly the computational cost for a 
mass-matrix multiplication; 
\item in our staggered Cartesian framework, each velocity component is defined on a different staggered dual control volume; consequently, the computation of convective and viscous 
terms on the main grid by interpolating from the dual grids to the main grid and vice versa is simpler and more natural than a discretization of these terms on the dual grids;  
\item the resulting numerical method achieves a \emph{spectral convergence} property, i.e. the computational error decreases exponentially when increasing the 
degree of the approximation polynomials in space and time.
\end{enumerate} 
The key-role of the proposed mesh-staggering, combined with the adopted semi-implicit or space-time DG time discretization, is to optimize the sparsity pattern of the  %the \emph{connectivity} of 
resulting pressure system. Furthermore, the pressure system is symmetric and only block five-diagonal for the 2D case, or only block seven-diagonal for the 3D case. 
In addition, we have presented a new way of evaluating the viscous terms in the DG framework, computing the velocity gradient (i.e. the stress tensor) on the staggered dual control volumes, 
\reds{which can be interpreted as the use of a Bassi-Rebay-type lifting operator that accounts for the jumps of the solution in the discrete gradients, but on the dual grid}. 
This allows to compute the viscous terms via an implicit discretization with essentially the same coefficient matrix that is already used in the discrete pressure system,  
with an additional diagonal term that further enforces the stability of the system. 
The resulting algorithm is shown to be arbitrary higher order accurate in space and time, robust, stable, and very efficient compared to other classical higher order DG methods for the 
incompressible Navier-Stokes equations on collocated grids, which either lead to a larger computational stencil or to a larger linear system with more unknowns. 
These features have been verified against a large set of test cases in two and three space dimensions. 
The designed space-time accuracy of our method has been verified up to $9$-th order through a series of numerical convergence tests in two and three space dimensions. 

For the simulation of turbulent flows, very high spatial and temporal resolution is needed for giving a correct and complete description of the flow physics. 
%In many cases, the smaller turbulent scales are spatially localized in some specific region of the computational domain, for example next to the wall boundaries. 
With the aim of improving the efficiency of our algorithm further, future work will concern the extension of the present high order staggered DG schemes to 
\emph{space-time adaptive meshes}, following the ideas outlined in \cite{Zanotti2015,AMR3DCL,AMR3DNC,Zanotti2015c,Zanotti2015d}. 
By introducing \emph{adaptive mesh refinement} (\textbf{AMR}) \emph{for staggered grids}, simulations of turbulent flows should become feasible.  
Moreover, it is a well known fact that discontinuous Galerkin schemes suffer of spurious oscillations when attempting to resolve shocks, 
because of Gibbs phenomenon. A novel \textit{a posteriori} approach of shock capturing for DG schemes, without losing the classical subcell resolution properties of the 
DG method, has been recently proposed for collocated grids in \cite{Dumbser2014,Zanotti2015c,Zanotti2015d}. The extension of this \textit{a posteriori} subcell limiter 
techniques to semi-implicit DG schemes on staggered meshes belongs to future investigations. 
%We think that this work, together with \cite{DumbserCasulli2013,TavelliDumbser2014,TavelliDumbser2014b,TavelliDumbser2015}, can give a good contribution to the development of a new class of %\emph{efficient} arbitrary higher order methods suitable for the simulations of turbulent flows or large scale fluid flow problems.
The possibility of extending the present semi-implicit staggered DG schemes to the context of the \emph{compressible Navier-Stokes equations}, following the ideas of 
\cite{klein,munz,KleinMach,MeisterMach,DumbserCasulli2016}, is also another topic of future research.

%\cite{DumbserCasulli,2STINS,2DSIUSW,BrugnanoCasulli,BrugnanoCasulli2,BrugnanoSestini,CasulliZanolli2010,CasulliZanolli2012,DumbserIbenIoriatti}. 

%=============================================================================
%==========    A C K N O W L E D G M E N T S
\section*{Acknowledgments}
The authors would like to thank Maurizio Tavelli and Vincenzo Casulli for the inspiring discussions on the topic. 

The research presented in this paper was financed by the European Research Council (ERC) under the
European Union's Seventh Framework Programme (FP7/2007-2013) with the research project \textit{STiMulUs}, 
ERC Grant agreement no. 278267.

The authors would like to acknowledge PRACE for awarding access to the SuperMUC supercomputer based in Munich, Germany at the Leibniz Rechenzentrum (LRZ).  

Last but not least, the authors would like to thank the two referees for their helpful comments and remarks.

\bibliographystyle{plain}

\bibliography{references} %biblio}

\begin{thebibliography}{100}

\bibitem{Albensoeder2005}
S.~Albensoeder and H.C. Kuhlmann.
\newblock Accurate three-dimensional lid-driven cavity flow.
\newblock {\em Journal of Computational Physics}, 206(2):536 -- 558, 2005.

\bibitem{Arakawa}
A.~Arakawa and V.R. Lamb.
\newblock {Computational design of the basic dynamical processes of the UCLA
  general circulation model}.
\newblock {\em Methods of Computational Physics}, 17:173--265, 1977.

\bibitem{Armaly1983}
B.~F. Armaly, F.~Durst, J.~C.~F. Pereira, and B.~Sch\"onung.
\newblock Experimental and theoretical investigation of backward-facing step
  flow.
\newblock {\em Journal of Fluid Mechanics}, 127:473--496, 2 1983.

\bibitem{Arnold2001}
D.N. Arnold, F.~Brezzi, B.~Cockburn, and L.D. Marini.
\newblock Unified analysis of discontinuous {Galerkin} methods for elliptic
  problems.
\newblock {\em SIAM J. Numer. Anal.}, 39(5):1749--1779, May 2001.

\bibitem{Arnold1984}
D.N. Arnold, F.~Brezzi, and M.~Fortin.
\newblock A stable finite element for the {Stokes} equations.
\newblock {\em Calcolo}, 21(4):337--344, 1984.

\bibitem{Arnold65}
V.I. Arnold.
\newblock Sur la topologic des \'ecoulements stationnaires des fluides
  parfaits.
\newblock {\em Comptes Rendus Hebdomadaires des S\'eances de l'Acad\'emie des
  Sciences}, 261:17--20, 1965.

\bibitem{BarthCharrier}
T.~Barth and P.~Charrier.
\newblock Energy stable flux formulas for the discontinuous {Galerkin}
  discretization of first-order nonlinear conservation laws.
\newblock Technical Report NAS-01-001, NASA, 2001.

\bibitem{Bassi2007}
F.~Bassi, A.~Crivellini, D.A.~Di Pietro, and S.~Rebay.
\newblock An implicit high-order discontinuous {Galerkin} method for steady and
  unsteady incompressible flows.
\newblock {\em Computers \& Fluids}, 36(10):1529 -- 1546, 2007.

\bibitem{BassiRebay}
F.~Bassi and S.~Rebay.
\newblock A high-order accurate discontinuous finite element method for the
  numerical solution of the compressible {Navier}-{Stokes} equations.
\newblock {\em Journal of Computational Physics}, 131:267--279, 1997.

\bibitem{BaumannOden1}
C.E. Baumann and J.T. Oden.
\newblock {A discontinuous hp finite element method for convection-diffusion
  problems}.
\newblock {\em Computer Methods in Applied Mechanics and Engineering},
  175:311--341, 1999.

\bibitem{BaumannOden2}
C.E. Baumann and J.T. Oden.
\newblock {A discontinuous hp finite element method for the Euler and
  Navier-Stokes equations}.
\newblock {\em International Journal for Numerical Methods in Fluids},
  31:79--95, 1999.

\bibitem{Becketal}
A.D. Beck, T.~Bolemann, D.~Flad, H.~Frank, G.J. Gassner, F.~Hindenlang, and
  C.-D. Munz.
\newblock High-order discontinuous {Galerkin} spectral element methods for
  transitional and turbulent flow simulations.
\newblock {\em International Journal for Numerical Methods in Fluids},
  76:522--548, 2014.

\bibitem{Bermudez1998}
A.~Bermudez, A.~Dervieux, J.A. Desideri, and M.E.~V\'azquez Cend\'on.
\newblock Upwind schemes for the two--dimensional shallow water equations with
  variable depth using unstructured meshes.
\newblock {\em Computer Methods in Applied Mechanics and Engineering},
  155:49--72, 1998.

\bibitem{Bermudez2014}
A.~Berm\'udez, J.L. Ferr\'in, L.~Saavedra, and M.E.~V\'azquez Cend\'on.
\newblock {A projection hybrid finite volume/element method for low-Mach number
  flows}.
\newblock {\em Journal of Computational Physics}, 271:360--378, 2014.

\bibitem{BoscheriDumbser}
W.~Boscheri, M.~Dumbser, and M.~Righetti.
\newblock A semi-implicit scheme for 3d free surface flows with high order
  velocity reconstruction on unstructured {Voronoi} meshes.
\newblock {\em International Journal for Numerical Methods in Fluids}, 72:607--
  631, 2013.

\bibitem{Brachet1983}
M.~E. Brachet, D.~I. Meiron, S.~A. Orszag, B.~G. Nickel, R.~H. Morf, and
  U.~Frisch.
\newblock Small-scale structure of the {Taylor-Green} vortex.
\newblock {\em Journal of Fluid Mechanics}, 130:411--452, 5 1983.

\bibitem{Brezzi1989}
F.~Brezzi, C.~Canuto, and A.~Russo.
\newblock A self-adaptive formulation for the {Euler/Navier-Stokes} coupling.
\newblock {\em Computer Methods in Applied Mechanics and Engineering},
  73(3):317--330, 1989.

\bibitem{Brooks1982}
A.~N. Brooks and T.~J.~R. Hughes.
\newblock Streamline upwind/{Petrov-Galerkin} formulations for convection
  dominated flows with particular emphasis on the incompressible navier-stokes
  equations.
\newblock {\em Computer Methods in Applied Mechanics and Engineering},
  32(1-3):199 -- 259, 1982.

\bibitem{BrugnanoCasulli}
L.~Brugnano and V.~Casulli.
\newblock Iterative solution of piecewise linear systems.
\newblock {\em SIAM Journal on Scientific Computing}, 30:463--472, 2008.

\bibitem{BrugnanoCasulli2}
L.~Brugnano and V.~Casulli.
\newblock Iterative solution of piecewise linear systems and applications to
  flows in porous media.
\newblock {\em SIAM Journal on Scientific Computing}, 31:1858--1873, 2009.

\bibitem{BrugnanoSestini}
L.~Brugnano and A.~Sestini.
\newblock Iterative solution of piecewise linear systems for the numerical
  solution of obstacle problems.
\newblock {\em Journal of Numerical Analysis, Industrial and Applied
  Mathematics}, 6:67--82, 2012.

\bibitem{Canuto1985}
C.~Canuto, S.I. Hariharan, and L.~Lustman.
\newblock Spectral methods for exterior elliptic problems.
\newblock {\em Numerische Mathematik}, 46(4):505--520, 1985.

\bibitem{Canuto1984}
C.~Canuto, Y.~Maday, and A.~Quarteroni.
\newblock Combined finite element and spectral approximation of the
  {Navier-Stokes} equations.
\newblock {\em Numerische Mathematik}, 44(2):201--217, 1984.

\bibitem{Canuto1998}
C.~Canuto, A.~Russo, and V.~Van~Kemenade.
\newblock Stabilized spectral methods for the {Navier-Stokes} equations:
  Residual-free bubbles and preconditioning.
\newblock {\em Computer Methods in Applied Mechanics and Engineering},
  166(1-2):65--83, 1998.

\bibitem{Canuto1996}
C.~Canuto and V.~Van~Kemenade.
\newblock Bubble-stabilized spectral methods for the incompressible
  {Navier-Stokes} equations.
\newblock {\em Computer Methods in Applied Mechanics and Engineering},
  135(1-2):35--61, 1996.

\bibitem{Casulli1990}
V.~Casulli.
\newblock {S}emi-implicit finite difference methods for the two-dimensional
  shallow water equations.
\newblock {\em J. Comp. Phys.}, 86:56--74, 1990.

\bibitem{Casulli1999}
V.~Casulli.
\newblock {A} semi-implicit finite difference method for non-hydrostatic,
  free-surface flows.
\newblock {\em Int. J. Numeric. Meth. Fluids}, 30:425--440, 1999.

\bibitem{Casulli2009}
V.~Casulli.
\newblock A high-resolution wetting and drying algorithm for free-surface
  hydrodynamics.
\newblock {\em International Journal for Numerical Methods in Fluids},
  60:391--408, 2009.

\bibitem{CasulliVOF}
V.~Casulli.
\newblock {A semi--implicit numerical method for the free--surface
  Navier--Stokes equations}.
\newblock {\em International Journal for Numerical Methods in Fluids},
  74:605--622, 2014.

\bibitem{CasulliCattani}
V.~Casulli and E.~Cattani.
\newblock {S}tability, accuracy and efficiency of a semi implicit method for
  three-dimensional shallow water flow.
\newblock {\em Comp. Math. Appl.}, 27:99--112, 1994.

\bibitem{CasulliCheng}
V.~Casulli and R.~T. Cheng.
\newblock {S}emi-implicit finite difference methods for three-dimensional
  shallow water flow.
\newblock {\em Int. J. Numeric. Meth. Fluids}, 15:629--648, 1992.

\bibitem{CasulliDumbserToro}
V.~Casulli, M.~Dumbser, and E.F. Toro.
\newblock {S}emi-implicit numerical modeling of axially symmetric flows in
  compliant arterial systems.
\newblock {\em Int. J. Numeric. Meth. Biomed. Engng.}, 28:257--272, 2012.

\bibitem{CasulliCompressible}
V.~Casulli and D.~Greenspan.
\newblock Pressure method for the numerical solution of transient, compressible
  fluid flows.
\newblock {\em International Journal for Numerical Methods in Fluids},
  4(11):1001--1012, 1984.

\bibitem{CasulliStelling2011}
V.~Casulli and G.~S. Stelling.
\newblock Semi-implicit subgrid modelling of three-dimensional free-surface
  flows.
\newblock {\em International Journal for Numerical Methods in Fluids},
  67:441--449, 2011.

\bibitem{CasulliWalters}
V.~Casulli and R.~A. Walters.
\newblock {A}n unstructured grid, three-dimensional model based on the shallow
  water equations.
\newblock {\em Int. J. Numeric. Meth. Fluids}, 32:331--348, 2000.

\bibitem{CasulliZanolli2002}
V.~Casulli and P.~Zanolli.
\newblock {S}emi-implicit numerical modeling of nonhydrostatic free-surface
  flows for environmental problems.
\newblock {\em Math. Comp. Model.}, 36:1131--1149, 2002.

\bibitem{CasulliZanolli2010}
V.~Casulli and P.~Zanolli.
\newblock A nested newton--type algorithm for finite volume methods solving
  {Richards'} equation in mixed form.
\newblock {\em SIAM Journal on Scientific Computing}, 32:2255--2273, 2009.

\bibitem{CasulliZanolli2012}
V.~Casulli and P.~Zanolli.
\newblock {I}terative solutions of mildly nonlinear systems.
\newblock {\em J. Comp. Appl. Math.}, 236:3937--3947, 2012.

\bibitem{ChungNS}
S.W. Cheung, E.~Chung, H.H. Kim, and Y.~Qian.
\newblock {Staggered discontinuous Galerkin methods for the incompressible
  Navier–Stokes equations}.
\newblock {\em Journal of Computational Physics}, 302:251--266, 2015.

\bibitem{Childress70}
S.~Childress.
\newblock New solutions of the kinematic dynamo problem.
\newblock {\em Journal of Mathematical Physics}, 11:3063--3076, 1970.

\bibitem{chorin1967}
A.J. Chorin.
\newblock A numerical method for solving incompressible viscous flow problems.
\newblock {\em Journal of Computational Physics}, 2(1):12--26, 1967.

\bibitem{chorin1968}
A.J. Chorin.
\newblock Numerical solution of the {Navier-Stokes} equations.
\newblock {\em Mathematics of Computation}, 22(104):745--762, 1968.

\bibitem{StaggeredDG2}
E.T. Chung, P.~Ciarlet, and T.F. Yu.
\newblock {Convergence and superconvergence of staggered discontinuous Galerkin
  methods for the three--dimensional Maxwell's equations on Cartesian grids}.
\newblock {\em Journal of Computational Physics}, 235:14--31, 2013.

\bibitem{StaggeredDG3}
E.T. Chung, H.H. Kim, and O.B. Widlund.
\newblock {Two--level overlapping Schwarz algorithms for a staggered
  discontinuous Galerkin method}.
\newblock {\em SIAM Journal on Numerical Analysis}, 51:47--67, 2013.

\bibitem{chung2012staggered}
E.T. Chung and C.S. Lee.
\newblock A staggered discontinuous {Galerkin} method for the
  convection--diffusion equation.
\newblock {\em Journal of Numerical Mathematics}, 20(1):1--32, 2012.

\bibitem{Cockburn1990}
B.~{Cockburn}, {S.} {How}, and {C.-W.} {Shu}.
\newblock {TVB Runge Kutta Local Projection Discontinuous {Galerkin} Finite
  Element Method for Conservation Laws IV: The Multidimensional Case}.
\newblock {\em Math. Comp.}, 54:545, 1990.

\bibitem{cockburn_2000_dg}
B.~Cockburn, G.~E. Karniadakis, and C.-W. Shu.
\newblock {\em Discontinuous Galerkin Methods: Theory, Computation and
  Applications}.
\newblock Lacture Notes on Computational Science and Engineering. Springer,
  2000.

\bibitem{Cockburn1989b}
B.~{Cockburn}, {S.-Y.} {Lin}, and {C.-W.} {Shu}.
\newblock {TVB Runge Kutta Local Projection Discontinuous {Galerkin} Finite
  Element Method for Conservation Laws III: One-Dimensional Systems}.
\newblock {\em Journal of Computational Physics}, 84:90, September 1989.

\bibitem{Cockburn1989a}
B.~{Cockburn} and {C.-W.} {Shu}.
\newblock {TVB Runge Kutta Local Projection Discontinuous {Galerkin} Finite
  Element Method for Scalar Conservation Laws II: General Framework}.
\newblock {\em Math. Comp.}, 52:411, 1989.

\bibitem{CockburnShu1998}
B.~Cockburn and C.W. Shu.
\newblock The local discontinuous {Galerkin} method for time-dependent
  convection-diffusion systems.
\newblock {\em SIAM Journal on Numerical Analysis}, 35(6):2440--2463, 1998.

\bibitem{CockburnShu98}
B.~Cockburn and C.W. Shu.
\newblock The {Runge--Kutta} discontinuous {Galerkin} method for conservation
  laws {V}: multidimensional systems.
\newblock {\em Journal of Computational Physics}, 141(2):199--224, 1998.

\bibitem{cockburn_2001_rkd}
B.~Cockburn and C.W. Shu.
\newblock {R}unge-{K}utta discontinuous {G}alerkin methods for
  convection-dominated problems.
\newblock {\em Journal of Scientific Computing}, 16(3):173, 2001.

\bibitem{Crivellini2013}
A.~Crivellini, V.~D'Alessandro, and F.~Bassi.
\newblock {High-order discontinuous Galerkin solutions of three-dimensional
  incompressible RANS equations}.
\newblock {\em Computers and Fluids}, 81:122--133, 2013.

\bibitem{Dolejsi2008}
V.~Dolejsi.
\newblock Semi-implicit interior penalty discontinuous {Galerkin} method for
  viscous compressible flows.
\newblock {\em Communications in Computational Physics}, 4(2):231--274, 2008.

\bibitem{Dolejsi2004}
V.~Dolejsi and M.~Feistauer.
\newblock A semi-implicit discontinuous {Galerkin} finite element method for
  the numerical solution of inviscid compressible flow.
\newblock {\em Journal of Computational Physics}, 198(2):727 -- 746, 2004.

\bibitem{Dolejsi2007}
V.~Dolejsi, M.~Feistauer, and J.~Hozman.
\newblock Analysis of semi-implicit {DGFEM} for nonlinear convection-diffusion
  problems on nonconforming meshes.
\newblock {\em Computer Methods in Applied Mechanics and Engineering},
  196(29-30):2813 -- 2827, 2007.

\bibitem{DumbserNSE}
M.~Dumbser.
\newblock Arbitrary high order {PNPM} schemes on unstructured meshes for the
  compressible {Navier--Stokes} equations.
\newblock {\em Computers \& Fluids}, 39:60--76, 2010.

\bibitem{Dumbser2008}
M.~{Dumbser}, D.~S. {Balsara}, E.~F. {Toro}, and C.-D. {Munz}.
\newblock {A unified framework for the construction of one-step finite volume
  and discontinuous Galerkin schemes on unstructured meshes}.
\newblock {\em Journal of Computational Physics}, 227:8209--8253, September
  2008.

\bibitem{DumbserCasulli2013}
M.~Dumbser and V.~Casulli.
\newblock A staggered semi-implicit spectral discontinuous {Galerkin} scheme
  for the shallow water equations.
\newblock {\em Applied Mathematics and Computation}, 219(15):8057 -- 8077,
  2013.

\bibitem{DumbserCasulli2016}
M.~Dumbser and V.~Casulli.
\newblock A conservative, weakly nonlinear semi-implicit finite volume scheme
  for the compressible {Navier--Stokes} equations with general equation of
  state.
\newblock {\em Applied Mathematics and Computation}, 272, Part 2:479 -- 497,
  2016.

\bibitem{AMR3DNC}
M.~Dumbser, A.~Hidalgo, and O.~Zanotti.
\newblock {High Order Space-Time Adaptive ADER-WENO Finite Volume Schemes for
  Non-Conservative Hyperbolic Systems}.
\newblock {\em Computer Methods in Applied Mechanics and Engineering},
  268:359--387, 2014.

\bibitem{DumbserIbenIoriatti}
M.~Dumbser, U.~Iben, and M.~Ioriatti.
\newblock An efficient semi-implicit finite volume method for axially symmetric
  compressible flows in compliant tubes.
\newblock {\em Applied Numerical Mathematics}, 89:24 -- 44, 2015.

\bibitem{dumbser_jsc}
M.~Dumbser and {C.D.} Munz.
\newblock Building blocks for arbitrary high order discontinuous {Galerkin}
  schemes.
\newblock {\em Journal of Scientific Computing}, 27:215--230, 2006.

\bibitem{HPRmodel}
M.~Dumbser, I.~Peshkov, and E.~Romenski.
\newblock {High order ADER schemes for a unified first order hyperbolic
  formulation of continuum mechanics: Viscous heat-conducting fluids and
  elastic solids}.
\newblock {\em Journal of Computational Physics}, 314:824--862, 2016.

\bibitem{AMR3DCL}
M.~Dumbser, O.~Zanotti, A.~Hidalgo, and D.S. Balsara.
\newblock {ADER-WENO Finite Volume Schemes with Space-Time Adaptive Mesh
  Refinement}.
\newblock {\em Journal of Computational Physics}, 248:257--286, 2013.

\bibitem{Dumbser2014}
M.~{Dumbser}, O.~{Zanotti}, R.~{Loub{\`e}re}, and S.~{Diot}.
\newblock {A posteriori subcell limiting of the discontinuous Galerkin finite
  element method for hyperbolic conservation laws}.
\newblock {\em Journal of Computational Physics}, 278:47--75, 2014.

\bibitem{Erturk2008}
E.~Erturk.
\newblock Numerical solutions of 2-d steady incompressible flow over a
  backward-facing step, part i: High reynolds number solutions.
\newblock {\em Computers and Fluids}, 37(6):633 -- 655, 2008.

\bibitem{Blood3D2014}
F.~Fambri, M.~Dumbser, and V.~Casulli.
\newblock An efficient semi-implicit method for three-dimensional
  non-hydrostatic flows in compliant arterial vessels.
\newblock {\em International Journal for Numerical Methods in Biomedical
  Engineering}, 30(11):1170--1198, 2014.

\bibitem{Fortin1981}
M.~Fortin.
\newblock Old and new finite elements for incompressible flows.
\newblock {\em International Journal for Numerical Methods in Fluids},
  1(4):347--364, 1981.

\bibitem{Gassner2011}
G.~Gassner and D.A. Kopriva.
\newblock A comparison of the dispersion and dissipation errors of {Gauss} and
  {Gauss-Lobatto} discontinuous {Galerkin} spectral element methods.
\newblock {\em SIAM Journal on Scientific Computing}, 33(5):2560--2579, 2011.

\bibitem{stedg2}
G.~Gassner, F.~L\"orcher, and C.~D. Munz.
\newblock A discontinuous {Galerkin} scheme based on a space-time expansion
  {II.} viscous flow equations in multi dimensions.
\newblock {\em Journal of Scientific Computing}, 34:260--286, 2008.

\bibitem{MunzDiffusionFlux}
G.~Gassner, F.~L\"orcher, and C.D. Munz.
\newblock A contribution to the construction of diffusion fluxes for finite
  volume and discontinuous {Galerkin} schemes.
\newblock {\em Journal of Computational Physics}, 224:1049--1063, 2007.

\bibitem{Gassner2013}
G.J. Gassner.
\newblock A skew-symmetric discontinuous {Galerkin} spectral element
  discretization and its relation to sbp-sat finite difference methods.
\newblock {\em SIAM Journal on Scientific Computing}, 35(3):A1233--A1253, 2013.

\bibitem{Gassner2016}
G.J. Gassner, A.R. Winters, and D.A. Kopriva.
\newblock A well balanced and entropy conservative discontinuous {Galerkin}
  spectral element method for the shallow water equations.
\newblock {\em Applied Mathematics and Computation}, 272:291--308, 2016.

\bibitem{Ghia1982}
U.~Ghia, K.N. Ghia, and C.T. Shin.
\newblock High-{Re} solutions for incompressible flow using the {Navier-Stokes}
  equations and a multigrid method.
\newblock {\em Journal of Computational Physics}, 48(3):387 -- 411, 1982.

\bibitem{GiraldoRestelli}
F.~X. Giraldo and M.~Restelli.
\newblock High-order semi-implicit time-integrators for a triangular
  discontinuous {Galerkin} oceanic shallow water model.
\newblock {\em International Journal for Numerical Methods in Fluids},
  63(9):1077--1102, 2010.

\bibitem{GSz}
U.~Grenander and G.~Szeg{\"o}.
\newblock {\em Toeplitz {F}orms and {T}heir {A}pplications}, volume 321.
\newblock Second Edition, Chelsea, New York, 1984.

\bibitem{HarlowWelch}
F.~H. Harlow and J.~E. Welch.
\newblock Numerical calculation of time‐dependent viscous incompressible flow
  of fluid with free surface.
\newblock {\em Physics of Fluids}, 8(12):2182--2189, 1965.

\bibitem{HartmannHouston1}
R.~Hartmann and P.~Houston.
\newblock Symmetric interior penalty {DG} methods for the compressible
  {Navier--Stokes} equations {I}: Method formulation.
\newblock {\em Int. J. Num. Anal. Model.}, 3:1--20, 2006.

\bibitem{HartmannHouston2}
R.~Hartmann and P.~Houston.
\newblock {An optimal order interior penalty discontinuous Galerkin
  discretization of the compressible Navier--Stokes equations}.
\newblock {\em Journal of Computational Physics}, 227:9670--9685, 2008.

\bibitem{Heywood1982}
J.G. Heywood and R.~Rannacher.
\newblock Finite element approximation of the nonstationary {Navier-Stokes}
  problem. {I}. regularity of solutions and second-order error estimates for
  spatial discretization.
\newblock {\em SIAM Journal on Numerical Analysis}, 19(2):275--311, 1982.

\bibitem{Heywood1988}
J.G. Heywood and R.~Rannacher.
\newblock Finite element approximation of the nonstationary {Navier-Stokes}
  problem {III}. smoothing property and higher order error estimates for
  spatial discretization.
\newblock {\em SIAM Journal on Numerical Analysis}, 25(3):489--512, 1988.

\bibitem{Hidalgo2011}
A.~Hidalgo and M.~Dumbser.
\newblock {ADER} schemes for nonlinear systems of stiff advection diffusion
  reaction equations.
\newblock {\em Journal of Scientific Computing}, 48:173--189, 2011.

\bibitem{HouLiu}
S.~Hou and X.~D. Liu.
\newblock Solutions of multi-dimensional hyperbolic systems of conservation
  laws by square entropy condition satisfying discontinuous {Galerkin} method.
\newblock {\em Journal of Scientific Computing}, 31:127--151, 2007.

\bibitem{Hughes1986}
T.~J.R. Hughes, M.~Mallet, and M.~Akira.
\newblock A new finite element formulation for computational fluid dynamics:
  {II. Beyond SUPG}.
\newblock {\em Computer Methods in Applied Mechanics and Engineering},
  54(3):341 -- 355, 1986.

\bibitem{Jiang1994}
G.~S. {Jiang} and C.-W. {Shu}.
\newblock {On a cell entropy inequality for discontinuous Galerkin methods}.
\newblock {\em Mathematics of Computation}, 62:531--538, 1994.

\bibitem{Karniadakis91}
G.~E. Karniadakis, M.~Israeli, and S.~A. Orszag.
\newblock High-order splitting methods for the incompressible {Navier-Stokes}
  equations.
\newblock {\em Journal of Computational Physics}, 97(2):414 -- 443, 1991.

\bibitem{KimMoin85}
J.~Kim and P.~Moin.
\newblock Application of a fractional-step method to incompressible
  {Navier--Stokes} equations.
\newblock {\em Journal of Computational Physics}, 59(2):308 -- 323, 1985.

\bibitem{KlaijVanDerVegt}
C.~Klaij, J.J.W.~Van der Vegt, and H.~Van der Ven.
\newblock {Space-time discontinuous Galerkin method for the compressible
  Navier-Stokes equations}.
\newblock {\em Journal of Computational Physics}, 217:589--611, 2006.

\bibitem{KleinKummerOberlack2013}
B.~Klein, F.~Kummer, and M.~Oberlack.
\newblock {A SIMPLE based discontinuous Galerkin solver for steady
  incompressible flows}.
\newblock {\em Journal of Computational Physics}, 237:235--250, 2013.

\bibitem{KleinMach}
R.~Klein.
\newblock Semi-implicit extension of a godunov-type scheme based on low mach
  number asymptotics {I}: one-dimensional flow.
\newblock {\em Journal of Computational Physics}, 121:213--237, 1995.

\bibitem{klein}
R.~Klein, N.~Botta, T.~Schneider, {C.D.} Munz, S.Roller, A.~Meister,
  L.~Hoffmann, and T.~Sonar.
\newblock Asymptotic adaptive methods for multi-scale problems in fluid
  mechanics.
\newblock {\em J. of Eng. Math.}, 39:261--343, 2001.

\bibitem{Kopriva2006}
D.A. Kopriva.
\newblock Metric identities and the discontinuous spectral element method on
  curvilinear meshes.
\newblock {\em Journal of Scientific Computing}, 26(3):301--327, 2006.

\bibitem{Kopriva2010}
D.A. Kopriva and G.~Gassner.
\newblock On the quadrature and weak form choices in collocation type
  discontinuous {Galerkin} spectral element methods.
\newblock {\em Journal of Scientific Computing}, 44(2):136--155, 2010.

\bibitem{Ku1987}
H.~C. Ku, R.~S. Hirsh, and T.~D. Taylor.
\newblock A pseudospectral method for solution of the three-dimensional
  incompressible {Navier-Stokes} equations.
\newblock {\em Journal of Computational Physics}, 70(2):439 -- 462, 1987.

\bibitem{Kurzweg85}
U.~H. Kurzweg.
\newblock Enhanced heat conduction in oscillating viscous flows within
  parallel-plate channels.
\newblock {\em Journal of Fluid Mechanics}, 156:291--300, 7 1985.

\bibitem{Landau-Lifshitz6}
L.~D. Landau and E.~M. Lifshitz.
\newblock {\em Fluid Mechanics, Course of Theoretical Physics, Volume 6}.
\newblock Elsevier Butterworth-Heinemann, Oxford, 2004.

\bibitem{Layton2008}
A.T. Layton.
\newblock On the choice of correctors for semi-implicit {Picard} deferred
  correction methods.
\newblock {\em Applied Numerical Mathematics}, 58(6):845--858, 2008.

\bibitem{LeeMateescu1998}
T.~Lee and D.~Mateescu.
\newblock Experimental and numerical investigation of 2-d backward-facing step
  flow.
\newblock {\em Journal of Fluids and Structures}, 12(6):703 -- 716, 1998.

\bibitem{levy2004}
D.~Levy, C.W. Shu, and J.~Yan.
\newblock Local discontinuous {Galerkin} methods for nonlinear dispersive
  equations.
\newblock {\em Journal of Computational Physics}, 196(2):751--772, 2004.

\bibitem{Liu2008}
C.~Liu, C.W. Shu, E.~Tadmor, and M.~Zhang.
\newblock L2 stability analysis of the central discontinuous galerkin method
  and a comparison between the central and regular discontinuous galerkin
  methods.
\newblock {\em ESAIM: Mathematical Modelling and Numerical Analysis},
  42(04):593--607, 2008.

\bibitem{Liu2007}
Y.~Liu, C.W. Shu, E.~Tadmor, and M.~Zhang.
\newblock {Central discontinuous Galerkin methods on overlapping cells with a
  nonoscillatory hierarchical reconstruction}.
\newblock {\em SIAM Journal on Numerical Analysis}, 45(6):2442--2467, 2007.

\bibitem{Loudon1998}
C.~Loudon and A.~Tordesillas.
\newblock The use of the dimensionless womersley number to characterize the
  unsteady nature of internal flow.
\newblock {\em Journal of Theoretical Biology}, 191(1):63 -- 78, 1998.

\bibitem{Marcus84}
Philip~S. Marcus.
\newblock Simulation of taylor-couette flow. part 1. numerical methods and
  comparison with experiment.
\newblock {\em Journal of Fluid Mechanics}, 146:45--64, 9 1984.

\bibitem{MeisterMach}
A.~Meister.
\newblock Asymptotic single and multiple scale expansions in the low mach
  number limit.
\newblock {\em SIAM Journal on Applied Mathematics}, 60(1):256--271, 1999.

\bibitem{minion2003}
M.~L. Minion.
\newblock Semi-implicit spectral deferred correction methods for ordinary
  differential equations.
\newblock {\em Commun. Math. Sci.}, 1(3):471--500, 2003.

\bibitem{minion2001}
M.~L. Minion.
\newblock Higher-order semi-implicit projection methods.
\newblock {\em in M. Hafez, editor}, Numerical Simulations of Incompressible
  Flows: Proceedings of a Conference Held at Half Moon Bay, CA (June 18-20,
  2001), January 2003.

\bibitem{Morf1980}
R.~H. Morf, S.~A. Orszag, and U.~Frisch.
\newblock Spontaneous singularity in three-dimensional inviscid, incompressible
  flow.
\newblock {\em Phys. Rev. Lett.}, 44:572--575, Mar 1980.

\bibitem{Mouza2005}
A.A. Mouza, M.N. Pantzali, S.V. Paras, and J.~Tihon.
\newblock Experimental and numerical study of backward-facing step flow.
\newblock {\em 5th National Chemical Engineering Conference, Thessaloniki,
  Greece}, 2005.

\bibitem{munz}
{C.D.} Munz, R.~Klein, S.~Roller, and {K.J.} Geratz.
\newblock The extension of incompressible flow solvers to the weakly
  compressible regime.
\newblock {\em Computers and Fluids}, pages 173--196, 2003.

\bibitem{Patankar1972}
S.~V. Patankar and D.~B. Spalding.
\newblock A calculation procedure for heat, mass and momentum transfer in
  three-dimensional parabolic flows.
\newblock {\em International Journal of Heat and Mass Transfer}, 15(10):1787 --
  1806, 1972.

\bibitem{patankar}
V.S. Patankar.
\newblock {\em Numerical {Heat} {Transfer} and {Fluid} {Flow}}.
\newblock series in computational methods in mechanics and thermal sciences.
  Hemisphere Publishing Corporation, 1980.

\bibitem{PeshRom2014}
I.~Peshkov and E.~Romenski.
\newblock A hyperbolic model for viscous {Newtonian} flows.
\newblock {\em Continuum Mechanics and Thermodynamics}, 28:85--104, 2016.

\bibitem{Prandtl1904}
L.~Prandtl.
\newblock {\"Uber Fl\"ussigkeitsbewegung bei sehr kleiner Reibung}.
\newblock {\em Verhandlg. {III.} Intern. Math. Kongr. Heidelberg}, pages
  484--491, 1904.

\bibitem{QiuDumbserShu}
J.~Qiu, M.~Dumbser, and {C.W.} Shu.
\newblock The discontinuous {Galerkin} method with {Lax}-{Wendroff} type time
  discretizations.
\newblock {\em Computer Methods in Applied Mechanics and Engineering},
  194:4528--4543, 2005.

\bibitem{Rani2007}
H.~P. Rani, Tony W.~H. Sheu, and Eric S.~F. Tsai.
\newblock Eddy structures in a transitional backward-facing step flow.
\newblock {\em Journal of Fluid Mechanics}, 588:43--58, 10 2007.

\bibitem{Reed:1973}
W.~H Reed and T.~R. Hill.
\newblock Triangular mesh methods for the neutron transport equation.
\newblock Technical report, Los Alamos Scientific Laboratory, 1973.

\bibitem{Rusanov1961a}
V.~V. Rusanov.
\newblock {Calculation of Interaction of Non--Steady Shock Waves with
  Obstacles}.
\newblock {\em J. Comput. Math. Phys. USSR}, 1:267--279, 1961.

\bibitem{saad:gmres}
Y.~Saad and M.~H. Schultz.
\newblock {GMRES}: {A} generalized minimum residual algorithm for solving
  nonsymmetric linear systems.
\newblock {\em SIAM J. Sci. Stat. Comput.}, 7(3):856--869, 1986.

\bibitem{BLTheory}
H.~Schlichting and K.~Gersten.
\newblock {\em Grenzschichttheorie}.
\newblock Springer Verlag, 2005.

\bibitem{serra1998}
S.~Serra-Capizzano.
\newblock Asymptotic results on the spectra of block {T}oeplitz preconditioned
  matrices.
\newblock {\em SIAM journal on matrix analysis and applications}, 20(1):31--44,
  1998.

\bibitem{glt}
S~Serra-Capizzano.
\newblock Generalized locally {T}oeplitz sequences: spectral analysis and
  applications to discretized partial differential equations.
\newblock {\em Linear Algebra Appl.}, 366:371--402, 2003.

\bibitem{Shu88}
C.~W. Shu and S.~J. Osher.
\newblock Efficient implementation of essentially non-oscillatory
  shock-capturing schemes.
\newblock {\em J. Comput. Phys.}, 77:439, 1988.

\bibitem{taube_jsc}
A.~Taube, M.~Dumbser, D.~Balsara, and {C.D.} Munz.
\newblock Arbitrary high order discontinuous {Galerkin} schemes for the
  magnetohydrodynamic equations.
\newblock {\em Journal of Scientific Computing}, 30:441--464, 2007.

\bibitem{TavelliDumbser2014}
M.~Tavelli and M.~Dumbser.
\newblock A high order semi-implicit discontinuous {Galerkin} method for the
  two dimensional shallow water equations on staggered unstructured meshes.
\newblock {\em Applied Mathematics and Computation}, 234:623 -- 644, 2014.

\bibitem{TavelliDumbser2014b}
M.~Tavelli and M.~Dumbser.
\newblock A staggered semi-implicit discontinuous {Galerkin} method for the two
  dimensional incompressible {Navier-Stokes} equations.
\newblock {\em Applied Mathematics and Computation}, 248:70 -- 92, 2014.

\bibitem{TavelliDumbser2015}
M.~Tavelli and M.~Dumbser.
\newblock A staggered space-time discontinuous {Galerkin} method for the
  incompressible {Navier-Stokes} equations on two-dimensional triangular
  meshes.
\newblock {\em Computers \& Fluids}, 119:235 -- 249, 2015.

\bibitem{TavelliDumbser2016}
M.~Tavelli and M.~Dumbser.
\newblock A staggered space-time discontinuous {Galerkin} method for the
  three-dimensional incompressible {Navier-Stokes} equations on unstructured
  tetrahedral meshes.
\newblock {\em Journal of Computational Physics}, 319:294 -- 323, 2016.

\bibitem{Taylor1973}
C.~Taylor and P.~Hood.
\newblock A numerical solution of the {Navier--Stokes} equations using the
  finite element technique.
\newblock {\em Computers \& Fluids}, 1(1):73 -- 100, 1973.

\bibitem{Toro99}
E.~F. Toro.
\newblock {\em Riemann Solvers and Numerical Methods for Fluid Dynamics}.
\newblock Springer-Verlag, 1999.

\bibitem{USFORCE}
E.F. Toro, A.~Hidalgo, and M.~Dumbser.
\newblock {FORCE} schemes on unstructured meshes {I}: Conservative hyperbolic
  systems.
\newblock {\em Journal of Computational Physics}, 228:3368--3389, 2009.

\bibitem{Tumolo2013}
G.~Tumolo, L.~Bonaventura, and M.~Restelli.
\newblock A semi-implicit, semi-lagrangian, p-adaptive discontinuous galerkin
  method for the shallow water equations.
\newblock {\em Journal of Computational Physics}, 232(1):46 -- 67, 2013.

\bibitem{Tylli2002}
N.~Tylli, L.~Kaiktsis, and B.~Ineichen.
\newblock Sidewall effects in flow over a backward-facing step: Experiments and
  numerical simulations.
\newblock {\em Physics of Fluids}, 14(11):3835--3845, 2002.

\bibitem{TyZ}
E.~Tyrtyshnikov and N.~Zamarashkin.
\newblock Spectra of multilevel toeplitz matrices: advanced theory via simple
  matrix relationships.
\newblock {\em Linear Algebra Appl.}, 270:15--27, 1998.

\bibitem{spacetimedg1}
J.~J.~W. van~der Vegt and H.~van~der Ven.
\newblock Space--time discontinuous {Galerkin} finite element method with
  dynamic grid motion for inviscid compressible flows {I}. general formulation.
\newblock {\em Journal of Computational Physics}, 182:546--585, 2002.

\bibitem{spacetimedg2}
H.~van~der Ven and J.~J.~W. van~der Vegt.
\newblock Space--time discontinuous {Galerkin} finite element method with
  dynamic grid motion for inviscid compressible flows {II}. efficient flux
  quadrature.
\newblock {\em Comput. Methods Appl. Mech. Engrg.}, 191:4747--4780, 2002.

\bibitem{vanKan1986}
J~van Kan.
\newblock A second-order accurate pressure correction scheme for viscous
  incompressible flow.
\newblock {\em SIAM J. Sci. Stat. Comput.}, 7(3):870--891, July 1986.

\bibitem{Verfuerth}
R.~Verf\"urth.
\newblock {Finite element approximation of incompressible Navier-Stokes
  equations with slip boundary condition II}.
\newblock {\em Numerische Mathematik}, 59:615--636, 1991.

\bibitem{Womersley}
J.~R. Womersley.
\newblock Method for the calculation of velocity, rate of flow and viscous drag
  in arteries when the pressure gradient is known.
\newblock {\em The Journal of Physiology}, 127(3):553--563, 1955.

\bibitem{yan2002}
J.~Yan and C.W. Shu.
\newblock A local discontinuous {Galerkin} method for {KdV} type equations.
\newblock {\em SIAM Journal on Numerical Analysis}, 40(2):769--791, 2002.

\bibitem{Zang86}
T.~A. Zang and M.~Y. Hussaini.
\newblock On spectral multigrid methods for the time-dependent {Navier--Stokes}
  equations.
\newblock {\em Applied Mathematics and Computation}, 19(1-4):359 -- 372, 1986.

\bibitem{Zanotti2015}
O.~Zanotti and M.~Dumbser.
\newblock A high order special relativistic hydrodynamic and
  magnetohydrodynamic code with space-time adaptive mesh refinement.
\newblock {\em Computer Physics Communications}, 188:110--127, 2015.

\bibitem{Zanotti2015d}
O.~{Zanotti}, F.~{Fambri}, and M.~{Dumbser}.
\newblock {Solving the relativistic magnetohydrodynamics equations with ADER
  discontinuous Galerkin methods, a posteriori subcell limiting and adaptive
  mesh refinement}.
\newblock {\em Mon. Not. R. Astron. Soc.}, 452:3010--3029, September 2015.

\bibitem{Zanotti2015c}
O.~{Zanotti}, F.~{Fambri}, M.~{Dumbser}, and A.~{Hidalgo}.
\newblock Space-time adaptive {ADER} discontinuous {Galerkin} finite element
  schemes with a posteriori sub-cell finite volume limiting.
\newblock {\em Computers and Fluids}, 118:204 -- 224, 2015.

\end{thebibliography}
%=============================================================================

%=============================================================================
%==========  A P P E N D I X
%\input{Appendix}
%=============================================================================

\clearpage

\listoffigures

\listoftables


\clearpage

\begin{table} %[!htbp] 
 \centering
 %\numerikNine
\scriptsize
 \begin{tabular}{c r cc cc c c    r cc cc c   c }
   \multicolumn{15}{c}{\textbf{2D Womersley problem --- Spectral-DG-$\mathbb{P}_N$}} \\
   \hline
	&$N_{\text{el}}$ &  $\epsilon_{L^2}$ & $\epsilon_{L^\infty}$ &  $\mathcal{O}_{L^2}$   & $\mathcal{O}_{L^\infty}$  &    & &$N_{\text{el}}$ &  $\epsilon_{L^2}$ & $\epsilon_{L^\infty}$ &  $\mathcal{O}_{L^2}$   & $\mathcal{O}_{L^\infty}$  &     &  \\
	%%%%%%%%%%%%%%%%%%%%%%%%%%%%%
		\hline
		\\
		&\multicolumn{4}{l}{$N=M=1$} & \multicolumn{2}{c}{ ($t_{\text{end}}$ = 0.5)} & & \multicolumn{7}{l}{$N=M=2$}\\
   \cline{2-7}    \cline{9-14}
	& 40$^2$	&    1.6288E-03  &  3.1562E-03    &  ---	   & ---		    &  		   	&  	&20$^2$	&  1.1014E-04  & 4.7994E-04  &  ---	  & ---	  & 	 	 \\
	& 50$^2$	&     1.3100E-03 &   2.5073E-03    &  0.98	  & 1.03    &  			  										&  &25$^2$	&  5.5799E-05  & 2.5608E-04  & 3.05  & 2.82  & 				             \\
	& 60$^2$	&    1.0947E-03  &  2.0868E-03     &  0.98	  & 1.00    &  														&  &30$^2$	&  3.1934E-05  & 1.5173E-04  & 3.06  & 2.87  & 					 	  \\
	& 70$^2$	&     9.3992E-04  &  1.7923E-03    &  0.99	 & 0.99	    &    														&  &35$^2$	&  1.9913E-05  & 9.6986E-05  & 3.06  & 2.90  &						     \\
	\\
		&\multicolumn{6}{l}{$N=M=3$} & & \multicolumn{7}{l}{$N=M=4$}\\
   \cline{2-7}    \cline{9-14}
	&16$^2$  &  2.0842E-05  & 7.7657E-05  &  ---	& ---	  		& 		 		&  & 15$^2$	  & 2.9044E-06  & 1.6983E-05  &  ---   & ---   &		 	 \\
	&20$^2$  &   8.9568E-06 &  3.1298E-05  &  3.78  &   4.07  & 													&  & 20$^2$	  & 7.3904E-07  & 4.9121E-06  &  4.76  &  4.31 &	  			\\
	&24$^2$  &  4.4439E-06  & 1.6255E-05  &  3.84   &   3.60  &  					   							&  & 25$^2$	  & 2.3730E-07  & 1.7422E-06  &   5.09 &   4.65 &	  	  	 	\\
	&28$^2$  &  2.4474E-06  & 9.0529E-06  &  3.87  	&   3.80  &   					 							&  & 30$^2$	  & 9.2319E-08  & 7.2933E-07  &  5.18  &   4.78 & 		 	 \\
	\\
		&\multicolumn{6}{l}{$N=M=5$} & & \multicolumn{4}{l}{$N=M=6$}&  \multicolumn{2}{l}{ ($t_{\text{end}}$ = 2.2)}  \\
   \cline{2-7}    \cline{9-14}
	& 8$^2$	    & 4.1408E-06  & 2.1569E-05   & ---	 &  ---	  &  	  				&  & 9$^2$	  & 1.7483E-07   &  1.0630E-06  &  ---	&  ---    &  			\\
	& 12$^2$	  & 7.1268E-07  &  4.5893E-06  &  4.34 &  3.82  &  														&  & 10$^2$	  & 8.4841E-08  &   5.1637E-07  &   6.86  &   6.85   &		\\
	& 16$^2$	  & 1.7093E-07  & 1.2088E-06   &  4.96 & 4.64   &  												  	&  & 11$^2$	  & 4.3911E-08   &  2.6415E-07   &   6.91 &   7.03   &		\\
	& 20$^2$	  & 5.0328E-08  & 3.6841E-07   &  5.48 &  5.32  &  														&  & 12$^2$	  & 2.3925E-08  &   1.4122E-07   &   6.98  &  7.20  &		\\	
	\\ 
%%%%%%%%%%%%%%%%%%%%%%%%%%%%%%%
   \hline
 \end{tabular}
\caption{Numerical convergence table for the two dimensional oscillatory flow between two flat plates computed with staggered spectral space-time DG schemes for $N=M=1,\ldots,6$.}
\label{tab:Womersley}
 \end{table}

		\begin{table} %[!htbp] 
 \centering
 %\numerikNine
\scriptsize
 \begin{tabular}{c r cc cc c c    r cc cc c   c }
   \multicolumn{15}{c}{\textbf{2D Taylor-Green vortex problem $L=2\pi$ --- Staggered spectral space-time DG-$\mathbb{P}_N$}} \\
   \hline
	&$N_{\text{el}}$ &  $\epsilon_{L^2}$ & $\epsilon_{L^\infty}$ &  $\mathcal{O}_{L^2}$   & $\mathcal{O}_{L^\infty}$  &     & &$N_{\text{el}}$ &  $\epsilon_{L^2}$ & $\epsilon_{L^\infty}$ &  $\mathcal{O}_{L^2}$   & $\mathcal{O}_{L^\infty}$  &     &  \\
	%%%%%%%%%%%%%%%%%%%%%%%%%%%%%
		\hline
		\\
		&\multicolumn{6}{l}{$N=M=1$} & & \multicolumn{7}{l}{$N=M=2$}\\
   \cline{2-7}    \cline{9-14}
	& 20$^2$	&  6.8094E-02 &  2.0375E-02  &  ---	 & ---	 &   	&   &10$^2$	&     2.9023E-03  & 1.3592E-03  & ---	  & ---	 &	 					\\
	& 25$^2$	&  5.2703E-02 &  1.5419E-02  &  1.15 &   1.25 &  										&   &15$^2$	&     6.6991E-04 &  3.2999E-04  &  3.62 &  3.49 & 							\\  
	& 30$^2$	&  4.2541E-02 &  1.2114E-02  &  1.17 &   1.32 &  										&  & 20$^2$	&     2.1652E-04 &  9.4734E-05  &   3.93 &  4.34& 									\\
	& 35$^2$	&  3.5663E-02 &  1.0062E-02  &  1.14 &   1.20 &  										&  &25$^2$	&     9.2092E-05 &  5.0350E-05  &  3.83 &  2.83&  							\\
	\\
		&\multicolumn{6}{l}{$N=M=3$} & & \multicolumn{7}{l}{$N=M=4$}\\
   \cline{2-7}    \cline{9-14}
& 	16$^2$	&  1.3626E-04 &  7.8999E-05   & ---	  & ---  &  	 		&  &16$^2$	&   1.0519E-06  & 6.1297E-07 & ---	  & ---  &      &				\\
&	20$^2$	&  6.5874E-05  & 3.7886E-05 & 3.26   &   3.29   &  							&  &20$^2$	&    2.7970E-07 &  1.6271E-07 &   5.94 &   5.94  &                   &				\\
&	24$^2$	&   3.6213E-05 &  2.0085E-05 &  3.28   &   3.48  & 	  						&  &24$^2$	&   9.1334E-08  & 5.3150E-08 &   6.14 &   6.14  &                  	&				\\
& 	28$^2$	& 2.1887E-05 &  1.1756E-05  &  3.27  &    3.47  &  							&  &28$^2$	&   3.4806E-08  & 2.0411E-08 &  6.26 &    6.21  &                	&				\\
		\\
		&\multicolumn{6}{l}{$N=M=5$} & & \multicolumn{7}{l}{$N=M=6$}\\
   \cline{2-7}    \cline{9-14}
&  12$^2$	&     2.9725E-07  & 2.1057E-07 &  ---	  & --- &  	  	&&12$^2$	&   3.9994E-09 &  2.4762E-09 & ---	  & --- &   					\\    
&  15$^2$&    8.4575E-08  & 6.3046E-08 &  5.63 &   5.40 &  		 						&&15$^2$	& 9.2671E-10  & 5.6882E-10 &  6.55 &  6.59 &                  						\\
&  18$^2$	&     2.9174E-08  & 2.2565E-08 &  5.84 &   5.64 &  	 						&&18$^2$	& 2.6783E-10  & 1.7230E-10 &  6.81 &  6.55 &                 							\\
&  21$^2$	&     1.1910E-08  & 9.4502E-09  &  5.81 &    5.65 & 						&&21$^2$	&  8.9312E-11 &  6.0292E-11 &  7.12 & 6.81 &                     						\\
		\\
			&\multicolumn{6}{l}{$N=M=7$} & & \multicolumn{7}{l}{$N=M=8$}\\
   \cline{2-7}    \cline{9-14}
&6$^2$	& 1.0586E-08 &  6.3735E-09  &  ---	  & --- &       	&&4$^2$	& 3.4616E-08 &   2.4849E-08  &   ---	  & ---	  &	 		  				\\ 
&	9$^2$	& 4.8791E-10 &  2.6452E-10 &   7.59 &   7.85 &						&&6$^2$	&  1.5605E-09 &   7.1039E-10  &    7.64 & 8.77  &	 	 								\\
&12$^2$	& 7.2738E-11 &  3.8053E-11 &   6.62 &   6.74 &       					&&8$^2$	&  6.1523E-11 &   3.3955E-11  &   11.24 &  10.57  & 	 								\\
&15$^2$	&  1.2830E-11 &  7.5665E-12 &   7.78 &    7.24 &  						&&10$^2$	&   5.7787E-12 &   4.9803E-12  &   10.60  & 8.60 &	 								\\
\\
%%%%%%%%%%%%%%%%%%%%%%%%%%%%%%%
   \hline
 \end{tabular}
\caption{Numerical convergence table computed for the two dimensional Taylor-Green vortex problem using staggered spectral space-time DG schemes with $N=M=1,\ldots,8$.}
\label{tab:TGV2D}
 \end{table}

		\begin{table} %[!htbp] 
 \centering
 %\numerikNine
\scriptsize
 \begin{tabular}{c r cc cc c c    r cc cc c   c }
   \multicolumn{15}{c}{\textbf{3D ABC flow problem  --- Staggered spectral space-time DG-$\mathbb{P}_N$}} \\
   \hline
	&$N_{\text{el}}$ &  $\epsilon_{L^2}$ & $\epsilon_{L^\infty}$ &  $\mathcal{O}_{L^2}$   & $\mathcal{O}_{L^\infty}$  &     & &$N_{\text{el}}$ &  $\epsilon_{L^2}$ & $\epsilon_{L^\infty}$ &  $\mathcal{O}_{L^2}$   & $\mathcal{O}_{L^\infty}$  &    &  \\
	%%%%%%%%%%%%%%%%%%%%%%%%%%%%%
		\hline
	\\
	&\multicolumn{6}{l}{$N=M=1$} & & \multicolumn{7}{l}{$N=M=2$}\\
   \cline{2-7}    \cline{9-14}
& 12$^2$	&   8.6905E-02 &   1.0949E-02  &    ---	  & ---	    &  				&&12$^2$	&    6.9061E-03 &  1.7183E-03  &  ---	  & ---  &  	 				\\
& 16$^2$	&   4.8615E-02 &   6.6433E-03   &    2.02  &   1.74  &   								&&15$^2$	&     3.3559E-03 &  8.2081E-04  & 3.23 &  3.31  &  								\\
& 	20$^2$	&   3.2626E-02  &  4.7890E-03  &   1.79 &   1.47 &  									&&18$^2$	&     1.7856E-03 &  4.2641E-04 &   3.46 &   3.59   & 							\\
& 	24$^2$	&   2.3886E-02  &  3.6992E-03  &   1.71 &   1.42 &  									&&21$^2$	&     1.0316E-03 &  2.4384E-04 &  3.56 &  3.63 &  								\\
	\\
	&\multicolumn{6}{l}{$N=M=3$} & & \multicolumn{7}{l}{$N=M=4$}\\
   \cline{2-7}    \cline{9-14}
& 4$^2$	&     1.2102E-02 &  2.3277E-03  &  ---	  & --- &    		   && 4$^2$	&     1.2102E-02 &  2.3277E-03  &  ---	  & --- &    		   							\\
& 6$^2$	&    2.1258E-03  & 4.9684E-04  &  4.29 &  3.81 &     		  					&& 6$^2$	&    2.1258E-03  & 4.9684E-04  &  4.29 &  3.81 &     		  							\\
& 8$^2$	&     6.4822E-04 &  1.6790E-04  &   4.13 &   3.77 &    		  					&& 8$^2$	&     6.4822E-04 &  1.6790E-04  &   4.13 &   3.77 &    		  							\\
& 10$^2$	&    2.6346E-04 &  7.5594E-05  &  4.03 &   3.58 &     		    			&& 10$^2$	&    2.6346E-04 &  7.5594E-05  &  4.03 &   3.58 &     		    				 			\\
	\\
	&\multicolumn{6}{l}{$N=M=5$} & & \multicolumn{7}{l}{$N=M=6$}\\
   \cline{2-7}    \cline{9-14}
&  2$^2$	&  4.6766E-03  & 9.2223E-04  & ---	  & ---  &    	&&2$^2$	&  2.4587E-04 &  8.9231E-05 &   ---	  & ---  &  	 		\\
& 4$^2$	&  6.6120E-05 &  1.9076E-05  & 6.14   & 5.60 &   						&&4$^2$	&  4.2976E-06 &  1.2496E-06 &   5.84 &   6.16 &  							 \\
& 6$^2$	&  5.7711E-06 &  2.1200E-06  &  6.01  & 5.42 &   						&&6$^2$	&  3.6205E-07 &  1.0417E-07 &  6.10 &   6.13 &  							\\
&  8$^2$	&  1.1153E-06  & 5.0069E-07  &  5.71  & 5.02 &   					&&8$^2$	&  5.6088E-08 & 1.6526E-08 &   6.48 &   6.40 &  							 \\
	\\
	&\multicolumn{6}{l}{$N=M=7$} & & \multicolumn{7}{l}{$N=M=8$}\\
   \cline{2-7}    \cline{9-14}		
& 2$^2$	&  5.4083E-05 &  1.4183E-05  & ---	  & ---  &    	  		&&1$^2$	&  1.5955E-03  &   5.0110E-04   &  ---	  & ---  &   								\\
& 3$^2$	&  2.1818E-06 &  7.6077E-07 &  7.92 &   7.22 &   							&&2$^2$	&  2.1017E-06  &   1.0666E-06  &   9.57  & 8.88 &  								 \\
& 4$^2$	&  2.1037E-07  & 6.2486E-08  &   8.13 &  8.69 &    							&&3$^2$	&  9.7717E-08  &   3.5484E-08  &  7.57 &    8.39 &								\\
& 5$^2$	&  3.5196E-08  & 1.3221E-08  &  8.01  &  6.96 &    							&&4$^2$	&  1.0666E-08  &   3.4279E-09  &  7.70 &    8.12&				\\
\\
%%%%%%%%%%%%%%%%%%%%%%%%%%%%%%
   \hline
 \end{tabular}
\caption{Numerical convergence table for the three dimensional Arnold-Beltrami-Childress (ABC) flow problem computed with staggered spectral space-time DG schemes for $N=M=1,\ldots,8$.}
\label{tab:ABC}
 \end{table}

\begin{table} %[!htbp] 
%\color{red} 
 \centering
 %\numerikNine
\scriptsize
\renewcommand{\arraystretch}{1.4} 
 \begin{tabular}{cr  ccc   c }

%   \multicolumn{6}{c}{\textbf{Spectral-DG-$\mathbb{P}_N$ - Stencil-size}} \\
   \hline
& space	& Collocated grid &  Vertex-based staggered grid &  Edge-based staggered grid & \\
& dimensions	&  (A--grid) &   (B--grid)&  (C--grid)& \\
	%%%%%%%%%%%%%%%%%%%%%%%%%%%%%
		\hline
	& \textbf{1D}&  $5$	  &   $\mathbf{3}$  &  $\mathbf{3}$    &  
	\\ 
	& \textbf{2D}&  $9$	  &   $9$   &  $\mathbf{5}$    &  
	\\ 
	& \textbf{3D}&  $13$	&   $27$  &  $\mathbf{7}$    &  
	\\ 
%%%%%%%%%%%%%%%%%%%%%%%%%%%%%%%
   \hline
 \end{tabular}
\caption{Total stencil-size for the resulting pressure systems for semi-implicit DG schemes using different grid types for different numbers of space dimensions. In all cases it is assumed 
that the discrete momentum equation is substituted into the discrete continuity equation, in order to yield one single equation system for the scalar pressure.}  
\label{tab:StencilSize}
 \end{table}

\end{document}